\documentclass[aop]{imsart}
\usepackage[T1]{fontenc}
\usepackage[utf8]{inputenc}
\usepackage{lmodern}
\usepackage{amsmath,amssymb,amsthm,mathtools}
\usepackage{microtype}
\usepackage{needspace}
\usepackage[round,authoryear]{natbib}
\usepackage[colorlinks=true,linkcolor=blue,citecolor=blue,urlcolor=blue]{hyperref}
\numberwithin{equation}{section}
\newtheorem{theorem}{Theorem}[section]
\newtheorem{proposition}[theorem]{Proposition}
\newtheorem{lemma}[theorem]{Lemma}
\newtheorem{corollary}[theorem]{Corollary}
\theoremstyle{definition}
\newtheorem{definition}[theorem]{Definition}
\newtheorem{problem}[theorem]{Open problem}
\theoremstyle{definition}
\newtheorem{remark}[theorem]{Remark}
\newcommand{\R}{\mathbb R}
\newcommand{\E}{\mathbb E}
\newcommand{\Prob}{\mathbb P}
\newcommand{\cH}{\mathcal H}
\newcommand{\bB}{\mathbf B}
\newcommand{\bX}{\mathbf X}
\newcommand{\dd}{\,\mathrm d}
\newcommand{\norm}[1]{\left\lVert#1\right\rVert}

\newcommand{\ind}{\mathbf1}

\DeclareMathOperator{\sgn}{sgn}

\DeclareMathOperator{\supp}{supp}

\makeatletter
\def\journal@name{}
\def\nocontentsline{%
  \let\@@addcontentsline\addcontentsline
  \def\addcontentsline##1##2##3{\let\addcontentsline\@@addcontentsline}}
\makeatother
\setpkgattr{author}{prefix}{}
\setpkgattr{copyright}{text}{}
\begin{document}
\begin{frontmatter}
\title{Stationary causal enhancements of fractional Brownian motion
at small Hurst indices}
\runtitle{Stationary causal fractional Brownian enhancements}
\runauthor{X. Xu}
\begin{aug}
\author[A]{\fnms{Xingcheng}~\snm{Xu}\ead[label=e1]{xingcheng.xu18@gmail.com}}
\address[A]{Shanghai Artificial Intelligence Laboratory\printead[presep={\\}]{e1}}
\end{aug}
\begin{abstract}
We construct stationary, complete-past causal geometric enhancements
of prescribed multidimensional fractional Brownian motion for every
\(0<H\leq1/4\). They retain the first level, are exactly self-similar,
and commute with linear coordinate maps.
The construction removes a remote-past cutoff from Volterra signatures.
At reciprocal indices \(H=1/N\), a lagged top-degree logarithmic-velocity
subtraction yields a quantitative joint cutoff and zero-lag limit
with sharp moment bounds.
Centered velocities and white-noise localization provide the critical cancellation.
The limit is stable under weak concentration at zero of deterministic
mass-one signed lag profiles on a common bounded lag range, with
uniformly bounded total variation.
Top-central perturbations give nonuniqueness at nonreciprocal indices,
preserving all lower levels and symmetries.
At reciprocal indices, square-integrable whole-line causal extensions
with fixed lower levels and common-noise shift and scaling identities
differ by a deterministic central drift. Logarithmic centering removes this drift.
The enhancements yield controlled integration, rough differential equations,
a finite-word regularity-structure model and an adapted realization
by ordinary signatures.
\end{abstract}
\end{frontmatter}
\thispagestyle{plain}

\hypersetup{pdfsubject={Stationary causal rough paths and stochastic integration}}
\setcounter{tocdepth}{1}
\tableofcontents

\section{The integration problem and the stationary selection}\label{sec:problem}
Let \(B^H=(B^{H,1},\ldots,B^{H,d})\) have independent centered Gaussian
components with covariance
\begin{equation}\label{eq:covariance}
 \E B^{H,i}_sB^{H,j}_t
 =\frac{\delta_{ij}}2\left(|s|^{2H}+|t|^{2H}-|t-s|^{2H}\right).
\end{equation}
We study geometric integration of smooth one-forms and rough differential
equations with noncommuting vector fields for \(d\geq2\) and \(0<H\leq1/4\).
The input is specified through its complete two-sided increment history.
For \(H>1/4\), smooth Gaussian approximations select a geometric
enhancement. At or below this threshold, off-diagonal second levels of such
approximations can fail to have a subsequential limit in probability
\citep[Section 4.4]{CoutinQian2002}, although selected geometric
enhancements exist for every positive Hurst index
\citep{Unterberger2010,NualartTindel2011,Unterberger2013}.
We seek a stationary enhancement measurable from the complete Gaussian
past and study how its choice affects integrals and differential equations.

We construct such a selection on the prescribed two-sided Gaussian
history and establish quantitative critical cutoff and lag-profile
stability. We use the signature algebra of the adapted finite-origin
Volterra lift of \citet[Theorem~1.1]{NualartTindel2011},
with a two-sided kernel and a remote-past cutoff.
Away from reciprocal indices, weighted signature estimates give the
stationary limit directly. At \(H=1/N\), a lagged degree-\(N\)
logarithmic-velocity subtraction permits joint removal of the remote
cutoff and lag. Centered auxiliary-time velocities and disjoint-window
white-noise localization give convergence below \(H\), with sharp
\(kH\)-moment scaling at every retained level.
For the finite-origin construction,
Propositions~\ref{prop:NT-receding} and \ref{prop:NT-endpoint}
identify logarithmic obstructions both to the uncorrected receding-origin
limit and to a uniform sharp critical moment bound.

\begin{definition}\label{def:requirements}
For a random enhancement \(\bB\), we distinguish the following properties.
\begin{enumerate}
\item \emph{Prescribed first level:}
 \(\pi_1\bB_{s,t}=B^H_t-B^H_s\).
\item \emph{Complete-past causality:}
 \(\bB_{s,t}\) is measurable with respect to
 the completed increment past
 \[
  \mathcal F_t^B=\overline{\sigma\{B_v^H-B_u^H:u<v\leq t\}}.
 \]
\item \emph{Interval-increment locality:}
 \(\bB_{s,t}\) is measurable with respect to
 \(\sigma\{B_r^H-B_s^H:r\in[s,t]\}\).
\item \emph{Stationarity and self-similarity:}
 increment families are invariant in law under time shifts, and
 level \(k\) scales with exponent \(kH\).
\item \emph{Approximation consistency:}
 a specified enhanced approximation converges in a topology supporting
 continuity of rough integration.
\item \emph{Coordinate naturality:}
 the enhancement commutes with fixed linear spatial maps.
\end{enumerate}
\end{definition}

Appendix~\ref{sec:filtration} identifies the complete fBm increment past
with the Wiener past of the moving-average representation.
For forward integration beginning at \(a\in\R\), the input is
recentered as \(B_r^H-B_a^H\), \(r\geq a\).
The complete past includes negative-time information beyond an observed
segment \([0,T]\). On that segment, the finite-origin construction of
\citet{NualartTindel2011} supplies an adapted alternative.
The interval-locality results of \citet{ChevyrevFerrucci2026} and
\citet{Han2026} explain the stronger information
restriction in item 3; Section~\ref{sec:locality} gives the
square-integrable obstruction used here.

\Needspace{7\baselineskip}
\begin{theorem}[Stationary causal calculus]\label{thm:main-scope}
Let \(d\geq2\), \(0<H\leq1/4\),
\(N=\lfloor1/H\rfloor\), and \(1/(N+1)<\alpha<H\).
\begin{enumerate}
\item The prescribed fBm has a complete-past causal, stationary,
 exactly self-similar, coordinate-natural geometric step-\(N\)
 enhancement \(\bB\), with
 \[
  \|\pi_k\bB_{s,t}\|_{L^p}\leq C_{k,p}|t-s|^{kH},
       \qquad 1\leq k\leq N,\quad 1\leq p<\infty.
 \]
\item This enhancement supplies adapted integration of adapted controlled
 jets, the geometric chain rule, and pathwise unique RDE solutions for
 bounded \(\operatorname{Lip}^{\gamma}\) vector fields with
 \(\gamma>1/\alpha\).
 Adapted piecewise linear controls realize it through their ordinary
 signatures. Their ordinary one-form integrals and ODE solutions converge
 in probability uniformly on compact forward time intervals.
\item If \(H=1/N\), the lagged cutoff models converge to a scale-free
 model \(Y^0\). On every compact interval, in the unrooted coordinate
 \(\alpha\)-H\"older distance,
 \[
  \|d_\alpha(Y^{R,\ell},Y^0)\|_{L^p}
       \leq C(R^{-H}+\ell^\varepsilon),
  \qquad 0<\varepsilon<\min\{H/2,1-N\alpha\},
 \]
 for \(0<\ell\leq1\), sufficiently large \(R\), and every finite
 \(p\geq1\). The parameters \(R,\ell\) are deterministic.
 The same joint limit holds for deterministic mass-one signed lag
 profiles on \((0,1]\), with uniformly bounded total variation,
 converging weakly to \(\delta_0\) as measures on \([0,1]\).
 Within a common bounded lag range, fixed nonunit mass prevents
 remote-cutoff tightness.
\end{enumerate}
\end{theorem}
\begin{proof}
Theorem~\ref{thm:noncritical-causal}, Theorem~\ref{zl:thm:zero}
and Corollary~\ref{zl:cor:symmetry} give the enhancements and moment
bounds. At reciprocal indices, use the centered representative of
Corollary~\ref{cor:critical-centering} for \(\bB\).
Theorems~\ref{thm:selected-calculus}, \ref{thm:controlled-integration}
and \ref{thm:adapted-BV} give the calculus and ordinary realization;
Corollary~\ref{cor:critical-ordinary-realization} combines the latter
with the critical cutoff limits.
The joint estimate is Theorem~\ref{zl:thm:joint}.
The profile assertions are Theorem~\ref{lp:thm:stability} and
Proposition~\ref{lp:prop:nonunit}.
\end{proof}

The choice of enhancement is visible in integrals and equations.
At nonreciprocal indices, Theorem~\ref{cg:thm:fullgauge} constructs
top-degree central perturbations preserving every lower level and the
properties in part 1. Their polynomial one-form integrals differ, and
an uncountable subfamily has mutually singular infinite-time laws
(Corollaries~\ref{cg:cor:integral} and \ref{cg:cor:singular}).
Theorem~\ref{thm:adapted-BV} applies to every member.
At reciprocal indices, positive lags give stationary causal models
whose scaling acts on the lag parameter according to
\eqref{eq:lag-scaling}; a fixed positive lag breaks exact
self-similarity. Distinct lags have different laws on every
nontrivial finite horizon and mutually singular laws on the infinite
half-line (Theorem~\ref{thm:lag-laws}).
Proposition~\ref{prop:calibration-calculus} identifies the resulting
changes in controlled integrals and noncommuting equations.
Central translation and its bracket interpretation are established
operations \citep{FrizOberhauser2009,BrunedEtAl2019,BellingeriEtAl2022};
the explicit perturbations use homogeneous Wiener kernels
\citep{BaiTaqqu2014}.

At a reciprocal index, we fix all lower levels and require the
functional shift and scaling identities to hold with the same noise.
Within this class, coordinatewise \(L^2\), whole-line causal extensions
differ only by a deterministic degree-\(N\) central drift. The normalization
\(\E[\pi_N\log\bB_{0,1}]=0\) selects a unique representative
(Theorem~\ref{thm:critical-rigidity} and
Corollary~\ref{cor:critical-centering}).
For square-integrable second levels, complete-past causality forces
whole-history time irreversibility (Theorem~\ref{thm:causal-irreversibility}).
Reversal gives a complete-future selection with the same first level,
stationarity, self-similarity and linear naturality. Its law differs
from that of the causal selection on every nontrivial finite horizon, and the
two half-line laws are mutually singular
(Corollaries~\ref{cor:forward-backward-models} and \ref{cor:directed-laws}).

Fourier constructions provide another useful comparison.
The homogeneous frequency restrictions of \citet{Unterberger2010}
are translation and dilation covariant wherever their kernel limits
are justified. Appendix~\ref{sec:fourier-comparison} gives an exact
angular integrability criterion for a real, simultaneously even class
of homogeneous resonance multipliers, including the hard cutoff and
smooth versions without an open zero cone.
Every square-integrable member at \(H\leq1/4\) has a noncausal
off-diagonal kernel in the prescribed innovation representation;
a full selection retaining it has a different complete-history law.
The comparison includes the two-vertex kernel-BPHZ prescription of
\citet{Unterberger2013}
(Section~\ref{sec:kernel-bphz}).
Section~\ref{sec:causal-Fourier} identifies the causal coordinate by a
complex Hermitian multiplier in a broader Fourier family.
The interacting-measure weak construction of
\citet[Theorem~0.1]{MagnenUnterberger2012}, for planar
area with \(H\in(1/8,1/4)\), retains a Gaussian limiting first-level
marginal under a different probability-limit convention from a causal
functional on the prescribed Wiener space.

Variance attenuation addresses the approximation obstruction by changing
the limiting enhanced driver. \citet{Hairer2025} obtains Lie-bracket
Brownian diffusions from attenuated smooth fBm controls, with a vanishing
first level. The stationary causal construction retains the prescribed
fBm first level.

Other approaches to integration against fBm use different domains and
algebraic conventions: extended divergence at \(H<1/2\)
\citep{CheriditoNualart2005}, higher-order regularization and chain rules
\citep{GradinaruEtAl2005}, covariance-based Riemann--Skorohod conversion
\citep{Liu2025}, and It\^o rough paths at \(H\in(1/3,1/2]\)
\citep{QianXu2025}.
The mixed-lift theorem of
\citet[Definition~2.5 and Theorem~2.6]{BonesiniEtAl2026} starts with a
given adapted rough block and Brownian motion. Its all-Hurst fBm
application uses a one-dimensional block; the multidimensional application
requires a separately supplied self-enhancement.
\citet[Section~6]{BellingeriFrizGerencser2022} reduce singular
finite-origin Volterra models to stationary reference noises in the
context of mixed rough-volatility integration against Brownian motion.

The finite-word model in Appendix~\ref{sec:word-model} expresses the
chosen enhancement in the regularity-structure framework
\citep{Hairer2014,BrunedHairerZambotti2019,ChandraHairer2018arxiv}.
Choose \(N=\lfloor1/H\rfloor\) and \(1/(N+1)<\alpha<H\).
The differentiated word symbols have degrees \(|w|\alpha-1\);
the symbol encoding \((B^i-B_s^i)\partial_tB^j\) has degree
\(2\alpha-1\). This finite truncation suffices for controlled
integration because \((N+1)\alpha>1\). The enhancement estimates
supply its random coordinates.

Section~\ref{sec:framework} introduces the signature algebra and kernel
estimates. Sections~\ref{sec:causal-construction} and
\ref{sec:critical-main} construct the enhancements, and
Section~\ref{sec:critical-rigidity} proves critical extension rigidity.
Sections~\ref{sec:integration-main}--\ref{sec:dependence} develop the
calculus, its adapted ordinary-signature realization and its dependence
on the enhancement. Section~\ref{sec:frontier} formulates the remaining
selection and stability questions. The appendices provide the supporting
estimates, information and Fourier comparisons, and the word model.

\section{Framework and Volterra kernels}\label{sec:framework}
\subsection{Rough enhancements and integration conventions}\label{sec:preliminaries}
Fix a nonreciprocal \(\alpha\in(0,H)\) and put \(N=\lfloor1/\alpha\rfloor\).
Write \(G^N(\R^d)\) for the step-\(N\) free nilpotent group, realized as
group-like elements of the truncated tensor algebra.
The homogeneous degree-\(N\) subspace of its Lie algebra is denoted by
\(\mathfrak L_N(\R^d)\), or \(\mathfrak L_N\) when the dimension is fixed.
A weakly geometric \(\alpha\)-H\"older enhancement of a path \(x\) is
a multiplicative family \(\bX_{s,t}\in G^N(\R^d)\) with
\[
 \pi_1\bX_{s,t}=x_t-x_s,\qquad
 \norm{\pi_k\bX_{s,t}}\leq C|t-s|^{k\alpha}.
\]
Geometric enhancements are limits of bounded-variation signatures,
with a possible arbitrarily small loss of H\"older exponent when
passing from weak geometricity to this closure. One can choose
\(\alpha<H\) away from all reciprocal exponents.
Abstract lifting is provided by \citet{LyonsVictoir2007}.

For a sufficiently smooth one-form \(F=(F_1,\ldots,F_d)\), rough integration
uses the compensated sums
\[
 \int_s^t F(x_r)\dd\bX_r
 =\lim_{|\pi|\to0}\sum_{[u,v]\in\pi}
       \sum_{k=1}^N D^{k-1}F(x_u)\,\pi_k\bX_{u,v}.
\]
The contraction uses the signature convention; no additional factorial
is inserted in this expression. In particular \(H=1/4\) requires at
least four levels. The existence of only a second level is insufficient.

Throughout the construction, \(W^1,\ldots,W^d\) are independent
real white noises on \(\R\). Write
\begin{align*}
 \mathcal F_t^W&=\overline{\sigma\{W^i(f):
          \operatorname{supp}f\subset(-\infty,t],\ 1\leq i\leq d\}},\\
 \mathcal F_t^B&=\overline{\sigma\{B_b^{H,i}-B_a^{H,i}:
                    a<b\leq t,\ 1\leq i\leq d\}}.
\end{align*}
The bars denote completion. For the moving-average realization below,
\(\mathcal F_t^W=\mathcal F_t^B\); a proof is given in
Proposition~\ref{prop:complete-past}. Causality is understood for
increments with respect to these fields. On a forward interval beginning
at \(a\), adapted one-form integration uses \(B_r^H-B_a^H\), \(r\geq a\).

All quantitative convergence rates use the unrooted coordinate distance
\begin{equation}\label{eq:coordinate-distance}
 d_{\alpha;I}(U,V)=\max_{1\leq j\leq N}
       \sup_{s<t\text{ in }I}
       \frac{\|\pi_j(U_{s,t}-V_{s,t})\|}{|t-s|^{j\alpha}}.
\end{equation}
We omit \(I\) when it is fixed. Numerical rates in this distance
need not agree with rates in a rooted homogeneous metric.

\begin{lemma}[The central fiber]\label{lem:central-fibre}
Let \(X,\widetilde X\) be multiplicative \(G^N(\R^d)\)-valued
increment families with the same levels below \(N\). Then
\[
 A_{s,t}=\pi_N(\widetilde X_{s,t}-X_{s,t})
\]
is an additive degree-\(N\) Lie-valued increment and
\(\widetilde X_{s,t}=X_{s,t}\exp(A_{s,t})\).
\end{lemma}
\begin{proof}
The group-like ratio \(X_{s,t}^{-1}\widetilde X_{s,t}\) has only
its top positive degree. Its logarithm is therefore the Lie element
\(A_{s,t}\), and its exponential is \(1+A_{s,t}\).
Centrality and Chen's identity give
\(A_{s,t}=A_{s,u}+A_{u,t}\).
\end{proof}

\Needspace{8\baselineskip}
\subsection{Volterra constructions and infinite-past kernels}\label{sec:volterra}
\subsubsection{The finite-origin construction}
For \(H<1/2\), write
\[
 B_t^{H,i}=\int_0^t K_H(t,u)\dd W_u^i .
\]
\citet[Theorem 1.1]{NualartTindel2011} construct all required levels throughout this range
by a sum over valley-ordered domains.
For distinct indices its second level is
\begin{align*}
 \mathbb B^{ij}_{s,t}
 ={}&\int_{0<u<v<t}
 [K_H(t,u)-K_H(s,u)]K_H(t,v)\dd W_u^i\dd W_v^j\\
 &-\int_{0<v<u<t}
 K_H(s,u)[K_H(t,v)-K_H(s,v)]\dd W_u^i\dd W_v^j .
\end{align*}
The underlying integrals extend to repeated indices with the
Stratonovich convention and its contraction terms.
The lift is adapted to the past. Its use of \(u<s\) explains why
it does not furnish interval-increment locality.

\subsubsection{An unanchored two-sided kernel}
Set \(\beta=H-1/2\) and
\[
 \kappa_t(u)=\mathfrak c_H(t-u)_+^\beta .
\]
Here \(\mathfrak c_H>0\) is chosen so that
\(\|\kappa_1-\kappa_0\|_{L^2(\R)}=1\); this fixes the moving-average normalization.
The kernel \(\kappa_t\) itself is not square-integrable at negative
infinity, whereas \(\kappa_t-\kappa_s\) is.
The associated first level is the two-sided moving-average
representation of \citet{MandelbrotVanNess1968}.
Let \(A_j^n\) consist of tuples with \(u_j\) the least coordinate,
with decreasing coordinates to its left and increasing coordinates
to its right. The formal \(n\)-th level is the sum over \(j\) of
\begin{equation}\label{eq:stationary-valley}
 (-1)^{j-1}\int_{A_j^n}
 \left(\prod_{\ell<j}\kappa_s(u_\ell)\right)
 [\kappa_t(u_j)-\kappa_s(u_j)]
 \left(\prod_{\ell>j}\kappa_t(u_\ell)\right)
 \circ\dd W_{u_1}^{i_1}\cdots\circ\dd W_{u_n}^{i_n}.
\end{equation}
The following sections prove cutoff convergence with all
repeated-index contractions when \(nH<1\), and introduce a
central subtraction at the critical highest degree.

For mutually distinct indices, integrate all variables other than
the minimum \(v=u_j\). The squared kernel norm is
\begin{align}
 &\frac{\mathfrak c_H^{2n}}
 {(j-1)!(n-j)!(2H)^{n-1}}\notag\\
 &\quad\times
 \int_{-\infty}^t
 \bigl[(t-v)_+^\beta-(s-v)_+^\beta\bigr]^2
 (s-v)_+^{2H(j-1)}(t-v)^{2H(n-j)}\dd v .
 \label{eq:distinct-kernel-norm}
\end{align}
The integrand behaves at negative infinity as a positive multiple
of \((t-s)^2|v|^{2nH-3}\).
Thus the raw infinite-past kernel has the strict condition \(nH<1\);
at \(nH=1\) it has a logarithmic infrared divergence.
For words with mutually distinct indices, the \(j\)-domains are disjoint,
so their contributions cannot cancel that divergence in the \(L^2\) kernel.

\begin{remark}
Repeated-index contractions are controlled by the estimates in
Section~\ref{sec:weighted}. A highest-chaos
obstruction can also be isolated for the word \(1^{n-1}2\):
restrict its unique color-\(2\) variable to be the least
integration time. After symmetrization of the color-\(1\)
variables only the last valley contributes there, with squared
norm bounded below by a positive multiple of
\(\int_{-R}^s|\kappa_t(v)-\kappa_s(v)|^2
                     (s-v)^{2H(n-1)}\dd v\).
This diverges logarithmically at \(nH=1\); lower-chaos
Stratonovich contractions cannot cancel it.
The finite-origin theorem uses a different kernel and remains
available on finite time intervals.
\end{remark}

\subsection{Weighted signature estimates}\label{sec:weighted}
All signatures in this section are Stratonovich signatures in an
auxiliary Brownian time variable. The terminal parameter is not that
integration variable.

\begin{lemma}[Mixed weights]\label{lem:mixed-weights}
Let \(f_1,\ldots,f_m\in L^2(I)\) be deterministic, \(q\geq2\), and fix
a Brownian coordinate word of length \(m\). The ordered Stratonovich
integral satisfies
\[
 \left\|\int_{u_1<\cdots<u_m}
       \prod_{\ell=1}^m f_\ell(u_\ell)
       \circ\dd W_{u_1}^{i_1}\cdots\circ\dd W_{u_m}^{i_m}
 \right\|_{L^q}
 \leq C_{m,q}\prod_{\ell=1}^m\norm{f_\ell}_2.
\]
The same estimate holds for the running supremum of the ordered integral.
\end{lemma}
\begin{proof}
Prove the assertion first for smooth weights and then pass to
\(L^2\) limits. The Stratonovich-to-It\^o recursion consists of
\[
 \int I_{m-1}(v)f_m(v)\dd W_v^{i_m}
 +\frac{\ind_{\{i_{m-1}=i_m\}}}{2}
       \int I_{m-2}(v)f_{m-1}(v)f_m(v)\dd v.
\]
Burkholder--Davis--Gundy and induction control the first term.
The induction hypothesis and Cauchy--Schwarz for
\(\int|f_{m-1}f_m|\) control the second. These arguments also give
the running-supremum estimate and continuity in each deterministic weight.
\end{proof}

\begin{lemma}[A minimum-variable estimate]\label{lem:minimum-weight}
Let \(m\geq2\), \(\psi,\phi\geq0\), and \(I=(a,b)\). Set
\(E(v)=\int_v^b\phi(u)^2\dd u\).
An ordered integral whose first weight is bounded by \(\psi\) in
absolute value and whose other weights are bounded by \(\phi\)
has \(L^2\) norm at most
\begin{equation}\label{eq:minimum-estimate}
 C_m\left[
 \left(\int_a^b\psi(v)^2E(v)^{m-1}\dd v\right)^{1/2}
 +\int_a^b\psi(v)\phi(v)E(v)^{(m-2)/2}\dd v
 \right].
\end{equation}
The first weight can be restricted to a subinterval.
The estimate allows arbitrary repeated Brownian coordinates.
\end{lemma}
\begin{proof}
For bounded smooth weights expand the second moment into Gaussian
pairings, including the contractions from Stratonovich integration.
Take absolute values of all deterministic factors. Removing the
coordinate-equality restrictions can only add nonnegative terms.
The resulting expression is bounded by the second moment of the
corresponding integral in a single Brownian coordinate, with weights
\(\psi,\phi,\ldots,\phi\). The latter integral equals
\[
 \frac1{(m-1)!}\int_a^b\psi(v)
       \left(\int_v^b\phi(u)\dd W_u\right)^{m-1}\circ\dd W_v.
\]
Reverse Brownian time at \(b\), writing
\(\overline W_r=W_b-W_{b-r}\) and
\(M_r=\int_0^r\phi(b-u)\dd\overline W_u\).
Apart from the fixed factor \(1/(m-1)!\), the scalar integral becomes
\[
 \int_0^{b-a}\psi(b-r)M_r^{m-1}\dd\overline W_r
 +\frac{m-1}{2}\int_0^{b-a}
       \psi(b-r)\phi(b-r)M_r^{m-2}\dd r.
\]
The adapted Gaussian martingale satisfies
\(\E M_{b-v}^2=E(v)\). It\^o isometry and Gaussian moments bound
the first integral by the square-root term in
\eqref{eq:minimum-estimate}; Minkowski's inequality bounds the
correction by its second term.
Approximation of weights, using the same estimates for differences,
gives the result whenever the displayed right side is finite.
\end{proof}

\Needspace{8\baselineskip}
\smallskip\noindent\emph{A degree-four example.}
For bounded nonnegative smooth weights on a finite interval, let
\(J_w\) denote the ordered
Stratonovich integral with weights \((\psi,\phi,\phi,\phi)\), and
\(I_w\) its uncontracted It\^o term. For \(a<x<y<b\), put
\[
 \begin{aligned}
 A(x,y)&=\phi(x)\phi(y)\int_a^x\psi(v)\phi(v)\dd v,\\
 B(x,y)&=\psi(x)\phi(y)E(y),\qquad
 D(x,y)=\psi(x)\phi(y)\int_x^y\phi(v)^2\dd v,\\
 c&=\int_a^b\psi(v)\phi(v)E(v)\dd v.
 \end{aligned}
\]
Writing \(I_{ii}(K)\) for the ordered double It\^o integral of
\(K\) in color \(i\), the adjacent contractions give
\[
 \begin{split}
 J_{1122}&=I_{1122}+\tfrac12 I_{22}(A)+\tfrac12 I_{11}(B)+\tfrac14c,\\
 J_{1111}&=I_{1111}+\tfrac12 I_{11}(A+B+D)+\tfrac14c.
 \end{split}
\]
The term \(D\) is the additional middle-pair contraction in the
scalar word. With \(Q=\|I_{1122}\|_2^2=\|I_{1111}\|_2^2\),
orthogonality and the It\^o isometry show
\[
 \E J_{1122}^2
 =Q+\tfrac14\int_{a<x<y<b}(A^2+B^2)\dd x\dd y+\tfrac1{16}c^2
 \leq \E J_{1111}^2,
\]
because \(A,B,D\geq0\). Signed weights bounded by
\((\psi,\phi,\phi,\phi)\) obey the same domination after taking
absolute values of their contracted kernels. In the reverse Brownian
time used above, the scalar term is exactly
\[
 J_{1111}=\frac16\int_0^{b-a}\psi(b-r)M_r^3\dd\overline W_r
 +\frac14\int_0^{b-a}\psi(b-r)\phi(b-r)M_r^2\dd r.
\]
Its two \(L^2\) bounds are
\(\frac{\sqrt{15}}6(\int_a^b\psi^2E^3)^{1/2}\) and
\(\frac{\sqrt3}4\int_a^b\psi\phi E\), as required by
\eqref{eq:minimum-estimate}.

\begin{lemma}[Valley tails and a finite-origin bound]\label{lem:valley-tails}
Put \(\kappa_t(u)=\mathfrak c_H(t-u)_+^{H-1/2}\).
Let \(X^R_{s,t}=(G_s^R)^{-1}G_t^R\), where \(G_t^R\) is the signature
of \(r\mapsto\int_{-R}^{r\wedge t}\kappa_t(u)\dd W_u\).
The following bounds hold for every \(1\leq q<\infty\).
For \(mH<1\), \(0\leq s<t\leq T\), and \(R'\geq R\geq2T+2\),
\begin{equation}\label{eq:valley-tail}
 \norm{\pi_m(X^{R'}_{s,t}-X^R_{s,t})}_{L^q}
 \leq C_{m,H,q,T}|t-s|R^{mH-1}.
\end{equation}
For \(mH\leq1\), fixed \(R<\infty\), and \(h=t-s\leq1\),
\begin{equation}\label{eq:finite-origin-bound}
 \norm{\pi_m X^R_{s,t}}_{L^q}
 \leq
 \begin{cases}
 C_{m,H,q,R,T}h^{mH},&mH<1,\\
 C_{m,H,q,R,T}h(1+|\log h|),&mH=1.
 \end{cases}
\end{equation}
\end{lemma}
\begin{proof}
Expand the inverse signature and telescope the product
\((G_s^R)^{-1}G_t^R\) at the least auxiliary integration time.
The result is the finite-cutoff valley formula
\eqref{eq:stationary-valley}. Each valley is a finite union of
ordered simplexes. On such a simplex the first weight is
\(\kappa_t-\kappa_s\) and every other weight is either
\(\kappa_s\) or \(\kappa_t\).
Apply Lemma~\ref{lem:minimum-weight} with
\(\psi=|\kappa_t-\kappa_s|\) and \(\phi=\kappa_s+\kappa_t\).

For \(v<-R\), the mean-value theorem gives
\[
 \psi(v)\leq C_H h|v|^{H-3/2},\quad
 \phi(v)\leq C_H|v|^{H-1/2},\quad
 E(v)\leq C_H|v|^{2H}.
\]
The two tail integrals in \eqref{eq:minimum-estimate} are thus
bounded by constants times
\[
 h\left(\int_R^\infty r^{2mH-3}\dd r\right)^{1/2},
 \qquad h\int_R^\infty r^{mH-2}\dd r.
\]
Both have order \(hR^{mH-1}\) when \(mH<1\).
Changing the cutoff changes exactly the terms whose minimum lies
in the added remote interval. This proves the \(L^2\) tail bound.

For the finite-origin bound, put \(v=s+hz\). Both integrals have
scaling factor \(h^{mH}\). The transformed lower endpoint is
at least \(-(T+R)/h\). The only unbounded remote powers are
\(|z|^{2mH-3}\) in the first squared integral and
\(|z|^{mH-2}\) in the second integral.
They are integrable if \(mH<1\). At \(mH=1\), they give at most
\(\sqrt{\log(2+(T+R)/h)}\) and \(\log(2+(T+R)/h)\), respectively.
The singularities at the finite endpoints are integrable because
\(H>0\). This gives \eqref{eq:finite-origin-bound}.
Each random variable is a sum of Wiener chaoses of degrees at most
\(m\); hypercontractivity gives every finite \(q\).
The case \(m=1\) follows directly from Wiener isometry.
\end{proof}

\section{The infinite-past construction away from critical levels}
\label{sec:causal-construction}
The moment estimates permit removal of the remote cutoff
when every retained degree satisfies \(kH<1\).

\begin{theorem}\label{thm:noncritical-causal}
Let \(0<H<1/2\) and \(NH<1\).
On the two-sided Wiener space, \(X^R_{s,t}\) converges in every finite
\(L^q\), through level \(N\), as \(R\to\infty\).
The limit on the whole real time line is group-like,
satisfies Chen's identity, and has first level
\[
 B^H_t-B^H_s=\int_{\R}(\kappa_t-\kappa_s)(u)\dd W_u.
\]
It has jointly stationary increments, degree-\(kH\) self-similarity,
and commutes with linear maps of the noise coordinates.
It is measurable from the Wiener noise up to time \(t\).
For every \(k\leq N\),
\begin{equation}\label{eq:causal-moments}
 \norm{\pi_k X_{s,t}}_{L^q}\leq C_{H,N,q}|t-s|^{kH}.
\end{equation}
It has a modification with group-metric H\"older exponent
\(\alpha<H\). If \(H^{-1}\notin\mathbb N\), the choice
\(N=\lfloor1/H\rfloor\) gives all required rough-path levels for
\(1/(N+1)<\alpha<H\).
\end{theorem}
\begin{proof}
The tail estimate \eqref{eq:valley-tail} proves the \(L^q\) Cauchy
property at each interval. All finite-cutoff objects are group-like
and multiplicative. Their polynomial Chen and shuffle identities
pass to the limit by H\"older's inequality, using higher moments.
The first level follows from Wiener isometry. Past measurability
passes to the \(L^q\) limit.

Write \((\theta_aW)(A)=W(A+a)\).
A change of the Brownian variable gives the functional identity
\[
 X^R_{s+a,t+a}(W)=X^{R+a}_{s,t}(\theta_aW).
\]
Cutoff convergence through real \(R\) passes this equality
to the limit. White-noise shift invariance proves joint
stationarity.
The construction initially proved on positive compact
intervals extends to any compact real interval by this
same finite-cutoff translation identity. Limits on
overlapping intervals agree, so no independent choice of
initial group value is introduced.
Brownian scaling and
\(\kappa_{\lambda t}(\lambda u)
=\lambda^{H-1/2}\kappa_t(u)\) similarly give
\[
 X^R_{\lambda s,\lambda t}(W)
       =\delta_{\lambda^H}X^{R/\lambda}_{s,t}(S_\lambda W),
\]
where \((S_\lambda W)(A)=\lambda^{-1/2}W(\lambda A)\)
and \(\delta\) denotes tensor dilation. The limiting
functional scaling identity follows as before.
Integer-shift identities hold on a common full-measure
set by rational endpoints, countable intersections and
continuity; in particular the stationary model is a
measurable factor of the unit white-noise shift.
Naturality under linear maps follows at finite cutoff from the
signature identity and then passes to the limit.

Finite moments on \([0,1]\), stationarity and scaling prove
\eqref{eq:causal-moments}. For a homogeneous metric on \(G^N\),
\(d_N(1,g)\leq C\sum_{k=1}^N|\pi_k g|^{1/k}\).
Arbitrarily high moments and Kolmogorov's criterion give every
H\"older exponent below \(H\). Define the continuous group path
from rational intervals; Chen's identity then holds everywhere.
Finally \(NH<1<(N+1)H\) for the stated nonreciprocal choice.
Weak geometricity gives geometricity at any slightly smaller
H\"older exponent by signature approximation.
\end{proof}

\begin{remark}
The cutoff enhancements in Theorem~\ref{thm:noncritical-causal}
are signatures in auxiliary Brownian time, rather than canonical
signatures of their first-level paths in terminal time.
The theorem therefore concerns this enhanced approximation scheme,
not the usual smooth approximations of fBm.
\end{remark}

\section{The critical construction}\label{sec:critical-main}
At reciprocal indices, the construction proceeds in three steps.
A fixed-reference subtraction gives a causal enhancement, and a lagged
calibration restores stationarity. The zero-lag limit then recovers
exact self-similarity and sharp moment bounds.

\subsection{A central subtraction at reciprocal Hurst indices}
\label{sec:critical-construction}
Let \(H=1/N\), \(N\geq4\), and work on the positive time half-line.
Fix the splitting point \(-1\). For \(R>1\), factor the frozen
signature as
\[
 G_t^R=H_t^R J_t,\qquad
 H_t^R=S_{\leq N}\left(r\mapsto\int_{-R}^{r\wedge(-1)}
                         \kappa_t(u)\dd W_u\right),
\]
where \(J_t\) is the corresponding signature over \([-1,t]\).
The notation \(H_t^R\) denotes a group element, not a Hurst parameter.
The path \(t\mapsto H_t^R\) is smooth for \(t\geq0\).
Set
\begin{equation}\label{eq:remote-logarithmic-derivative}
 \Omega_t^R=(H_t^R)^{-1}\partial_tH_t^R,\qquad
 Z_R=\pi_N\Omega_0^R.
\end{equation}
The logarithmic derivative belongs to the free Lie algebra.
In particular \(Z_R\) lies in its central homogeneous degree \(N\).
It is measurable with respect to Wiener noise before time \(-1\).

\begin{lemma}[Remote derivatives]\label{lem:remote-derivatives}
For \(k<N\), the process \(\pi_k\Omega^R\) converges in
the space \(L^q(\Omega;C^r([0,T]))\) for all finite \(q,r,T\).
At degree \(N\),
\begin{equation}\label{eq:remote-difference}
 \pi_N\Omega_t^R-\pi_N\Omega_0^R
\end{equation}
converges in the same spaces.
The cutoff tails of \(\partial_t\pi_N\Omega^R\) have size
\(O_{q,r,T}(R^{-H})\).
\end{lemma}
\begin{proof}
Differentiating the finite-domain valley identity at equal
terminal times gives
\[
 \pi_k\Omega_t^R
 =\sum_{j=1}^k(-1)^{j-1}
 \int_{A_j^k\cap[-R,-1]^k}
       \kappa'_t(u_j)\prod_{\ell\ne j}\kappa_t(u_\ell)
       \circ\dd W_{u_1}\cdots\circ\dd W_{u_k}.
\]
Coordinate labels are retained in each tensor word.
On this remote interval all kernels are smooth in the terminal
parameter, and differentiation follows from
Lemma~\ref{lem:mixed-weights}.
For \(k<N\), Lemma~\ref{lem:minimum-weight} with minimum weight
\(|\kappa'_t|\) gives a tail \(O(R^{kH-1})\).

At degree \(N\) differentiate once more in \(t\).
Terms with \(\kappa''_t\) at the minimum have tail \(O(R^{-1})\).
All other terms have \(\kappa'_t\) at the minimum and one
\(\kappa'_t\) at a later integration time.
Taking absolute values in Gaussian pairings, removing color
restrictions and summing the possible later positions reduces the
bound to a scalar Brownian integral with tail factor
\[
 \frac1{(N-2)!}U_v^{N-2}V_v,\qquad
 U_v=\int_v^{-1}\kappa_t(u)\dd W_u,\quad
 V_v=\int_v^{-1}|\kappa'_t(u)|\dd W_u .
\]
Let \(E(v)=\E U_v^2\), \(F(v)=\E V_v^2\).
For \(v<-2T-2\),
\(E(v)\leq C_T|v|^{2H}\) and \(F(v)\leq C_T\).
Backward It\^o conversion and Gaussian moments bound the remote
tail by a constant times
\begin{align*}
 &\left(\int_{-\infty}^{-R}|\kappa'_t(v)|^2
                 E(v)^{N-2}F(v)\dd v\right)^{1/2}\\
 &\quad+\int_{-\infty}^{-R}|\kappa'_t(v)|\kappa_t(v)
                 E(v)^{(N-3)/2}F(v)^{1/2}\dd v\\
 &\quad+\int_{-\infty}^{-R}|\kappa'_t(v)|^2
                 E(v)^{(N-2)/2}\dd v.
\end{align*}
The first two terms have order
\(R^{(N-1)H-1}=R^{-H}\); the last has order
\(R^{NH-2}=R^{-1}\).
The derivatives of \(U_v^{N-2}\) and \(V_v\) both contribute
to these bounds. The later derivative requires the mixed Gaussian
moment estimate above; a uniform remote factor \(R^{-1}\)
is unavailable.

Higher terminal-time derivatives have the same bound. Indeed,
on \(u\leq-1\), a derivative of order at least two is bounded by
a constant times \(|\kappa''_t|\), and a derivative of positive
order by a constant times \(|\kappa'_t|\), on compact \(t\)
intervals. Every differentiated term is thus controlled by
one of the two cases just considered.
Finite-chaos hypercontractivity gives all finite moments.
One-dimensional Sobolev estimates in the terminal parameter,
using one additional derivative, give the \(C^r\) bounds.
Finally integrate \(\partial_t\pi_N\Omega^R\) from zero to \(t\)
to obtain convergence of \eqref{eq:remote-difference}.
\end{proof}

\begin{lemma}[Transfer from a smooth remote factor]
\label{lem:remote-rough-transfer}
Let \(I\) be compact with an anchor \(t_*\in I\).
Let \(A_R,A\) be continuous random Lie-coefficient paths,
uniformly bounded in every finite \(L^q(C^0(I))\), with
\(\|A_R-A\|_{L^q(C^0)}\leq C_q\rho_R\), \(\rho_R\to0\).
Solve \(K_R'=K_RA_R\), \(K'=KA\), with common value \(1\)
at the anchor. Let \(J:I\to G^N(\R^d)\) be continuous,
with all finite moments of the coordinate suprema of \(J,J^{-1}\).
Put
\[
 U^R_{s,t}=J_s^{-1}K_R(s)^{-1}K_R(t)J_t,\qquad
 U_{s,t}=J_s^{-1}K(s)^{-1}K(t)J_t.
\]
For \(0<\alpha<1/N\), the unrooted coordinate distance
\eqref{eq:coordinate-distance} satisfies
\(\|d_{\alpha;I}(U^R,U)\|_{L^p}\leq C_p\rho_R\)
for every finite \(p\). No independence assumption is needed.
This is an estimate of the difference; individual rough regularity
also requires regularity of the common local factor.
\end{lemma}
\begin{proof}
The nilpotent ODE has coordinates that are finite sums of
ordinary ordered integrals of products of at most \(N\)
coefficient components. Expanding differences one factor at
a time and applying H\"older's inequality gives uniform finite
\(C^1\) moments of \(K_R,K\), and
\(\|K_R-K\|_{L^q(C^1)}\leq C_q\rho_R\).
Inversion is polynomial and preserves these estimates.
The increment difference
\(V^R_{s,t}-V_{s,t}=K_R(s)^{-1}K_R(t)-K(s)^{-1}K(t)\)
vanishes at \(t=s\). Integrating its terminal derivative gives
\[
 \|\pi_j(V^R_{s,t}-V_{s,t})\|
 \leq |t-s|\,P_R\|K_R-K\|_{C^1},\qquad j\geq1,
\]
where \(P_R\) is a polynomial in the displayed path and
inverse norms, with uniform finite moments of every order.
Multiplication by \(J_s^{-1},J_t\) introduces only their
coordinate suprema: the middle difference has zero degree-zero
component. Divide by \(|t-s|^{j\alpha}\), take the supremum,
and use \(j\alpha<1\) and high finite-moment H\"older estimates.
\end{proof}

\begin{theorem}[Forward-causal critical enhancement]\label{thm:critical-causal}
Let \(H=1/N\), \(N\geq4\), and define
\begin{equation}\label{eq:central-renormalised-model}
 \widehat X^R_{s,t}=X^R_{s,t}\exp(-(t-s)Z_R),
 \qquad 0\leq s<t.
\end{equation}
On each compact positive time interval this converges in probability
in every \(\alpha\)-H\"older rough-path topology with
\(1/(N+1)<\alpha<H\).
The limit is a geometric enhancement of the prescribed fBm first
level and is measurable from Wiener noise up to time \(t\).
It commutes with linear coordinate maps. Its lower levels satisfy
\eqref{eq:causal-moments}; its highest level satisfies, for
\(0<t-s\leq1\),
\begin{equation}\label{eq:critical-top-moment}
 \norm{\pi_N\widehat X_{s,t}}_{L^q}
 \leq C_{q,T}|t-s|(1+|\log|t-s||).
\end{equation}
This theorem does not assert stationarity or exact self-similarity
of the critical highest level.
\end{theorem}
\begin{proof}
Since \(Z_R\) is central, \eqref{eq:central-renormalised-model}
preserves the finite-cutoff group and Chen identities.
Define the remote group path
\[
 \widehat H_t^R=(H_0^R)^{-1}H_t^R\exp(-tZ_R).
\]
It starts at the identity, and its logarithmic derivative is
\(\Omega_t^R-Z_R\).
Lemma~\ref{lem:remote-derivatives} gives convergence of this Lie
algebra-valued coefficient in every finite \(C^r\)-moment topology.
The step-\(N\) group ODE
\(\partial_t\widehat H=\widehat H(\Omega-Z)\) is triangular in
homogeneous coordinates. Its solution is a finite polynomial of
iterated ordinary time integrals of the coefficients.
Consequently \(\widehat H^R\) converges in the same smooth
moment topologies to a group-valued path \(\widehat H\).

The local factor \(J_t\) is independent of \(R\), and centrality gives
\begin{equation}\label{eq:local-remote-product}
 \widehat X^R_{s,t}
 =J_s^{-1}(\widehat H_s^R)^{-1}\widehat H_t^R J_t .
\end{equation}
The finite-origin estimate \eqref{eq:finite-origin-bound}, applied
at cutoff \(1\), controls \(J_s^{-1}J_t\):
the degree \(k<N\) is \(O_{L^q}(h^{kH})\) and degree \(N\) is
\(O_{L^q}(h(1+|\log h|))\), where \(h=t-s\).
Subtract \(J_s^{-1}J_t\) from \eqref{eq:local-remote-product}.
Each positive-degree coordinate of the smooth middle increment
is \(O_{L^q,T}(h)\); coordinates of \(J_s^{-1}\) and \(J_t\)
on the compact interval have all finite supremum moments.
These follow from \eqref{eq:finite-origin-bound}, high-moment
group Kolmogorov estimates and the finite Gaussian-signature
moments of \(J_0\).
Every coordinate of the difference is therefore
\(O_{L^q,T}(h)\).
Since \(kH<1\) for \(k<N\), this proves the stated lower estimates
and \eqref{eq:critical-top-moment}.

Lemma~\ref{lem:remote-rough-transfer} now applies to the
normalized smooth remote ODE and the common local factor.
It gives convergence after division by \(h^{k\alpha}\), since
\(k\alpha<1\) for all \(k\leq N\), with the same coefficient
cutoff rate whenever that rate is specified.
Thus \eqref{eq:local-remote-product} converges in the coordinate
\(\alpha\)-H\"older rough-path metric, in probability (in fact in
every finite moment after fixing a slightly larger local exponent).
The logarithm in the highest-level estimate is absorbed by an
arbitrarily small H\"older loss. Group Kolmogorov estimates give
the continuous modification.
All polynomial group identities pass to the limit.
First-level values are unchanged by a degree-\(N\) subtraction.
Past measurability and linear covariance hold at every cutoff and
pass to the limit. Geometricity at a specified exponent follows by
first working at a slightly larger exponent below \(H\).
\end{proof}

Appendix~\ref{sec:slope-obstructions} classifies the convergent
time-linear central subtractions of this cutoff family and proves
their stationarity and scaling obstructions. The next section
uses a time-dependent calibration.

\subsection{Stationary critical enhancements by lagged calibration}
\label{sec:moving-lag}
A time-dependent, past-measurable calibration restores stationarity
to the critical construction of Section~\ref{sec:critical-construction}.
Fix \(H=1/N\), \(N\geq4\), and a deterministic lag \(\ell>0\).
For \(-R<a<x\), define
\[
 H_x^{R,a}=S_{\leq N}\left(u\mapsto
                \int_{-R}^{u\wedge a}\kappa_x(v)\dd W_v\right),
 \qquad
 \Omega_x^{R,a}=(H_x^{R,a})^{-1}\partial_x H_x^{R,a}.
\]
The derivative holds the auxiliary endpoint \(a\) fixed.
Only after differentiation do we put
\begin{equation}\label{eq:lag-density}
 z_R^\ell(r)=\pi_N\Omega_r^{R,r-\ell},\qquad
 Y^{R,\ell}_{s,t}=X^R_{s,t}
       \exp\left(-\int_s^t z_R^\ell(r)\dd r\right).
\end{equation}
On a fixed compact interval, \(R\) is large enough that
\(-R<r-\ell\). The density is Lie-valued of homogeneous degree
\(N\), hence central. It is not the total derivative of the
curve \(r\mapsto H_r^{R,r-\ell}\), whose endpoint has Brownian
variation. Its ordinary time integral in
\eqref{eq:lag-density} is well defined.
Logarithmic velocities and multiplication by a central path
are established rough-path operations
\citep{FrizOberhauser2009,BellingeriEtAl2022}.
We prove convergence of the lagged causal calibration as the
remote cutoff is removed.

\begin{theorem}[Stationary causal critical model]\label{thm:moving-lag}
For every fixed \(\ell>0\), \(Y^{R,\ell}\) converges through
arbitrary real cutoffs \(R\to\infty\), locally on the whole real
time line, in each \(\alpha\)-H\"older rough-path topology with
\[
 \frac1{N+1}<\alpha<H=\frac1N.
\]
The convergence holds in probability and in every finite moment
of the coordinate rough seminorms.
The limit \(Y^\ell\) is a geometric enhancement of the prescribed
two-sided fBm increments. It has jointly stationary increments,
is measurable from Wiener noise up to the right endpoint, and
commutes with linear coordinate maps.
Its levels below \(N\) are the infinite-past limits already
constructed. For every finite \(q\), uniformly in \(s\in\R\),
\begin{align}
 \|\pi_kY^\ell_{s,s+h}\|_{L^q}
       &\leq C_{H,q}h^{kH},\qquad k<N,\label{eq:lag-low-moments}\\
 \|\pi_NY^\ell_{s,s+h}\|_{L^q}
       &\leq C_{H,\ell,q}h(1+|\log h|),\quad0<h\leq1.
       \label{eq:lag-top-moments}
\end{align}
Every degree-\(k\) coordinate belongs to the sum of Wiener
chaoses of orders at most \(k\).
\end{theorem}

\subsubsection{Finite-band continuity and cancellation of the remote divergence}
Fix a compact interval \(I=[S,T]\), a reference \(r_*\in I\),
and \(a<S-\ell-1\). For \(a<b<x\), let \(A_x^b\) be the
frozen signature over the finite auxiliary interval \([a,b]\).
Put
\[
 A^\ell(r)=A_r^{r-\ell},\qquad
 q^\ell(r)=
 \left.(A_x^{r-\ell})^{-1}\partial_xA_x^{r-\ell}\right|_{x=r}.
\]

\begin{lemma}\label{lem:lag-band}
The curves \(A^\ell,(A^\ell)^{-1},q^\ell\) belong locally
to \(L^p(\Omega;C^\eta(I))\) for all finite \(p\) and
every \(0<\eta<1/2\).
\end{lemma}
\begin{proof}
The gap between terminal parameter and auxiliary endpoint is
\(\ell\); the weights and all needed terminal derivatives are
therefore bounded on the finite region under consideration.
For \(r'>r\), first change the terminal parameter at the fixed
endpoint \(r-\ell\), then extend the endpoint to \(r'-\ell\).
The first change is \(O_{L^p}(|r'-r|)\), by
Lemma~\ref{lem:mixed-weights} applied to one extra parameter
derivative. Concatenation at the smaller endpoint expresses the
second change using signatures over an interval of length
\(|r'-r|\). Each positive-degree term on that added interval
contains a weight with \(L^2\) norm \(O(|r'-r|^{1/2})\);
all remaining factors have uniformly bounded moments.
Thus this change is \(O_{L^p}(|r'-r|^{1/2})\).
The same argument applies to the terminal partial derivative
of the signature: differentiate the finite concatenation
identity before substituting the moving endpoint.
Polynomial inversion and multiplication preserve these
estimates. High-moment Kolmogorov estimates give the asserted
H\"older versions and all finite moments of their norms.
\end{proof}

Let \(\Omega^R(r)=\Omega_r^{R,a}\).
Concatenation at the fixed endpoint \(a\), followed by partial
differentiation at a fixed \(b\), gives
\[
 \Omega_x^{R,b}=(A_x^b)^{-1}\Omega_x^{R,a}A_x^b
                          +(A_x^b)^{-1}\partial_xA_x^b .
\]
Since degree \(N\) is central, evaluation at
\((x,b)=(r,r-\ell)\) yields the exact decomposition
\begin{equation}\label{eq:lag-decomposition}
 z_R^\ell(r)=\pi_N\Omega^R(r)
   +P_N(A^\ell(r),\pi_{<N}\Omega^R(r))+\pi_Nq^\ell(r),
\end{equation}
where \(P_N(A,u)=\pi_N(A^{-1}uA)\) for \(u\) of degrees below
\(N\). In particular no divergent top component occurs in
the polynomial \(P_N\).

The proof of Lemma~\ref{lem:remote-derivatives} applies with
fixed endpoint \(a\) and reference \(r_*\). Its remote powers
are unchanged, and the positive gap \(S-a\) controls the finite
end. It gives convergence in every finite smooth moment norm
of
\[
 \pi_{<N}\Omega^R(r),\qquad
 \pi_N\Omega^R(r)-\pi_N\Omega^R(r_*).
\]
The slowest cutoff rate is \(R^{-H}\).
Consequently Lemma~\ref{lem:lag-band} and
\eqref{eq:lag-decomposition} show that
\begin{equation}\label{eq:lag-centred-density}
 b_R^\ell(r):=z_R^\ell(r)-\pi_N\Omega^R(r_*)
       \longrightarrow b^\ell(r)
 \quad\text{in }L^p(\Omega;C^\eta(I)),\quad\eta<1/2.
\end{equation}
The convergence concerns the centered density, rather than
\(z_R^\ell\) itself.

\begin{proof}[Proof of Theorem~\ref{thm:moving-lag}]
Write \(Z_R^*=\pi_N\Omega^R(r_*)\).
The fixed-end proof of Theorem~\ref{thm:critical-causal},
with reference \(r_*\) and endpoint \(a\), gives convergence
of \(\widetilde X^R_{s,t}
=X^R_{s,t}\exp(-(t-s)Z_R^*)\) on \(I\).
The same group ODE integrates from \(r_*\) in either direction,
and Lemma~\ref{lem:remote-rough-transfer} applies on an
arbitrary compact real interval.
Centrality gives
\[
 Y^{R,\ell}_{s,t}=\widetilde X^R_{s,t}
             \exp\left(-\int_s^t b_R^\ell(r)\dd r\right).
\]
By \eqref{eq:lag-centred-density}, the last central correction
converges in every finite \(C^1\)-moment norm. Its degree-\(N\)
rough seminorm difference is bounded by
\[
 |I|^{1-N\alpha}\|b_R^\ell-b_{R'}^\ell\|_{C(I)}.
\]
This proves rough convergence, since \(N\alpha<1\).
Only the highest level changes in the last multiplication.
The finite-origin bound and the smooth remote-factor argument
in Theorem~\ref{thm:critical-causal} therefore give
\eqref{eq:lag-top-moments} locally, while lower levels are
unchanged. Continuous versions on the countable exhaustion
\([-m,m]\) agree on overlaps as limits of the same cutoff
family. This gives a whole-line model.
Chen and shuffle pass to the limit; geometricity follows by
first working at a slightly larger exponent below \(H\).

At every cutoff \(z_R^\ell(r)\) is measurable from noise
before \(r-\ell\). Thus the time integral on \([s,t]\) is
\(\mathcal F^W_{t-\ell}\)-measurable, and the whole increment
is \(\mathcal F_t^W\)-measurable.
Coordinate naturality holds for signatures and their
partial logarithmic derivatives before taking limits.
The finite-chaos assertion follows from the Stratonovich
contraction formula and closure in \(L^2\).

For a fixed shift \(c\), let \((\theta_cW)(A)=W(A+c)\).
The exact cutoff identities are
\[
 X^R_{s+c,t+c}(W)=X^{R+c}_{s,t}(\theta_cW),\qquad
 z_R^\ell(r+c;W)=z_{R+c}^\ell(r;\theta_cW).
\]
Changing variables in the time integral and using convergence
through real cutoffs gives
\begin{equation}\label{eq:lag-shift}
 Y^\ell_{s+c,t+c}(W)=Y^\ell_{s,t}(\theta_cW).
\end{equation}
For each fixed shift this is an equality of random models
on compact intervals. White-noise invariance proves joint
stationarity. Shifting the moment estimate on \([0,1]\)
then gives its uniform-in-position version.
\end{proof}

\subsubsection{Scaling and change of lag}
Let \((S_\lambda W)(A)=\lambda^{-1/2}W(\lambda A)\).
Kernel homogeneity and Brownian scaling give
\[
 H_{\lambda x}^{R,\lambda a}(W)
       =\delta_{\lambda^H}H_x^{R/\lambda,a}(S_\lambda W).
\]
Partial differentiation contributes a factor \(\lambda^{-1}\).
Since \(NH=1\), the density \(z_R^\ell\) has scaling degree
zero; its time integral has scaling degree one. Hence
\begin{equation}\label{eq:lag-scaling}
 Y^\ell_{\lambda s,\lambda t}(W)
       =\delta_{\lambda^H}Y^{\ell/\lambda}_{s,t}(S_\lambda W).
\end{equation}
This is covariance of a family with a lag parameter, rather
than self-similarity of a fixed member.

\begin{proposition}\label{prop:lag-change}
For \(\ell,L>0\),
\[
 K^{\ell,L}(r)=\lim_{R\to\infty}
                         (z_R^\ell(r)-z_R^L(r))
\]
exists locally in every finite \(L^p(C^\eta)\), \(\eta<1/2\).
It is a stationary, continuous, degree-\(N\) Lie-valued process,
measurable from noise up to \(r-\min(\ell,L)\), with
coordinates in Wiener chaoses of orders at most \(N\).
It satisfies
\begin{equation}\label{eq:lag-change}
 Y^\ell_{s,t}=Y^L_{s,t}
                  \exp\left(-\int_s^tK^{\ell,L}(r)\dd r\right),
 \qquad K^{\ell,L}+K^{L,M}=K^{\ell,M}.
\end{equation}
For \(d\geq2\), \(w=1^{N-1}2\), and \(0<\ell<L\),
\begin{equation}\label{eq:lag-nondegenerate}
 \|(K^{\ell,L}(0))^w\|_2^2
 \geq\frac1{(N-1)!}\int_{-L}^{-\ell}|\kappa'_0(v)|^2
       \left(\int_v^{-\ell}\kappa_0(u)^2\dd u\right)^{N-1}\dd v>0.
\end{equation}
\end{proposition}
\begin{proof}
Use the common endpoint \(a\) in
\eqref{eq:lag-decomposition}. Subtracting the two formulas
cancels the top remote term at the same time \(r\).
All remaining terms converge by the lower remote estimates
and Lemma~\ref{lem:lag-band}.
The stated symmetries, measurability, chaos orders and
identities pass from their cutoff counterparts.
For the last assertion project onto colored chaos
\((N-1,1)\), with the color-\(2\) variable
\(-L<v<-\ell\) and every color-\(1\) variable between \(v\)
and \(-\ell\). The density \(z_R^L(0)\) has no noise
variables there. Only the last valley of \(z_R^\ell(0)\)
contributes. Its symmetrized kernel is
\[
 \frac{(-1)^{N-1}}{(N-1)!}\kappa'_0(v)
                  \prod_{j=1}^{N-1}\kappa_0(u_j).
\]
Colored isometry gives \eqref{eq:lag-nondegenerate}.
The domain stays away from zero, so its norm is finite and
strictly positive. Lower-chaos contractions cannot cancel it.
\end{proof}

\subsubsection{Different lag selections give different laws}
Let \(\mathcal E_+\) be the Polish space of continuous
\(G^N\)-valued Chen increment families on
\(\{0\leq s\leq t<\infty\}\), with local uniform topology.
Its Borel sigma-field is generated by increment evaluations.
Write \(T_1\mathbf y=(\mathbf y_{s+1,t+1})_{s,t}\).

\begin{theorem}\label{thm:lag-laws}
Assume \(d\geq2\). For any two distinct positive lags,
the laws of \(Y^\ell,Y^L\) on every interval \([0,T]\),
\(T>0\), are different.
Their laws on \(\mathcal E_+\), and consequently on the
whole real line, are mutually singular.
Every such law is mixing under \(T_1\); all have the same
first-level fBm law.
No singularity on a fixed compact horizon is asserted.
\end{theorem}
\begin{proof}
Put \(w=1^{N-1}2\). For \(0<h<\ell/4\), project onto
colored chaos \((N-1,1)\) generated by noise after
\(-\ell/2\). The counterterm in \eqref{eq:lag-density}
is measurable before \(h-\ell<-\ell/2\), so this projection
kills it. Restrict further to a color-\(2\) minimum
\(v\in(-\ell/2,-h)\), with all other variables in \((v,0)\).
The last valley and colored isometry give
\begin{align}
 \|(Y^\ell_{0,h})^w\|_2^2
 &\geq\frac1{(N-1)!}\int_{-\ell/2}^{-h}
     |\kappa_h(v)-\kappa_0(v)|^2
       \left(\int_v^0\kappa_0(u)^2\dd u\right)^{N-1}\dd v\notag\\
 &\geq c\,h^2\log\frac{\ell}{2h}.
 \label{eq:lag-UV}
\end{align}
The last estimate uses
\(|\kappa_h(-r)-\kappa_0(-r)|\geq c\,h r^{H-3/2}\)
for \(r\geq h\), and \(NH=1\).

Relabel so that \(L=\lambda\ell\), \(\lambda>1\), and
let \(v_L(h)=\E|(Y^L_{0,h})^w|^2\).
Equation~\eqref{eq:lag-scaling} gives
\(v_L(\lambda h)=\lambda^2v_\ell(h)\).
If the two laws on \([0,T]\) agreed, then for
\(h_m=T\lambda^{-m}\),
\[
 v_L(h_m)/h_m^2=v_L(h_{m-1})/h_{m-1}^2=v_L(T)/T^2.
\]
This contradicts \eqref{eq:lag-UV}.
Some \(w\)-coordinate second moment on this sequence
therefore differs. Truncating its square yields a bounded
continuous cylinder function \(F\) with different expectations.

The map from white noise to \(Y^\ell\) is measurable as a
probability limit of cutoff models.
Equation~\eqref{eq:lag-shift}, first for integer shifts and
rational endpoints and then by continuity, makes each law
a factor of the unit white-noise shift.
Lemma~\ref{lem:ergodic-singularity} therefore gives mixing
and ergodicity. Applied to the separating cylinder \(F\),
the same lemma proves mutual singularity on the infinite
half-line and hence on the whole-line space.
\end{proof}

The calibration has a precise finite-cutoff normalization.
If all noise after \(a=r-\ell\) is frozen out, the
degree-\(N\) infinitesimal logarithmic velocity of
\[
 (H_r^{R,a})^{-1}H_{r+h}^{R,a}\exp(-hc)
\]
vanishes exactly when \(c=z_R^\ell(r)\).
This condition fixes the lagged response convention within the
cutoff construction; it does not determine a unique enhancement
among all possible selections.

\subsection{A scale-free critical limit}\label{sec:zero-lag}
Fix $H=1/N$, $N\geq4$. The positive-lag models of
Section~\ref{sec:moving-lag} admit a zero-lag limit.
The convergence concerns integrated central corrections;
it does not require convergence of their densities.
All rates below refer to the unrooted coordinate distance
\eqref{eq:coordinate-distance} on the compact interval \(I\)
under consideration; we abbreviate \(d_{\alpha;I}\) by \(d_\alpha\).
The quantitative proof of Theorem~\ref{thm:moving-lag}, at
reference lag $1$, gives for every finite $p$
\begin{equation}\label{zl:eq:basecutoff}
 \|d_\alpha(Y^{R,1},Y^1)\|_{L^p}\leq C_{I,\alpha,p}R^{-H}.
\end{equation}
Indeed the slowest remote coefficient tail in
Lemma~\ref{lem:remote-derivatives} is $R^{-H}$.
The triangular smooth group ODE and local--remote multiplication
preserve that rate in every finite moment of the coordinate
rough seminorms. We now prove cancellation between comparable
lags, uniformly in their ratio.

\smallskip\noindent\emph{Proof strategy.}
Lemma~\ref{zl:lem:mean} supplies centering, and
Lemma~\ref{zl:lem:band} and Proposition~\ref{zl:prop:localization}
give relative-window approximations to comparable-lag differences.
Disjoint centered windows yield Lemma~\ref{zl:lem:covariance};
time integration gives Proposition~\ref{zl:prop:integrated}.
Interpolation and Lemma~\ref{zl:lem:GRR} turn this into summable
little-H\"older bounds, used in Theorem~\ref{zl:thm:zero} for every
deterministic positive lag. The functional symmetries of this limit
then give sharp moments in Corollary~\ref{zl:cor:symmetry}.
The absolute-cutoff bound in
Lemma~\ref{zl:lem:joint-band}, together with
\eqref{zl:eq:basecutoff}, yields Theorem~\ref{zl:thm:joint}.

\subsubsection{Zero mean of the logarithmic partial velocity}

For $a<b<x$, let $H_x^{[a,b]}$ denote the auxiliary-time
signature with weight $\kappa_x$ on $[a,b]$.
Write $\Omega_x^{[a,b]}=(H_x^{[a,b]})^{-1}\partial_xH_x^{[a,b]}$,
with both endpoints held fixed under differentiation.

\begin{lemma}[Mean-zero Maurer--Cartan velocity]\label{zl:lem:mean}
For fixed finite endpoints $a<b$ and terminal parameter $x>b$,
every homogeneous coordinate of $\Omega_x^{[a,b]}$ has mean zero.
Consequently $\E z_R^\ell(r)=0$.  The same assertion holds when
the finite auxiliary interval is $[r-M,r-\ell]$, provided its
endpoints are fixed before taking the partial derivative.
\end{lemma}

\begin{proof}
First use a finite piecewise-linear auxiliary Brownian approximation
with scalar deterministic weights $f_j(x)$ on its segments.  Its
signature is
\[
 G_x=\prod_{j=1}^m\exp(V_j(x)),\qquad
 V_j(x)=f_j(x)\Delta W_j,
\]
where $\Delta W_j$ are independent centered vector increments.
The common scalar weight across all spatial coordinates is essential:
$V'_j=f'_j\Delta W_j$ commutes with $V_j$.
Thus $\exp(-V_j)\partial_x\exp(V_j)=V'_j$ and the exact
product rule gives
\[
 G_x^{-1}\partial_xG_x
   =\sum_{j=1}^m
       \operatorname{Ad}\left((\exp(V_{j+1})\cdots\exp(V_m))^{-1}\right)
                       (f'_j(x)\Delta W_j).
\]
The tail in the $j$th summand is independent of $\Delta W_j$.
The adjoint action is linear in its displayed argument and a finite
polynomial in the truncated tail coordinates.  Each summand is
therefore integrable and centered, in every homogeneous degree.

It remains to justify the partial-derivative limit, rather than
differentiate a moving endpoint.  On any compact neighborhood of
the given $x>b$, the weights $\kappa_x$ and $\partial_x\kappa_x$ are
bounded on $[a,b]$ and smooth in $x$.  Step approximations in
auxiliary time converge for both weights in deterministic $L^2$.
Each coordinate of $\partial_xH_x$ is the finite sum of mixed
Stratonovich integrals obtained by inserting one
$\partial_x\kappa_x$ weight.  The mixed-weight $L^p$ estimate makes
these derivatives continuous under the step approximation.
For step weights, polygonal approximation gives convergence of
these mixed integrals as well: they are coordinates of the joint
signature of the Gaussian martingale with weights
$(\kappa_x,\partial_x\kappa_x)$.  This is ordinary auxiliary Brownian
signature approximation and follows from the finite
It\^o--Stratonovich recursion, including its contractions.
Polynomial signature inversion and multiplication then give
$G_x^{-1}\partial_xG_x\to\Omega_x^{[a,b]}$ in every finite
$L^p$.  Its expectation is zero.
Evaluating afterwards at $x=r$ and endpoints depending on the
deterministic value of $r$ changes none of this argument.
\end{proof}

\subsubsection{Uniform lag-band tails and relative localization}

For $1\leq\vartheta\leq2$ define
\[
 K_R^{\ell,\vartheta}(r)
       =z_R^\ell(r)-z_R^{\vartheta\ell}(r),
 \qquad
 K^{\ell,\vartheta}(r)
       =\lim_{R\to\infty}K_R^{\ell,\vartheta}(r).
\]
The estimates also include the case $\vartheta=1$.
The following tail estimate proves existence of the coefficient
limit directly from subcritical-degree remote estimates.

\begin{lemma}[Uniform normalized band estimate]\label{zl:lem:band}
For every finite $p$, the limits $K^{1,\vartheta}(0)$ exist in
$L^p$, and, uniformly for $1\leq\vartheta\leq2$ and $A\geq4$,
\begin{equation}\label{zl:eq:normalizedtail}
 \norm{K_A^{1,\vartheta}(0)-K^{1,\vartheta}(0)}_{L^p}
                  \leq C_p A^{-H},\qquad
 \norm{K^{1,\vartheta}(0)}_{L^p}\leq C_p.
\end{equation}
The finite and limiting coefficients are centered and belong to
the sum of Wiener chaoses of orders at most $N$.
\end{lemma}

\begin{proof}
Fix the remote/local split $a=-3$ and let
$L_b=H_0^{[a,b]}$, $b\in[-2,-1]$, and
$\omega_b=L_b^{-1}\partial_xH_x^{[a,b]}|_{x=0}$.
The logarithmic derivative product identity gives
\[
 \Omega_0^{[-A,b]}
    =L_b^{-1}\Omega_0^{[-A,a]}L_b+\omega_b.
\]
In the difference between $b=-1$ and $b=-\vartheta$, the
degree-$N$ component of the remote logarithmic derivative cancels,
because it is central.  Put
$\Lambda_A=\pi_{<N}\Omega_0^{[-A,a]}$.
Then the exact remaining expression is
\begin{align*}
 K_A^{1,\vartheta}(0)=\pi_N\bigl[
 L_{-1}^{-1}\Lambda_A L_{-1}
 -L_{-\vartheta}^{-1}\Lambda_A L_{-\vartheta}
                     +\omega_{-1}-\omega_{-\vartheta}\bigr].
\end{align*}
Every degree of $\Lambda_A$ is below $N$, so Lemma~\ref{lem:remote-derivatives} gives its
convergence, with worst error $A^{(N-1)H-1}=A^{-H}$.
The local signatures, inverses and partial derivatives have all
finite running-supremum moments for $b\in[-2,-1]$; all weights and
their derivatives are bounded on $[-3,-1]$.
Polynomial multiplication and H\"older's inequality therefore give
\eqref{zl:eq:normalizedtail}, with constants independent of
$\vartheta$.
Lemma~\ref{zl:lem:mean} shows that each finite coefficient and its
$L^2$ limit are centered. At finite cutoff the degree-$N$ velocity is a sum
of chaoses of orders at most $N$.  Their sum is closed in $L^2$,
so the same assertion holds for the limit.
\end{proof}

For $M\geq4\ell$ define a relative-window coefficient by
\begin{equation}\label{zl:eq:localization}
 K_{\mathrm{loc}}^{\ell,\vartheta;M}(r)
   =\left.\pi_N\left(
      \Omega_x^{[r-M,r-\ell]}
       -\Omega_x^{[r-M,r-\vartheta\ell]}\right)\right|_{x=r}.
\end{equation}
The lower endpoint $r-M$, like both upper endpoints, is held fixed
under $\partial_x$ before $x=r$ is substituted.

\begin{proposition}[Relative localization]\label{zl:prop:localization}
The coefficient in \eqref{zl:eq:localization} is centered and measurable
from white noise on $[r-M,r-\ell]$ alone.  For every finite $p$,
uniformly in $r\in\R$, $\ell>0$, $M\geq4\ell$, and
$1\leq\vartheta\leq2$,
\begin{equation}\label{zl:eq:localization-tail}
 \norm{K^{\ell,\vartheta}(r)
             -K_{\mathrm{loc}}^{\ell,\vartheta;M}(r)}_{L^p}
                  \leq C_p(\ell/M)^H,
 \qquad
 \norm{K^{\ell,\vartheta}(r)}_{L^p}\leq C_p.
\end{equation}
For an absolute lower cutoff $-R$, whenever $R+r\geq4\ell$,
\begin{equation}\label{zl:eq:absolutetail}
 \norm{K_R^{\ell,\vartheta}(r)
                   -K^{\ell,\vartheta}(r)}_{L^p}
                  \leq C_p\ell^H(R+r)^{-H}.
\end{equation}
\end{proposition}

\begin{proof}
Translation to $r=0$ followed by scaling the auxiliary time by
$\ell$ gives the exact identity
\begin{equation}\label{zl:eq:band-scaling}
 K_R^{\ell,\vartheta}(r;W)
  =K_{(R+r)/\ell}^{1,\vartheta}
                  (0;S_\ell\theta_rW).
\end{equation}
The degree-$N$ logarithmic partial derivative has scaling degree
$NH-1=0$; this is why there is no amplitude factor.
The same calculation applied to the relative interval
$[r-M,r-\ell]$ gives
\[
 K_{\mathrm{loc}}^{\ell,\vartheta;M}(r;W)
       =K_{M/\ell}^{1,\vartheta}(0;S_\ell\theta_rW).
\]
Passing to the coefficient limit in
\eqref{zl:eq:band-scaling} gives the corresponding identity for
$K^{\ell,\vartheta}$.  White-noise shift and dilation preserve
law.  Lemma~\ref{zl:lem:band} therefore gives
\eqref{zl:eq:localization-tail} and \eqref{zl:eq:absolutetail}.
The support assertion follows from the two finite auxiliary
intervals; centering follows from Lemma~\ref{zl:lem:mean}.
\end{proof}

\subsubsection{Covariance decay and an integrated band estimate}

The coefficient processes have stationary laws, as follows either
from their finite-cutoff shift identities and real-cutoff convergence
or directly from \eqref{zl:eq:band-scaling}.
All norms and covariances below may be taken in any fixed coordinates
of the finite-dimensional space $\mathfrak L_N$.

\begin{lemma}[A uniform covariance bound]\label{zl:lem:covariance}
Uniformly in $r,s\in\R$, $\ell>0$ and $1\leq\vartheta\leq2$,
\begin{equation}\label{zl:eq:covariance}
 \left\|\E\left[K^{\ell,\vartheta}(r)
                      \otimes K^{\ell,\vartheta}(s)\right]\right\|
       \leq C\min\left(1,\left(\frac\ell{|r-s|}\right)^H\right),
\end{equation}
with the right-hand side interpreted as $C$ when $r=s$.
\end{lemma}

\begin{proof}
Let $d=|r-s|$.  If $d\geq12\ell$, use $M=d/3$ in
\eqref{zl:eq:localization}.  The two relative auxiliary windows are
disjoint, so their localized coefficients are independent.
Both are centered, so their cross moment is zero.
Subtract these localized variables from the full coefficients.
Cauchy--Schwarz, their uniform $L^2$ bounds, and
\eqref{zl:eq:localization-tail} bound the remaining cross moment by
$C(\ell/M)^H=C3^H(\ell/d)^H$.
For $d<12\ell$, use the uniform $L^2$ bound.  Enlarging the
constant by $12^H$ gives \eqref{zl:eq:covariance} throughout.
\end{proof}

Set
\[
 A^{\ell,\vartheta}_{st}
              =\int_s^tK^{\ell,\vartheta}(r)dr.
\]
The fixed-positive-lag construction gives continuous coefficient
versions, but only their pointwise moment estimates will be used.

\begin{proposition}[Integrated band cancellation]\label{zl:prop:integrated}
For every finite $p$, uniformly in $s<t$, $\ell>0$ and
$1\leq\vartheta\leq2$, with $h=t-s$,
\begin{equation}\label{zl:eq:integratedmin}
 \norm{A^{\ell,\vartheta}_{st}}_{L^p}
       \leq C_p\min\bigl(h,\ell^{H/2}h^{1-H/2}\bigr).
\end{equation}
Consequently, for $0\leq\varepsilon\leq H/2$,
\begin{equation}\label{zl:eq:interpolated}
 \norm{A^{\ell,\vartheta}_{st}}_{L^p}
                     \leq C_p\ell^\varepsilon h^{1-\varepsilon}.
\end{equation}
\end{proposition}

\begin{proof}
At $p=2$, integrate \eqref{zl:eq:covariance} over the square
$[s,t]^2$.  For $h\leq\ell$ the result is $Ch^2$.
For $h\geq\ell$ it is bounded by
\[
 C\int_0^h(h-u)\min(1,(\ell/u)^H)du
          \leq C\ell^Hh^{2-H},
\]
because $0<H<1$.  Taking square roots gives
\eqref{zl:eq:integratedmin} at $p=2$.
Every coordinate of $K(r)$ is in the sum of chaoses of orders
at most $N$.  Bochner time integration stays in that closed
finite-chaos subspace.  Finite-chaos hypercontractivity gives
\eqref{zl:eq:integratedmin} for all finite $p$ with constants depending
on $p,N$ and no lag loss.
For $h\leq\ell$ the first term in the minimum is at most
$\ell^\varepsilon h^{1-\varepsilon}$; for $h\geq\ell$ the
second term is at most that expression when $\varepsilon\leq H/2$.
This proves \eqref{zl:eq:interpolated}.
\end{proof}

\subsubsection{The H\"older function-space step}

Let $I=[-T,T]$, $T\geq1$.  Write $c_0^a(I;\mathfrak L_N)$ for the
anchored little-H\"older space, the closure of smooth paths $F$ with
$F(0)=0$ in the $a$-H\"older norm, $0<a<1$.
It is a separable Banach space. Band primitives are $C^1$ and hence
belong to it. We use this separable space for Bochner $L^p$
integration.

\begin{lemma}[Moment-to-H\"older estimate]\label{zl:lem:GRR}
Suppose continuous anchored random paths $F^\ell$ satisfy, for every
sufficiently large finite $q$,
\[
 \norm{F^\ell(t)-F^\ell(s)}_{L^q}
             \leq C_q b_\ell |t-s|^\nu,
\]
uniformly on $I$.  If $0<a<\nu$ and $p<\infty$, then
\[
 \norm{F^\ell}_{L^p(C^a(I))}\leq C_{p,a,\nu,I}b_\ell.
\]
It suffices to use a $q\geq p$ with $q(\nu-a)>1$.
\end{lemma}

\begin{proof}
Choose $a+1/q<\sigma<\nu$.
The expected $q$th power of the fractional Sobolev seminorm
with smoothness $\sigma$ is bounded by
\[
 C_qb_\ell^q\int_I\int_I
          |t-s|^{q(\nu-\sigma)-1}ds\,dt<\infty.
\]
The one-dimensional fractional Sobolev/GRR embedding into
$C^{\sigma-1/q}$ gives the $C^a$ bound.  Anchoring at zero
controls the supremum norm, and passing from $L^q$ to $L^p$
proves the assertion.
\end{proof}

Fix $1/(N+1)<\alpha<H$ and put $a=N\alpha<1$.
Choose
\begin{equation}\label{zl:eq:epsilon}
 0<\varepsilon<\min(H/2,1-N\alpha).
\end{equation}
Apply Lemma~\ref{zl:lem:GRR} to
$F^{\ell,\vartheta}(t)=A^{\ell,\vartheta}_{0t}$ using
\eqref{zl:eq:interpolated}, with $\nu=1-\varepsilon>a$.
It gives
\begin{equation}\label{zl:eq:Cband}
 \norm{F^{\ell,\vartheta}}_{L^p(c_0^a(I))}
                      \leq C_{p,I,\alpha,\varepsilon}\ell^\varepsilon,
 \qquad 1\leq\vartheta\leq2.
\end{equation}
This interpolation is necessary near $N\alpha=1$.
Using only the exponent $1-H/2$ in
\eqref{zl:eq:integratedmin} would not cover every $\alpha<H$.

\subsubsection{Zero-lag convergence of the limiting models}

The fixed-lag comparison identity is
\begin{equation}\label{zl:eq:lag-comparison}
 Y^\ell_{st}=Y^{\vartheta\ell}_{st}
                    \exp(-A^{\ell,\vartheta}_{st}).
\end{equation}
It holds at every finite cutoff, then in the fixed-lag limit.

\begin{theorem}[Zero-lag limit]\label{zl:thm:zero}
There is a whole-line enhancement $Y^0$ such that, on every compact
interval, for every $1/(N+1)<\alpha<H$, finite $p$, and
$\varepsilon$ as in \eqref{zl:eq:epsilon},
\begin{equation}\label{zl:eq:zerorate}
 \norm{d_\alpha(Y^\ell,Y^0)}_{L^p}
                         \leq C\ell^\varepsilon,
                         \qquad0<\ell\leq1.
\end{equation}
The parameter $\ell$ is deterministic in this bound.  The limit is
geometric for every such $\alpha$, has all required $N$ levels, and
retains the same levels below $N$ as the original critical model.
Its coordinates have every finite moment on compact intervals.
\end{theorem}

\begin{proof}
Let $\ell_j=2^{-j}$.  All lag differences are additive central
increments.  By \eqref{zl:eq:lag-comparison}, their degree-$N$
coordinates are the negatives of the increments of
$F^{\ell_j,2}$, and all lower coordinates agree.
By \eqref{zl:eq:Cband},
\[
 \sum_{j=1}^\infty
       \norm{F^{\ell_j,2}}_{L^p(c_0^a(I))}
             \leq C\sum_{j=1}^\infty2^{-j\varepsilon}<\infty.
\]
Thus the sum converges in this Banach space.  If $F$ is its limit,
define
\[
 Y^0_{st}=Y^1_{st}\exp(-(F(t)-F(s))).
\]
The tail of the series is $O(\ell_j^\varepsilon)$ in
$L^p(c_0^a)$, proving \eqref{zl:eq:zerorate} along dyadic lags.
The $N\alpha$-H\"older estimate applies to the additive central
differences summed in this series.

For any $0<\ell\leq1$, choose $j$ with
$\ell\leq\ell_j<2\ell$.  Equation
\eqref{zl:eq:lag-comparison}, now with
$\vartheta=\ell_j/\ell\in[1,2)$, and the uniform estimate
\eqref{zl:eq:Cband} give
$\norm{d_\alpha(Y^\ell,Y^{\ell_j})}_{L^p}
\leq C\ell^\varepsilon$.
The dyadic tail proves \eqref{zl:eq:zerorate} for every real positive
lag.

Every difference summed above is Lie-valued in degree $N$ and
additive, so the resulting multiplication preserves Chen, shuffle,
and all lower levels.  Its $N\alpha$-H\"older estimate is exactly
what is required at the top level.  Geometricity at a given exponent
follows by first doing the argument at a slightly larger exponent
below $H$ and then taking the usual smaller-exponent geometric
closure.  The limits on different compact intervals and different
admissible exponents agree by uniqueness of convergence in
probability at each increment.  A countable exhaustion gives a
consistent whole-line version.  The fixed reference model has all
finite coordinate moments, and the Banach-space series has all
finite moments, proving the last assertion.
\end{proof}

\subsubsection{The zero-lag model's symmetries and sharp scaling moments}

\begin{corollary}\label{zl:cor:symmetry}
The model $Y^0$ is causal for the complete Wiener past, jointly
stationary, exactly $H$-self-similar in the graded group sense,
and natural under fixed linear coordinate maps.  For fixed $b\in\R$
and $c>0$ the functional identities are
\[
 Y^0_{s+b,t+b}(W)=Y^0_{st}(\theta_bW),\qquad
 Y^0_{cs,ct}(W)=\delta_{c^H}Y^0_{st}(S_cW).
\]
They imply, for every $1\leq j\leq N$ and finite $p$,
\begin{equation}\label{zl:eq:sharp-moments}
 \norm{\pi_jY^0_{st}}_{L^p}=C_{j,p}|t-s|^{jH},
\end{equation}
with the equality interpreted with the corresponding fixed
$[0,1]$-coordinate norm.  In particular the top bound is $C_{N,p}h$,
without a logarithm.
\end{corollary}

\begin{proof}
For each fixed lag, $Y^\ell_{st}$ is $\mathcal F_t^W$-measurable.
Its pointwise convergence in probability proves the same assertion
for $Y^0$ in the completed filtration.
Pass to the lag limit in \eqref{eq:lag-shift}.
For scaling, pass to the limit in \eqref{eq:lag-scaling}:
both $\ell$ and $\ell/c$ tend to zero, and white-noise dilation
preserves the probability law used in the $L^p$ estimates.
This proves the displayed functional identities.  Spatial
naturality similarly passes from the fixed-lag models.
Unit-interval moments are finite by Theorem~\ref{zl:thm:zero}.
Stationarity and scaling therefore give
\eqref{zl:eq:sharp-moments}.  Its top-degree consequence uses $NH=1$.
Proposition~\ref{prop:complete-past} identifies the Wiener past with
the complete two-sided fBm increment past. Causality therefore holds
for this complete past, which includes information beyond a finite
observed segment $[0,t]$.
\end{proof}

\subsubsection{Uniform absolute-cutoff comparison and the joint limit}

The following estimate uses pointwise coefficient control followed by
time integration.

\begin{lemma}[Integrated absolute-cutoff error]\label{zl:lem:joint-band}
On $I=[-T,T]$, for $R\geq2T+8$, $0<\ell\leq1$,
$1\leq\vartheta\leq2$, and $a=N\alpha<1$,
\begin{equation}\label{zl:eq:joint-band}
 \left\|t\longmapsto\int_0^t
       (K_R^{\ell,\vartheta}-K^{\ell,\vartheta})(r)dr
                 \right\|_{L^p(c_0^a(I))}
             \leq C_{p,I,a}\ell^H R^{-H}.
\end{equation}
\end{lemma}

\begin{proof}
For $r\in I$, $R+r\geq R/2\geq4\ell$.
Equation \eqref{zl:eq:absolutetail} gives, for every finite $q$,
\[
 \norm{K_R^{\ell,\vartheta}(r)-K^{\ell,\vartheta}(r)}_{L^q}
                       \leq C_q\ell^H R^{-H}.
\]
Minkowski's inequality gives increment bounds
$C_q h\ell^H R^{-H}$ for their time integrals.
Apply Lemma~\ref{zl:lem:GRR} with $\nu=1$ and sufficiently high $q$.
This proves \eqref{zl:eq:joint-band}.
\end{proof}

\begin{theorem}[Joint cutoff and zero-lag convergence]\label{zl:thm:joint}
On each compact interval, for every finite $p$ and admissible
$\alpha,\varepsilon$ as above, there is a constant independent of
$0<\ell\leq1$ and sufficiently large $R$ such that
\begin{align}
 \sup_{0<\ell\leq1}
        \norm{d_\alpha(Y^{R,\ell},Y^\ell)}_{L^p}
                         &\leq C R^{-H},\label{zl:eq:uniform-R}\\
 \norm{d_\alpha(Y^{R,\ell},Y^0)}_{L^p}
                         &\leq C(R^{-H}+\ell^\varepsilon).
                                      \label{zl:eq:joint}
\end{align}
The supremum in \eqref{zl:eq:uniform-R} is outside the $L^p$ norm.
These bounds prove joint convergence along every deterministic pair
$(R,\ell)\to(\infty,0)$; they do not assert a pathwise uniform
statement for a lag chosen after observing the noise.
\end{theorem}

\begin{proof}
For dyadic lags $\ell_j=2^{-j}$, telescope the central differences
relative to lag $1$ both at finite cutoff and at the limiting
cutoff.  Their difference in degree $N$ is the reference-cutoff
error plus increments of the path
\[
 -\sum_{j=1}^m\int_0^\cdot
       (K_R^{\ell_j,2}-K^{\ell_j,2})(r)dr.
\]
Lemma~\ref{zl:lem:joint-band} and
$\sum_j\ell_j^H<\infty$ bound this sum uniformly in $m$ by
$CR^{-H}$ in $L^p(c_0^a)$.
The lower degrees do not depend on lag, so their error is already
contained in the reference-lag bound \eqref{zl:eq:basecutoff}.
This proves \eqref{zl:eq:uniform-R} on dyadic lags.

For a general lag $\ell$, compare it to
$\ell_j\in[\ell,2\ell)$ as in Theorem~\ref{zl:thm:zero}.
Subtract the two finite/infinite comparison identities.
Their additional error is the integral of
$K_R^{\ell,\ell_j/\ell}-K^{\ell,\ell_j/\ell}$, bounded by
$C\ell^H R^{-H}$ by Lemma~\ref{zl:lem:joint-band}.
The constants are uniform over the ratio in $[1,2]$.
This proves \eqref{zl:eq:uniform-R} for all deterministic lags.
Finally combine it with \eqref{zl:eq:zerorate} by the triangle
inequality to obtain \eqref{zl:eq:joint}.
\end{proof}

\begin{corollary}[The reverse order of limits]\label{zl:cor:reverse}
For every fixed finite $R$, $Y^{R,\ell}$ has a zero-lag limit
$Y^{R,0}$ locally in $\alpha$ rough topology on the open domain
$(-R,\infty)$.  On each fixed compact interval and for sufficiently
large $R$,
\[
 \norm{d_\alpha(Y^{R,0},Y^0)}_{L^p}\leq C R^{-H}.
\]
Thus the iterated limits in either order agree with the joint limit.
At finite $R$ its first level is the finite-cutoff Gaussian path,
not the limiting fBm path.
\end{corollary}

\begin{proof}
Let $I\Subset(-R,\infty)$, enlarged inside this domain to include
zero if necessary, and
$\delta_R=\inf_{r\in I}(R+r)>0$.
For all sufficiently small lags $\ell<\delta_R/4$,
\eqref{zl:eq:absolutetail} and the proof of
Lemma~\ref{zl:lem:joint-band} give a bound
$C\ell^H\delta_R^{-H}$ for each integrated band error.
The limiting band has bound $C\ell^\varepsilon$ by
\eqref{zl:eq:Cband}.  The sum of the two bounds is dyadically
summable, so the same central Banach-space series gives
$Y^{R,0}$, with arbitrary-lag convergence supplied by ratios
in $[1,2]$.  Group identities, rough regularity and geometry pass
as before.  For large $R$ on a fixed $I$, take $\ell\downarrow0$
in \eqref{zl:eq:uniform-R}; convergence in the indicated Banach
coordinate spaces gives the asserted bound.
\end{proof}

\section{Critical extension rigidity with fixed lower levels}
\label{sec:critical-rigidity}
At a reciprocal index, the scaling degree of an additive
top-level difference is one. Combined with whole-line causality,
this determines the difference after the lower levels and
common-noise symmetries are fixed.

\begin{theorem}[Critical top-degree rigidity]\label{thm:critical-rigidity}
Let \(H=1/N\). Suppose that \(\bX,\widetilde{\bX}\) are
continuous whole-line geometric step-\(N\) enhancements on the
same two-sided Wiener space. Assume that all their tensor
coordinates at fixed increments are square-integrable, their
levels \(1,\ldots,N-1\) agree, and both are causal for
\(\mathcal F_t^W\). For each fixed \(b\in\R,c>0\), assume
the noise-functional identities
\[
 X_{s+b,t+b}(W)=X_{s,t}(\theta_bW),\qquad
 X_{cs,ct}(W)=\delta_{c^H}X_{s,t}(S_cW),
\]
and the same identities for \(\widetilde{\bX}\).
Then there is a deterministic degree-\(N\) Lie element \(D\)
such that, simultaneously for all intervals,
\begin{equation}\label{eq:critical-rigidity}
 \widetilde{\bX}_{s,t}=\bX_{s,t}\exp((t-s)D).
\end{equation}
The shift and scaling identities hold jointly with the same noise;
separate marginal symmetries do not suffice for the conclusion.
\end{theorem}
\begin{proof}
By Lemma~\ref{lem:central-fibre}, the increment
\[
 A_{s,t}:=\pi_N(\widetilde{\bX}_{s,t}-\bX_{s,t})
\]
is an additive degree-\(N\) Lie-valued increment and
\(\widetilde{\bX}_{s,t}=\bX_{s,t}\exp(A_{s,t})\).
The assumed common-noise identities imply joint stationarity
and degree-one self-similarity of \(A\), since \(NH=1\).
Put \(a_t=A_{0,t}\) using oriented increments, and
\(M=\E\|a_1\|^2<\infty\). For all \(s,t\in\R\),
\[
 \E\|a_t\|^2=t^2M,\qquad
 \E\|a_t-a_s\|^2=(t-s)^2M.
\]
For negative \(t\), use \(a_t=-A_{t,0}\) and stationarity.
Polarization gives the raw second-moment identity
\[
 \E\langle a_s,a_t\rangle=stM,\qquad
 \E\|a_t-ta_1\|^2=0.
\]
Rational times and continuity give a common event on which \(A_{s,t}=(t-s)D\), with \(D=A_{0,1}\).

For every positive integer \(m\), causality gives
\(D=A_{-m-1,-m}\in\mathcal F_{-m}^W\).
Thus \(D\) is independent of the noise on \((-m,\infty)\).
These future sigma-fields increase to the full white-noise
sigma-field. Since \(D\) is itself measurable from that noise,
it is independent of itself and hence deterministic.
Independence is unchanged by completion of the sigma-fields.
The drift lies in the closed degree-\(N\) Lie space.
\end{proof}

\begin{corollary}[A centered critical representative]
\label{cor:critical-centering}
For a model in Theorem~\ref{thm:critical-rigidity}, define
\begin{equation}\label{eq:critical-centering}
 m_X=\E[\pi_N\log\bX_{0,1}],\qquad
 \overline{\bX}_{s,t}=\bX_{s,t}\exp(-(t-s)m_X).
\end{equation}
This is well defined under the stated coordinate \(L^2\)
assumption. It preserves the lower levels, causality,
geometric rough regularity, and the functional shift and
critical scaling identities. Moreover
\[
 \E[\pi_N\log\overline{\bX}_{s,t}]=0.
\]
Any two extensions satisfying the theorem have the same
centered representative. In particular, the lower levels of \(Y^0\)
have a unique centered extension within this class.
\end{corollary}
\begin{proof}
For a group-like tensor and a word \(w\) of length \(N\),
\[
 (\log g)^w=\sum_{k=1}^N\frac{(-1)^{k+1}}k
   \sum_{\substack{w=w_1\cdots w_k\\|w_j|>0}}
                          g^{w_1}\cdots g^{w_k}.
\]
Shuffle identities express each scalar product on the right
as a finite linear combination of level-\(N\) coordinates.
Thus \(\pi_N\log g\) is a universal linear function of
\(\pi_Ng\), and \(m_X\) exists with no extra product-moment
assumption. Its value belongs to the degree-\(N\) Lie space.

The correction is deterministic, additive and central.
It preserves the rough estimates because \(N\alpha<1\),
and geometricity follows by approximation at a slightly larger
exponent. Its shift identity is immediate. Critical scaling
is preserved since
\[
 \delta_{c^H}\exp(-(t-s)m_X)=\exp(-c(t-s)m_X).
\]
Stationarity and scaling imply
\(\E[\pi_N\log\bX_{s,t}]=(t-s)m_X\).
Centrality gives
\[
 \pi_N\log\overline{\bX}_{s,t}
                  =\pi_N\log\bX_{s,t}-(t-s)m_X.
\]
Finally \eqref{eq:critical-rigidity} implies
\(m_{\widetilde X}=m_X+D\), so the centered models agree.
\end{proof}

Centering is natural for transported model laws. If \(L\) is
a deterministic linear spatial map, logarithms commute with
its tensor action and
\[
 m_{L_*X}=L^{\otimes N}m_X,\qquad
 \overline{L_*\bX}=L_*\overline{\bX}.
\]
The calibrating mean must be transported with the noise law.
For an orthogonal map on standard Wiener noise, the law is
unchanged and \(m_X\) is invariant whenever the original model
is orthogonally equivariant.
The normalization subtracts the possibly nonzero mean \(m_X\).
Uniqueness concerns extensions with fixed lower levels and the
stated common-noise identities; the first-level covariance alone
does not determine the integral.

\section{Stochastic integration and differential equations}\label{sec:integration-main}
\subsection{Controlled integration}\label{sec:controlled-integration}
Let \(\bX\) be a geometric step-\(N\) \(\alpha\)-H\"older rough path
over a path \(B\), where \(1/(N+1)<\alpha<1/N\), and set \(\rho=(N+1)\alpha>1\).
We specify its integration domain through controlled coefficients
\citep{Gubinelli2004,Hairer2014}. This includes integrands that depend
on the enhanced history, as well as smooth one-forms of the first level.

\begin{definition}\label{def:integration-jet}
An integration jet over \(\bX\) is a collection of locally
bounded continuous coefficients \(c_w(s)\), indexed by nonempty
words of length at most \(N\), such that, locally in \(s,t\),
\begin{equation}\label{eq:integration-jet}
 \left|c_v(t)-\sum_{0\leq|u|\leq N-|v|}
                        c_{uv}(s)X_{s,t}^u\right|
 \leq C|t-s|^{(N+1-|v|)\alpha},\qquad1\leq|v|\leq N.
\end{equation}
Here \(uv\) denotes concatenation with \(u\) as the prefix.
The empty prefix has \(X_{s,t}^{\varnothing}=1\).
The first-degree coefficients \(c_i\) are the integrand;
the higher coefficients are included in the jet data.
For a random integration jet, each \(c_w(s)\) is also
required to be measurable at every fixed time, with
continuous sample paths; the displayed bounds hold almost
surely with locally finite random constants.
\end{definition}

\begin{theorem}[Controlled stochastic integration]
\label{thm:controlled-integration}
For each integration jet, the compensated sums
\begin{equation}\label{eq:general-compensated-sums}
 \sum_{[u,v]\in\pi}\ \sum_{1\leq|w|\leq N}
                              c_w(u)X_{u,v}^w
\end{equation}
converge as the partition mesh tends to zero, to the unique
additive increment \(J\) satisfying
\begin{equation}\label{eq:word-integral-jet}
 J_{s,t}=Q_{s,t}+O(|t-s|^\rho),\qquad
 Q_{s,t}=\sum_{1\leq|w|\leq N}c_w(s)X_{s,t}^w.
\end{equation}
The integral depends continuously on the rough model and
on the coefficient and controlled-remainder norms in
\eqref{eq:integration-jet}.
For random integration jets over the enhancements constructed
here, the integral is measurable;
if every \(c_w(s)\) is adapted to the complete increment past,
then \(J_{s,t}\) is adapted to its right endpoint.
\end{theorem}
\begin{proof}
Chen's identity gives
\[
 Q_{s,t}-Q_{s,r}-Q_{r,t}
 =\sum_{1\leq|v|\leq N}
 \left(\sum_{0\leq|u|\leq N-|v|}
         c_{uv}(s)X_{s,r}^u-c_v(r)\right)X_{r,t}^v.
\]
Thus its defect is \(O(|t-s|^\rho)\).
Sewing, including its estimate for differences, proves
existence, uniqueness, convergence of the sums, and
continuity in the stated norms. The same estimates apply
to inverse increments and reversed endpoints.
For a random jet, deterministic partitions make every
compensated sum measurable, hence also its pathwise limit.
Continuity and rational endpoints give measurability as a
random continuous additive functional. If the jet is adapted,
every term in \eqref{eq:general-compensated-sums} over
\([s,t]\) is measurable at time \(t\); the limit has the same
property in the completed filtration.
\end{proof}

\begin{corollary}[One-forms]\label{cor:one-form-jet}
Assume that \(\bX\) enhances \(B\). For \(F\in C_b^{N+1}\) and \(w=i_1\cdots i_m\), set
\[
 c_w(s)=(\partial_{i_1}\cdots
                  \partial_{i_{m-1}}F_{i_m})(B_s).
\]
These are an integration jet and give the canonical one-form
rough integral for the specified enhancement.
For forward adaptedness on an interval beginning at \(a\),
use the recentered path \(B_r-B_a\), \(r\geq a\).
Locally smooth one-forms are also allowed by localization to
the compact sample-path range.
\end{corollary}
\begin{proof}
Symmetry of the derivative indices and the geometric shuffle
identity give, for every fixed nonempty word \(v\),
\begin{equation}\label{eq:shuffle-Taylor}
 \sum_{|u|=k}c_{uv}(s)X_{s,t}^u
     =\frac1{k!}D^kc_v(B_s)[B_{s,t}^{\otimes k}].
\end{equation}
Taylor's remainder through order \(N-|v|\) proves
\eqref{eq:integration-jet}. This also explains why no extra
factorial occurs in the coefficients \(c_w\).
The sums are the usual one-form rough sums.
\end{proof}

\subsubsection{Regularity-structure interpretation}
The controlled coefficients also describe a modelled distribution.
Attach to each nonempty word \(w\) a derivative symbol \(\dot\tau_w\)
of degree \(|w|\alpha-1\), and put \(\gamma=\rho-1>0\).
For the word model induced by \(\bX\), the function
\[
 f(s)=\sum_{1\leq|w|\leq N}c_w(s)\dot\tau_w
\]
belongs to \(\mathcal D^\gamma\), with reconstruction
\(\mathcal Rf=\partial_tJ\).
Thus sewing and reconstruction give the same integral on this domain.
Appendix~\ref{sec:word-model} constructs the positive and negative word
sectors and proves this identification in
Proposition~\ref{prop:word-reconstruction}. The analytic input is the
chosen enhancement and its rough-path estimates.

\subsection{The selected stochastic calculus}\label{sec:calculus}
The constructed enhancement now gives integration and equations on
its complete two-sided input history. The controlled domain is
specified in Section~\ref{sec:controlled-integration}.

\begin{theorem}\label{thm:selected-calculus}
Let \(0<H\leq1/4\), \(N=\lfloor1/H\rfloor\), and
\(1/(N+1)<\alpha<H\). Use the infinite-past enhancement at a
nonreciprocal index and the centered zero-lag enhancement
\(\overline Y^0\) of Corollary~\ref{cor:critical-centering} otherwise.
This is a whole-line geometric enhancement \(\bB^H\), adapted to
\(\mathcal F_t^B\), with jointly stationary increments and exact
self-similarity. For every finite \(p\) and \(k\leq N\),
\[
 \|\pi_k\bB^H_{s,t}\|_{L^p}\leq C_{k,p}|t-s|^{kH}.
\]
The following pathwise assertions hold on every finite time interval.
For adaptedness, integrate forward from an anchor \(a\in\R\),
using \(B_r^H-B_a^H\), \(r\geq a\), and initial data measurable
with respect to \(\mathcal F_a^B\). With anchor zero:
\begin{enumerate}
\item For every \(F\in C_b^{N+1}(\R^d;(\R^d)^*)\), the rough
 integral \(\int_0^\cdot F(B_t^H)\dd\bB_t^H\) is well-defined,
 measurable and adapted. Its controlled coefficients are induced
 by the spatial derivatives of \(F\), as in
 Corollary~\ref{cor:one-form-jet}.
\item The geometric chain rule holds. If \(F=\dd U\), its integral
 is \(U(B_t^H)-U(B_0^H)\) whenever the required derivatives exist
 on the sample-path range.
\item For bounded \(\operatorname{Lip}^{\gamma}\) vector fields
 \(V_1,\ldots,V_d\), with \(\gamma>1/\alpha\), the equation
 \[
  \dd Y_t=\sum_iV_i(Y_t)\dd\bB_t^{H,i},\qquad Y_0=y_0,
 \]
 has a pathwise unique global rough solution on every finite
 time interval. The solution is adapted for deterministic initial
 data or for initial data measurable at time zero.
\item The enhanced cutoff approximations in the chosen
 renormalization converge in probability in the rough-path topology.
 Their one-form integrals and equation solutions converge in
 probability uniformly on compact time intervals.
\end{enumerate}
At reciprocal indices, Theorem~\ref{zl:thm:joint} quantifies the joint
cutoff and zero-lag convergence. Each positive-lag model \(Y^\ell\)
is also stationary and causal; its scaling law is
\eqref{eq:lag-scaling}, with the lag rescaled together with time.
\end{theorem}
\begin{proof}
Theorem~\ref{thm:noncritical-causal}, Theorem~\ref{zl:thm:zero} and
Corollary~\ref{cor:critical-centering} give the chosen drivers,
moments and symmetries. Proposition~\ref{prop:complete-past}
identifies the complete fBm and Wiener pasts.

Theorem~\ref{thm:controlled-integration} and
Corollary~\ref{cor:one-form-jet} give the one-form integral and its
adaptedness. The geometric chain rule and rough-path solution theory
\citep{Lyons1998} give the first three assertions pathwise.
On an interval beginning at \(a\), the solution map uses only the
initial data and enhanced increments through the current time.
This proves adaptedness with the stated recentering convention.

For nonreciprocal cutoff convergence, use the fixed remote/local
splitting of Section~\ref{sec:critical-construction} without a
central subtraction. Every retained remote logarithmic derivative
converges because \(kH<1\). The smooth group ODE and the
local--remote product give rough-topology convergence, as in the
critical proof. At reciprocal indices, use
Theorem~\ref{zl:thm:joint} and apply the same limiting deterministic
centering to the approximations. Continuity of rough integration and
of the solution map proves the last assertion. These probability
convergence statements alone do not imply solution moment bounds.
\end{proof}

\begin{remark}[Finite observed history]\label{rem:finite-history}
Observations of fBm on \([0,T]\) do not determine a coupling to
an infinite negative-time history. The finite-origin enhancement of
\citet{NualartTindel2011}, in its canonical Volterra representation,
gives an adapted selection on that observed interval.
The same controlled integration and RDE theory apply to that selection.
Corollary~\ref{cor:finite-history-stationarity} shows that enhanced
increments cannot also be jointly stationary under this
finite-observation adaptation requirement when \(H\leq1/4\).
\end{remark}

\begin{corollary}[Stationary driver and solution cocycle]
\label{cor:rough-cocycle}
Use the stationary enhanced law on its canonical whole-line increment
space, with shift
\((\vartheta_s\mathbf x)_{u,v}=\mathbf x_{u+s,v+s}\).
For autonomous vector fields satisfying
Theorem~\ref{thm:selected-calculus}, the solution maps satisfy
\[
 \varphi(t+s,\mathbf x,y)
   =\varphi(t,\vartheta_s\mathbf x,
                        \varphi(s,\mathbf x,y)),\qquad s,t\geq0.
\]
They form a continuous solution cocycle for each chosen model.
\end{corollary}
\begin{proof}
The canonical shifts are defined pathwise, and stationarity preserves
their law. Restriction and concatenation of enhanced paths, together
with uniqueness of the rough equation, give the identity.
Continuous dependence gives continuity in the initial state.
This is the stationary-rough-path construction of
\citet[Theorem~5]{BailleulRiedelScheutzow2017}.
\end{proof}

\section{Adapted bounded-variation realization}\label{sec:causal-approximation}
The selected calculus can also be realized by ordinary
Riemann--Stieltjes integrals and ODEs. The approximating
controls encode the chosen higher increments through small
horizontal excursions. We use established rough-path
approximation geometry \citep{LyonsVictoir2007,FrizOberhauser2009}
with a one-cell delay to retain past adaptedness.

\begin{lemma}[Finite horizontal realization]\label{lem:horizontal-selector}
There is a Borel family of piecewise linear paths
\(\gamma_g:[0,1]\to\R^d\), \(g\in G^N(\R^d)\), with
a fixed finite number of coordinate-axis segments, such that
\[
 \gamma_g(0)=0,\qquad S_N(\gamma_g)_{0,1}=g,\qquad
 \operatorname{Lip}(\gamma_g)\leq C_{N,d}
                  \max_{1\leq j\leq N}\|\pi_jg\|^{1/j}.
\]
\end{lemma}
\begin{proof}
Choose in each homogeneous Lie degree a basis of bracket
words in \(e_1,\ldots,e_d\), and fix their order.
For a bracket word \(b\), define a horizontal group word
\[
 Q_{e_i}(u)=\exp(ue_i),\qquad
 Q_{[b_1,b_2]}(u)=Q_{b_1}(u)Q_{b_2}(u)
                         Q_{b_1}(u)^{-1}Q_{b_2}(u)^{-1},
 \quad u\geq0 .
\]
Concatenation and reversed traversal realize it by finitely
many axis segments. Truncated BCH gives no logarithmic
degree below \(|b|\), and degree \(|b|\) equal to
\(u^{|b|}b\).

Start with \(h=1\). Once all degrees below \(k\) match \(g\),
write the degree-\(k\) residual as
\[
             \pi_k\log(h^{-1}g)=\sum_{|b|=k}a_b b.
\]
Append \(Q_b(|a_b|^{1/k})^{\operatorname{sgn}(a_b)}\)
in the fixed basis order. These corrections preserve all
lower degrees and add exactly the residual in degree \(k\).
Their higher-degree terms are corrected at later stages.
After degree \(N\) the signature equals \(g\) exactly.
Keep zero-length segments in their fixed slots.

All segment coefficients are Borel functions of \(g\),
obtained by finite group operations, linear coordinates,
signs and roots. They are bounded on bounded target sets.
The construction commutes with positive tensor dilation:
\(g\mapsto\delta_rg\) multiplies every segment length by \(r\).
The set where \(\max_j\|\pi_jg\|^{1/j}=1\) is compact.
Homogeneity therefore bounds total segment length by the
displayed norm. Give the fixed number of segments equal
time slots. This gives the Lipschitz bound and a Borel
vertex list, hence a Borel path family.
\end{proof}

The infimum of variation over horizontal realizations defines
the Carnot--Carath\'eodory distance \(d_{\rm CC}\).
The lemma and the factorial signature estimate give its
equivalence to the homogeneous norm in the lemma.
Concatenation and reversal give its metric triangle
inequality.

\begin{theorem}[Adapted ordinary-signature approximation]
\label{thm:adapted-BV}
Let \(0<H\leq1/4\), \(N=\lfloor1/H\rfloor\),
\(1/(N+1)<\alpha<\beta<H\), and \(0<T<\infty\).
Let \(\bX\) be a whole-line geometric step-\(N\) enhancement
of \(B^H\), causal for \(\mathcal F_t^B\), such that, on every
compact real interval \(I\),
\[
 \|\pi_jX_{s,t}\|_{L^q}\leq C_{I,j,q}|t-s|^{jH},
 \qquad 1\leq j\leq N,\quad 1\leq q<\infty.
\]
The selected model of Theorem~\ref{thm:selected-calculus}
satisfies these assumptions. For \(0<\varepsilon\leq1\)
put \(g_k=X_{(k-1)\varepsilon,k\varepsilon}\), \(k\geq0\).
Starting at \(b^\varepsilon(0)=0\), concatenate
\begin{equation}\label{eq:delayed-BV}
 b^\varepsilon(t)=b^\varepsilon(k\varepsilon)
    +\gamma_{g_k}\bigl((t-k\varepsilon)/\varepsilon\bigr),
 \qquad k\varepsilon\leq t\leq(k+1)\varepsilon .
\end{equation}
These ordinary piecewise linear paths are adapted to
\(\mathcal F_t^B\) for \(t\geq0\). Their canonical signatures
satisfy
\begin{equation}\label{eq:delayed-grid-signature}
       S_N(b^\varepsilon)_{0,k\varepsilon}
                              =X_{-\varepsilon,(k-1)\varepsilon}.
\end{equation}
For the unrooted coordinate distance \eqref{eq:coordinate-distance},
one has, for every finite \(p\),
\begin{equation}\label{eq:BV-approximation-rate}
 \|d_{\alpha;[0,T]}(S_N(b^\varepsilon),\bX)\|_{L^p}
                         \leq C_{p,T,\alpha,\beta}\varepsilon^{\beta-\alpha}.
\end{equation}
The constant also depends on \(H,N,d\) and the displayed model moment bounds.
Convergence also holds almost surely as
\(\varepsilon\downarrow0\).
For \(F\in C_b^{N+1}\), the ordinary integrals
\(\int_0^\cdot F(b^\varepsilon_r)\dd b^\varepsilon_r\)
converge almost surely, uniformly on compacts, to the
rough integral for \(\bX\). The ordinary ODEs driven by
\(b^\varepsilon\), with the vector fields and initial data
of Theorem~\ref{thm:selected-calculus}, converge in the same
sense to its RDE solutions.
\end{theorem}
\begin{proof}
Each \(g_k\) is known at time \(k\varepsilon\), so its
entire next-cell control is measurable at that time.
The Borel family and concatenation give adaptedness.
The first cell uses \(X_{-\varepsilon,0}\), available
at time zero. Chen's identity proves
\eqref{eq:delayed-grid-signature}. In particular the
first-level grid values are
\(B_{(k-1)\varepsilon}-B_{-\varepsilon}\).
Their adaptedness follows from measurability of these increments
with respect to the complete past.

Fix a sample with finite
\[
 K=\sup_{-1\leq u<v\leq T+1}
            \frac{d_{\rm CC}(1,X_{u,v})}{|v-u|^\beta}.
\]
Inside a cell, horizontal variation bounds the signature
distance by \(CK\varepsilon^{\beta-1}|t-s|\).
For \(|t-s|\leq\varepsilon\) this is at most
\(CK|t-s|^\beta\). If such an interval crosses a cell
boundary, split it into two pieces and use the same
variation estimate and the metric triangle inequality.
Split any larger interval into its
two partial cells and full-grid middle. The latter is
exactly a delayed increment of \(X\), by
\eqref{eq:delayed-grid-signature}. The metric triangle
inequality gives a uniform \(\beta\)-rough bound \(CK\).

Write \(X^{\rm del}_{s,t}=X_{s-\varepsilon,t-\varepsilon}\).
For each \(j\leq N\), both the interpolation error
\(S_N(b^\varepsilon)-X^{\rm del}\) and the shift error
\(X^{\rm del}-X\), with \(h=t-s\), obey
\begin{equation}\label{eq:BV-coordinate-error}
 \|\pi_j(U_{s,t}-V_{s,t})\|
 \leq CK^j
 \begin{cases}
     h^{j\beta},&h\leq\varepsilon,\\
     \varepsilon^\beta h^{(j-1)\beta},&h\geq\varepsilon.
 \end{cases}
\end{equation}
The first case follows from the uniform rough bounds.
For the second, interpolation and the delayed path
have identical full-grid middle increments. Expand their
three Chen factors: every difference term contains a
positive-degree boundary error of size
\(CK^i\varepsilon^{i\beta}\), while the remaining
degree is bounded by \(CK^{j-i}h^{(j-i)\beta}\).
For the shift use the exact oriented-increment identity
\[
 X_{s-\varepsilon,t-\varepsilon}
           =X_{s-\varepsilon,s}X_{s,t}X_{t,t-\varepsilon}.
\]
Subtract \(X_{s,t}\); every remaining term has positive
outer degree and the same bound.

Since
\[
 (j-1)\beta-j\alpha
 <\frac{j-1}{N}-\frac{j}{N+1}<0,
\]
division of \eqref{eq:BV-coordinate-error} by \(h^{j\alpha}\)
gives \(CK^j\varepsilon^{j(\beta-\alpha)}\) in both cases.
Thus the unrooted distance is bounded by
\(C(1+K)^N\varepsilon^{\beta-\alpha}\), deterministically.
The model has all finite moments of \(K\):
its increment estimates give
\(\|d_{\rm CC}(1,X_{s,t})\|_{L^q}\leq C_q|t-s|^H\);
high-moment metric Kolmogorov or dyadic chaining gives
\(K\in L^p\) for every finite \(p\), choosing
\(q(H-\beta)>1\).
This proves \eqref{eq:BV-approximation-rate} and
pathwise convergence. Continuity of rough one-form
integration and the RDE solution map proves the ordinary
integral and ODE conclusions.
\end{proof}

\begin{corollary}[Finite-innovation realization of the critical limit]
\label{cor:critical-ordinary-realization}
Let \(H=1/N\), \(N\geq4\), and fix
\(1/(N+1)<\alpha<\beta<H\), \(0<T<\infty\), and
\(1\leq p<\infty\). Use
\[
 A^{R,\ell}_{s,t}
       =Y^{R,\ell}_{s,t}\exp(-(t-s)m),\qquad
 A^0=\overline Y^0,\qquad
 m=\E[\pi_N\log Y^0_{0,1}].
\]
Apply \eqref{eq:delayed-BV} to \(A^{R,\ell}\) to obtain
ordinary adapted controls \(b^{\varepsilon,R,\ell}\).
For any
\(0<\sigma<\min\{H/2,1-N\alpha\}\), large \(R\), and
\(0<\ell,\varepsilon\leq1\),
\begin{equation}\label{eq:ordinary-critical-joint}
 \|d_{\alpha;[0,T]}(S_N(b^{\varepsilon,R,\ell}),A^0)\|_{L^p}
           \leq C\bigl(\varepsilon^{\beta-\alpha}
                                  +R^{-H}+\ell^\sigma\bigr).
\end{equation}
At time \(t\), each finite-parameter control uses only
the innovation noise on \([-R,t]\).
The choice \(\varepsilon_m=\ell_m=2^{-m}\),
\(R_m=4+2^m\) gives almost sure rough convergence on every
compact interval, and ordinary integral and ODE limits
as in Theorem~\ref{thm:adapted-BV}.
\end{corollary}
\begin{proof}
Apply Theorem~\ref{zl:thm:joint} on \([-1,T+1]\).
At exponent \(\beta<H\), choose any positive
\(\sigma_\beta<\min\{H/2,1-N\beta\}\).
Its bound gives uniform finite moments of the
\(\beta\)-rough norms of \(A^{R,\ell}\) for large \(R\)
and \(0<\ell\leq1\). The same deterministic centering
is used in both models and cancels from their
coordinate difference.
Indeed, homogeneous norm equivalence bounds these CC seminorms by
\[
 C\left(\|A^0\|_{\beta,\mathrm{CC}}
       +\sum_{j=1}^N d_{\beta}(A^{R,\ell},A^0)^{1/j}\right).
\]
These seminorms enter only through their uniform moment bounds.
The deterministic argument in the proof of
Theorem~\ref{thm:adapted-BV} therefore gives
the first term of \eqref{eq:ordinary-critical-joint}
with a uniform constant.
At exponent \(\alpha\), the joint cutoff theorem gives
the other two terms with the stated full range of
\(\sigma\), independently of \(\sigma_\beta\).
For sufficiently large \(R\) all intervals
\([-\varepsilon,0]\) and moving endpoints lie in the
finite-cutoff domain. Every auxiliary integral involved
then uses noise in \([-R,t]\).
For the displayed sequence the \(L^1\) errors are summable.
Tonelli's theorem gives almost sure convergence for
each compact interval and admissible exponent.
A countable exhaustion and cofinal sequence of
exponents give the simultaneous version.
\end{proof}

These controls converge to the original fBm first level while
encoding the selected higher increments. This dependence on the
enhancement distinguishes them from unmodified polygonal or
Gaussian mollifier approximations of fBm. The finite-mesh controls
need not preserve the symmetries of the limiting model.

The critical realization also allows deterministic mass-one
profiles \(\nu\) on \((0,1]\) with
\(\|\nu\|_{\mathrm{TV}}\leq M\). Replace
\(Y^{R,\ell}\) by \(Y^{R,\nu}\) in the corollary, with the
same centering. Theorem~\ref{lp:thm:stability} gives
\[
 \|d_\alpha(S_N(b^{\varepsilon,R,\nu}),A^0)\|_{L^p}
       \leq C_M\left(\varepsilon^{\beta-\alpha}
                                  +R^{-H}+M_\sigma(\nu)\right).
\]
Indeed its bound at \(\beta\), with any admissible positive
\(\sigma_\beta\), supplies the required uniform rough-norm
moments because \(M_{\sigma_\beta}(\nu)\leq M\); its bound at
\(\alpha\) gives the stated full range of \(\sigma\).
For signed profiles converging weakly to \(\delta_0\), the
qualitative part of that theorem gives ordinary-signature
convergence in finite \(L^p\) as \(R\to\infty\) and
\(\varepsilon\downarrow0\), even without vanishing absolute
lag moments. Continuity then gives ordinary one-form and ODE
convergence in probability uniformly on compacts.

\section{Dependence of integrals and equations on the enhancement}
\label{sec:dependence}
A top-degree calibration preserves the lower-level controlled domain
and changes the integral through an explicit Young term.
The resulting equation correction is the central-translation principle
\citep{FrizOberhauser2009,BrunedEtAl2019} applied to the models above.
For comparison, Proposition~\ref{prop:closed-forms} characterizes the
smooth one-forms whose integrals are independent of every enhancement.

\begin{proposition}[Top-degree calibration and calculus]
\label{prop:calibration-calculus}
Let \(\bX,\widetilde{\bX}\) be geometric step-\(N\)
\(\alpha\)-H\"older rough paths with the same levels below
\(N\), where \((N+1)\alpha>1\).
Then
\[
 A_{s,t}:=\pi_N\widetilde{\bX}_{s,t}-\pi_N\bX_{s,t}
\]
is an additive degree-\(N\) Lie-valued increment, of
H\"older order \(N\alpha\), and
\(\widetilde{\bX}_{s,t}=\bX_{s,t}\exp(A_{s,t})\).
Their integration-jet domains coincide.
For a common jet \(c\), their integrals satisfy
\begin{equation}\label{eq:top-integral-change}
 \widetilde J_{s,t}-J_{s,t}
      =\sum_{|w|=N}\int_s^t c_w(r)\dd A_r^w ,
\end{equation}
where the right-hand side is a Young integral.
For the lagged critical models, with jet fixed in both,
this becomes
\begin{equation}\label{eq:lag-integral-change}
 J^\ell_{s,t}-J^L_{s,t}
       =-\sum_{|w|=N}\int_s^t
                         c_w(r)(K^{\ell,L}(r))^w\dd r .
\end{equation}
In particular, the polynomial integrals with recentered integrands
\[
 I^\ell_{s,t}
   :=\int_s^t(B_r^1-B_s^1)^{N-1}\dd_{Y^\ell} B_r^2
       =(N-1)!(Y^\ell_{s,t})^{1^{N-1}2}
\]
are adapted to the right endpoint even at negative times
and depend nontrivially on the lag.

For \(V_i\in C_b^{N+2}\), the RDE driven by
\(\widetilde{\bX}\) is equivalently the equation driven by
\(\bX\) with the additional Young term
\(\Phi_V(\dd A)\).
Here \(\Phi_V\) evaluates words in the untruncated free Lie
algebra at the vector fields \(V_i\), and is then restricted
to homogeneous degree \(N\).
For the change from \(Y^L\) to \(Y^\ell\), the additional
ordinary time drift is
\begin{equation}\label{eq:lag-equation-change}
                -\Phi_V(K^{\ell,L}(r))\dd r.
\end{equation}
\end{proposition}
\begin{proof}
Lemma~\ref{lem:central-fibre} gives the additive top-degree
Lie increment \(A\) and the identity
\(\widetilde{\bX}_{s,t}=\bX_{s,t}\exp(A_{s,t})\).
Its \(N\alpha\)-H\"older bound follows by subtracting the
top coordinates of the two rough paths.
The jet relation \eqref{eq:integration-jet} uses only
degrees of the driver below \(N\), so its domain is the
same for both models. For \(|w|=N\), that relation also
gives \(c_w\in C^\alpha\).
Subtract the compensated sums: only
\(\sum_{[u,v]}\sum_{|w|=N}c_w(u)A_{u,v}^w\) remains.
Since \(\alpha+N\alpha>1\), these sums converge to the
Young integral in \eqref{eq:top-integral-change}.
Equation~\eqref{eq:lag-change} proves the negative sign in
\eqref{eq:lag-integral-change}.

The repeated-letter signature identity gives
\[
 \int_s^t(B_r^1-B_s^1)^{N-1}\dd_{\bX}B_r^2
                   =(N-1)!X_{s,t}^{1^{N-1}2}.
\]
It follows either from the geometric shuffle formula or
by applying Chen's identity and the binomial expansion to
the compensated sums. For each fixed \(s\), its jet coefficients depend on
\(B_r-B_s\), \(r\geq s\). They are measurable from the
complete increment past at \(r\), including when \(s<0\).
By \eqref{eq:lag-nondegenerate} and \(L^2\) continuity,
\[
 h^{-1}\int_0^h(K^{\ell,L}(r))^w\dd r
                  \longrightarrow(K^{\ell,L}(0))^w
                  \quad\text{in }L^2 .
 \]
Thus the polynomial integral maps are different.

For the equation, put \(V_i f=Df\,V_i\) and
\(V_w=V_{i_1}\cdots V_{i_m}\) for a word \(w\).
The solution jet has coefficients \(V_w\mathrm{id}(U)\).
Applying \eqref{eq:top-integral-change} to this jet
adds \(\sum_{|w|=N}V_w\mathrm{id}(U)\dd A^w\).
For a Lie-valued \(A\), its word evaluation is a first-order
operator, precisely \(\Phi_V(\dd A)(U)\).
This proves the equation correction and the sign in
\eqref{eq:lag-equation-change}.
One must evaluate in the untruncated free Lie algebra;
general vector fields do not define a Lie homomorphism
from the truncated nilpotent Lie algebra.
\end{proof}

\section{Selection and stability questions}\label{sec:frontier}
The remote-past construction gives a stationary causal calculus with
sharp scaling moments at every \(0<H\leq1/4\).
Existence and symmetry do not by themselves determine a unique
enhancement. The remaining questions concern
normalization beyond fixed lower levels and stability beyond the
specified approximation schemes.

\begin{problem}[Selection and enhancement dependence]
Which additional normalization principles select a stationary causal
model independently of its explicit construction?
At reciprocal indices, logarithmic centering determines the top level
once all lower levels and the common-noise symmetries are fixed
(Corollary~\ref{cor:critical-centering}).
A broader selection principle would also have to constrain changes
in the lower levels. The top-central gauges at nonreciprocal indices
and the positive-lag critical family give concrete tests for such a
principle. Time reversal changes the information direction, so
comparison within a fixed causal class is a separate question.
\end{problem}

\begin{problem}[Stability of the selection]
Determine which other causal kernel approximations or subtractions
converge to the selected enhanced law.
The joint cutoff limit, adapted ordinary-signature realization and
normalized lag-profile stability give three classes for comparison
(Theorems~\ref{zl:thm:joint}, \ref{thm:adapted-BV} and
\ref{lp:thm:stability}).
For schemes selecting another law, identify whether the change occurs
in lower levels or through a central correction, and how it affects
integrals and equations through
Proposition~\ref{prop:calibration-calculus}.
\end{problem}

\appendix
\section{The complete fBm past and the Wiener past}\label{sec:filtration}
For the moving-average model define completed Gaussian sigma-fields
\[
 \mathcal F_t^W=\sigma\{W(f):f\in L^2(\R),\
                                  \supp f\subset(-\infty,t]\},
 \quad
 \mathcal F_t^B=\sigma\{B_b^H-B_a^H:a<b\leq t\}.
\]
These are increment-past filtrations on the whole real line.
For \(t<0\), the anchored process values \(B_r\), \(r\leq t\),
with \(B_0=0\), must not replace the increment family in the
second definition.
Equality of these complete increment filtrations is established in
\citet[Theorem~4.1 and its following discussion]{Picard2011}.
We include a proof in the kernel convention used here.

\begin{proposition}\label{prop:complete-past}
For \(0<H<1/2\), \(\mathcal F_t^B=\mathcal F_t^W\).
\end{proposition}
\begin{proof}
Put \(\delta=1/2-H\in(0,1/2)\).
Each kernel \(\kappa_b-\kappa_a\), \(a<b\leq t\), is supported in
\((-\infty,t]\), so one inclusion is immediate.
To prove the other, suppose that \(f\in L^2((-\infty,t])\),
extended by zero, is orthogonal to all these increment kernels.
For \(f\in L^2(\R)\), define the causal fractional derivative
\(D_+^\delta f\in\mathcal S'(\R)\) by its Fourier representative
\((i\xi)^\delta\widehat f(\xi)\). Write \(D_-^\delta\) for its
adjoint on the \(L^2\) domain; smooth compactly supported functions
and the interval indicators below belong to that domain.
The Fourier transform here has negative exponent and the principal
branches are used.
In \(L^2(\R)\),
\[
 D_-^\delta\ind_{(a,b)}
 =\frac{(b-u)_+^{-\delta}-(a-u)_+^{-\delta}}
        {\Gamma(1-\delta)}
 =\frac{\kappa_b-\kappa_a}{\mathfrak c_H\Gamma(1-\delta)}.
\]
For a smooth test function \(\varphi\) compactly supported before
\(t\), choose \(a<\inf\supp\varphi\). Its representation as a
Bochner integral of \(\ind_{(a,b)}\), weighted by
\(-\varphi'(b)\), shows that
\(\langle f,D_-^\delta\varphi\rangle=0\).
The Bochner integral is valid because the displayed increment
kernels have finite \(L^2\) norm of order \(|b-a|^H\).
More precisely the indicator-valued representation holds
in \(L^2\); its fractional derivative is a compactly supported
Bochner integral in the integration parameter, with values in
\(L^2(\R)\). Its distributional value is
\(D_-^\delta\varphi\), justifying the adjoint pairing.
Consequently the distribution \(g=D_+^\delta f\) vanishes on
\((-\infty,t)\), so \(\supp g\subset[t,\infty)\).

For \(\varepsilon>0\), convolve with the causal resolvent kernel
\[
 K_{\delta,\varepsilon}(r)
 =\ind_{\{r>0\}}\frac{r^{\delta-1}e^{-\varepsilon r}}{\Gamma(\delta)}.
\]
Its multiplier is \((\varepsilon+i\xi)^{-\delta}\).
Convolution with this rapidly decreasing causal distribution
preserves support in \([t,\infty)\). To justify this for the
possibly distributional \(g\), multiply the kernel by a smooth
cutoff equal to one on \([-L,L]\) and supported on \([-2L,2L]\).
The truncated convolution is defined and has support in
\([t,\infty)\). The omitted tail is smooth and tends to zero
in every Schwartz seminorm as \(L\to\infty\), by exponential
decay. Its convolution with a Schwartz test function tends to
zero in the Schwartz topology, so the convolutions converge
in tempered distributions. Equivalently, the full multiplier
\((\varepsilon+i\xi)^{-\delta}\) is smooth with all derivatives
polynomially bounded, defining the same convolution on
tempered distributions.
The resulting distribution \(h_\varepsilon\) is an \(L^2\) function,
with
\[
 \widehat h_\varepsilon(\xi)
 =\frac{(i\xi)^\delta}{(\varepsilon+i\xi)^\delta}\widehat f(\xi).
\]
The multiplier has modulus at most one and converges to one
almost everywhere. Plancherel and dominated convergence imply
\(h_\varepsilon\to f\) in \(L^2\).
Hence \(f\) is also supported in \([t,\infty)\).
Its original support was in \((-\infty,t]\), so \(f=0\).
The increment kernels therefore span the whole Wiener-past
Hilbert space. Equality of the completed Gaussian sigma-fields
follows.
\end{proof}

\begin{corollary}
The enhancements constructed in
Theorems~\ref{thm:noncritical-causal} and
\ref{thm:critical-causal} are adapted to the complete fBm
increment-past filtration.
The stationary critical enhancements of
Theorem~\ref{thm:moving-lag}, their zero-lag limit in
Theorem~\ref{zl:thm:zero}, and its centered representative in
Corollary~\ref{cor:critical-centering} have the same adaptedness.
For fBm specified only on a finite positive interval, its coupling
to a complete negative-time history is additional model data.
\end{corollary}

\section{Deterministic lag profiles}\label{app:lag-profiles}
\label{sec:lag-profiles}
We continue at a reciprocal Hurst index, with the notation of
Section~\ref{sec:zero-lag}. A deterministic signed measure of lags
can replace the single lag in the critical calibration.
For profiles accumulating at zero, averaging acts on the integrated
central correction, rather than its velocity.
Fix \(I=[-T,T]\), \(T\geq1\), and an admissible \(\alpha<H\).
Put \(a=N\alpha\) and
\(\mathcal B_a=c_0^a(I;\mathfrak L_N)\), with anchor zero.
For sufficiently large \(R\) define
\[
 C_R^\lambda(t)=\int_0^t z_R^\lambda(r)\dd r,
                         \qquad 0<\lambda\leq1.
\]

\begin{lemma}[Continuity in the lag]\label{lp:lem:continuity}
For every \(1\leq p<\infty\), the map
\(\lambda\mapsto C_R^\lambda\) extends continuously to
\([0,1]\) with values in \(L^p(\Omega;\mathcal B_a)\).
These extensions are consistent across compact intervals,
moment orders and admissible exponents.
\end{lemma}
\begin{proof}
On a positive compact lag range \([\delta,1]\), change
the upper auxiliary endpoint at fixed terminal parameter.
The finite-band concatenation argument of
Lemma~\ref{lem:lag-band}, including one inserted terminal
derivative, gives
\[
 \|z_R^\lambda(r)-z_R^{\lambda'}(r)\|_{L^q}
       \leq C_{R,I,\delta,q}|\lambda-\lambda'|^{1/2},
                 \qquad r\in I.
\]
All weights on the added band are bounded by the positive gap;
its positive-degree factors have an \(L^2\) weight of order
\(|\lambda-\lambda'|^{1/2}\). Polynomial inversion and
multiplication preserve the estimate. Time integration and
Lemma~\ref{zl:lem:GRR}, with time exponent one, give continuity
in \(L^p(\mathcal B_a)\). Each positive-lag primitive is
\(C^1\), so it belongs to the little-H\"older space.

All finite-cutoff lag models have the same lower levels.
Centrality therefore gives the increment identity
\[
 \delta(C_R^\lambda-C_R^\mu)_{s,t}
       =\pi_N(Y^{R,\mu}_{s,t}-Y^{R,\lambda}_{s,t}),
 \qquad \delta C_{s,t}=C(t)-C(s).
\]
Corollary~\ref{zl:cor:reverse} makes these anchored differences
Cauchy in \(L^p(\mathcal B_a)\) as \(\lambda\downarrow0\).
This gives the continuous extension, which is bounded on the
compact parameter interval. Limits in the common
\(L^1(C(I))\) space prove consistency. In particular the lag
map is strongly measurable.
\end{proof}

For a deterministic finite real signed measure \(\nu\) on
\((0,1]\), write
\[
 \|\nu\|_{\mathrm{TV}}=|\nu|((0,1]),\qquad
 M_\sigma(\nu)=\int_{(0,1]}\lambda^\sigma|\nu|(\dd\lambda).
\]
The preceding lemma defines a Bochner integral and a model by
\begin{equation}\label{lp:eq:definition}
 C_R^\nu=\int_{(0,1]}C_R^\lambda\,\nu(\dd\lambda),
 \qquad
 Y^{R,\nu}_{s,t}=X^R_{s,t}\exp(-\delta C_R^\nu{}_{s,t}).
\end{equation}
The Bochner integrals agree in all the indicated spaces.

\begin{theorem}[Normalized lag-profile stability]\label{lp:thm:stability}
Suppose \(\nu((0,1])=1\). The model
\eqref{lp:eq:definition} is geometric, retains every level below
\(N\) of \(X^R\), is causal for innovation noise on \([-R,t]\)
at the right endpoint \(t\), and commutes with linear spatial maps.
For every \(1\leq p<\infty\) and
\(0<\sigma<\min\{H/2,1-N\alpha\}\),
\begin{equation}\label{lp:eq:rate}
 \|d_{\alpha;I}(Y^{R,\nu},Y^0)\|_{L^p}
 \leq C\left(\|\nu\|_{\mathrm{TV}}R^{-H}+M_\sigma(\nu)\right),
\end{equation}
with \(C\) independent of the deterministic profile and cutoff.
For each fixed \(\nu\), there is a stationary causal limit
\(Y^\nu\), natural under linear maps, with
\[
 \|d_\alpha(Y^{R,\nu},Y^\nu)\|_{L^p}
       \leq C\|\nu\|_{\mathrm{TV}}R^{-H},\qquad
 \|d_\alpha(Y^\nu,Y^0)\|_{L^p}\leq C M_\sigma(\nu).
\]
If \(R_n\to\infty\), the profiles \(\nu_n\) have mass one
and uniformly bounded total variation, and
\(\nu_n\) converge weakly to \(\delta_0\) as signed measures
on \([0,1]\), then \(Y^{R_n,\nu_n}\to Y^0\) in every finite
\(L^p\) coordinate \(\alpha\)-rough topology.
\end{theorem}
\begin{proof}
The averaged path is additive and Lie-valued in the central
degree \(N\). Its exponential preserves Chen and shuffle.
Working first at a slightly larger exponent below \(H\)
proves geometricity at \(\alpha\). Mass one gives the exact
affine-fiber identity
\begin{equation}\label{lp:eq:affine}
 \pi_jY^{R,\nu}_{s,t}
       =\int\pi_jY^{R,\lambda}_{s,t}\,\nu(\dd\lambda),
                         \qquad 1\leq j\leq N.
\end{equation}
Indeed the lower coordinates are independent of \(\lambda\),
and the top coordinate is \(\pi_NX^R_{s,t}-\delta C_R^\lambda{}_{s,t}\).
This uses the central fiber, rather than an averaging operation
on arbitrary group elements. Subtract \(\pi_jY^0_{s,t}\), take
the coordinate rough seminorm, and apply Minkowski and
Theorem~\ref{zl:thm:joint}. This proves \eqref{lp:eq:rate}.
The inequality can first be taken on rational endpoints,
then extended by continuity.

For negative \(t\), use oriented increments
\(Y_{0,t}=(Y_{t,0})^{-1}\). At infinite cutoff put
\[
 D^\lambda(t)=\pi_N(Y^\lambda_{0,t}-Y^0_{0,t}),\qquad
 Y^\nu_{s,t}=Y^0_{s,t}
      \exp\left(\delta\!\left[\int D^\lambda\nu(\dd\lambda)\right]_{s,t}\right).
\]
The central path \(D^\lambda\) is continuous into
\(L^p(\mathcal B_a)\) on \([0,1]\), with \(D^0=0\).
For positive lags this follows by uniformly approximating it
with \(\pi_N(Y^{R,\lambda}_{0,\cdot}-Y^{R,1}_{0,\cdot})\)
and adding the fixed reference difference;
\eqref{zl:eq:uniform-R} bounds the approximation error
uniformly in the lag, with the supremum outside the \(L^p\) norm. At zero use
Theorem~\ref{zl:thm:zero}. The affine identity and the same
uniform-cutoff estimate give the two bounds for \(Y^\nu\).

For the weak-concentration assertion, approximate the continuous
Banach-valued function \(D:[0,1]\to L^p(\mathcal B_a)\)
uniformly by its piecewise linear interpolants. These are
finite sums of fixed Banach vectors times scalar continuous
hat functions. Weak convergence handles each finite sum;
the uniform total-variation bound controls the approximation
error. Hence \(\int D^\lambda\nu_n(\dd\lambda)\to D^0=0\)
in this Banach space. Add the uniform-cutoff error
\(C\|\nu_n\|_{\mathrm{TV}}R_n^{-H}\).

For causality, each increment \(\delta C_R^\lambda{}_{s,t}\)
is measurable from noise on \([-R,t]\). This completed
measurable subspace of \(L^p\) is closed, so its Bochner
average has the same property. The assertion concerns
increments; a primitive anchored at zero need not be adapted
at negative times. Naturality commutes with scalar integration.
The finite-cutoff shift identity passes to the limit for a fixed
absolute-lag profile, proving stationarity of \(Y^\nu\).
\end{proof}

\begin{remark}[A common bound on the lags]\label{lp:rem:unbounded}
The common support restriction in the profile theorem cannot be
discarded merely because the profiles concentrate weakly at zero.
For example, the probability measures
\[
 \nu_n=(1-1/n)\delta_{1/n}+(1/n)\delta_{\exp(n^4)}
\]
converge weakly to \(\delta_0\) on \([0,\infty)\).
Each finite mixture still defines an infinite-cutoff model by the
central affine-fiber formula. However, for \(w=1^{N-1}2\),
\eqref{eq:lag-UV} at \(h=1\) gives
\(\|(Y^L_{0,1})^w\|_2\geq c\sqrt{\log(L/2)}\) for large \(L\).
The triangle inequality, the finite moments of \(Y^0\), and the
zero-lag convergence at \(1/n\) therefore give
\[
 \|\pi_N(Y^{\nu_n}_{0,1}-Y^0_{0,1})\|_2\geq cn-C.
\]
Thus even probability profiles with escaping lag support need not
converge in the asserted \(L^2\) topology. A single fixed bounded
lag interval is essential to the stated weak-profile conclusion.
\end{remark}

For relative profiles obtained by pushing a fixed mass-one
measure \(\eta\) on \([b_-,b_+]\Subset(0,\infty)\) forward
under \(\theta\mapsto\ell\theta\), with
\(0<\ell\leq1/b_+\), the lag term is
\(\ell^\sigma\int\theta^\sigma|\eta|(\dd\theta)\).
On this positive compact range the primitive definition agrees
with averaging the velocities by ordinary Fubini.
All comparison bounds survive a common deterministic critical
centering. The weak signed-profile conclusion is stronger than
the sufficient condition \(M_\sigma(\nu_n)\to0\): for fixed
\(0<q<1\), the profiles
\(\delta_{1/n}+\delta_q-\delta_{q+1/n}\), for large \(n\),
converge weakly to \(\delta_0\), although their absolute
\(\sigma\)-moments tend to \(2q^\sigma\).

\subsection{The normalization of a lag profile}
\label{sec:lag-mass}
Any fixed nonunit mass prevents remote-past tightness
within the bounded-lag class.

\begin{proposition}[Nonunit mass]\label{lp:prop:nonunit}
Let \(d\geq2\), fix \(B<\infty\) and \(m\ne1\), and let
\(\nu_R\) be arbitrary deterministic finite signed measures on
\((0,B]\) of total mass \(m\). Their total variations need
not be bounded. Define \(Y^{R,\nu_R}\) by primitive averaging.
For \(w=1^{N-1}2\), constants independent of the profiles give
\[
 \E\bigl|(Y^{R,\nu_R}_{0,1})^w\bigr|^2
                         \geq c_m\log R-C_m .
\]
These coordinates are not tight, even after arbitrary
deterministic central centerings depending on the cutoff.
\end{proposition}
\begin{proof}
For fixed \(B\), the continuity lemma and primitive construction
apply on \((0,B]\) after enlarging the finite-cutoff domain.
Project on the highest colored chaos with \(N-1\) variables
of color \(1\), denoted \(u_j\), and one of color \(2\),
denoted \(v=-x\). Restrict to
\[
 x_0\leq x\leq R,\qquad -x/2<u_j<-x/4,
                         \qquad x_0>4(B+1).
\]
All variables precede every endpoint \(r-\lambda\), for
\(0\leq r\leq1\) and \(\lambda\leq B\). The color-\(2\)
variable is the unique minimum. Thus only the last valley
in \eqref{eq:stationary-valley} contributes. After
symmetrizing the color-\(1\) variables, the restricted kernel is
\[
 \frac{(-1)^{N-1}}{(N-1)!}
 \left[(\kappa_1-\kappa_0)(v)\prod_j\kappa_0(u_j)
       -m\int_0^1\kappa'_r(v)\prod_j\kappa_r(u_j)\dd r\right].
\]
The derivative kernel on this region is independent of the
auxiliary endpoint, so the profile enters only through its
mass. Kernel restriction is bounded on the highest chaos and
commutes with the primitive Bochner integral, including for
profiles accumulating at zero.
Uniformly on this region the bracket equals
\[
 (1-m)L+O((1+|m|)x^{-1}|L|),\qquad
                L=\kappa'_0(-x)\prod_j\kappa_0(u_j).
\]
Increase \(x_0\) so the remainder is at most
\(|1-m||L|/2\). Colored Wiener isometry bounds the squared
norm below by a positive constant times
\[
 (1-m)^2\int_{x_0}^R |\kappa'_0(-x)|^2
    \left(\int_{-x/2}^{-x/4}\kappa_0(u)^2\dd u\right)^{N-1}\dd x
               =c(1-m)^2\int_{x_0}^R x^{2NH-3}\dd x.
\]
Since \(NH=1\), this is logarithmic. Other valleys and
lower-chaos contractions cannot cancel this restricted norm.

The full coordinate \(Q_R\), after any deterministic centering,
still belongs to chaoses of orders at most \(N\), and its
highest-chaos lower bound is unchanged. Hypercontractivity gives
\(\|Q_R\|_4\leq3^{N/2}\|Q_R\|_2\). Paley--Zygmund applied
to \(Q_R^2\) gives
\(\Prob\{|Q_R|\geq\|Q_R\|_2/\sqrt2\}\geq(4\,3^{2N})^{-1}\),
which precludes tightness. This argument does not concern
noise-dependent additional counterterms or masses tending to one.
\end{proof}

\section{Receding finite origins at a reciprocal index}
\label{sec:receding-origin}
The original finite-origin lift cannot simply be made
stationary by sending its initial time to negative infinity
at the critical degree. The following calculation uses the
canonical kernel in \citet[equation~(5)]{NualartTindel2011},
which is different from the unanchored kernel \(\kappa\).
Write \(\delta=1/2-H\). In its positive-time variables it is
\begin{equation}\label{eq:NT-canonical-kernel}
 K_H(t,u)=c_H\left[
  (u/t)^\delta(t-u)^{-\delta}
  +\delta u^\delta\int_u^t
          v^{-\delta-1}(v-u)^{-\delta}\dd v\right],
                                  \qquad 0<u<t,
\end{equation}
and is zero elsewhere. The constant \(c_H>0\) here is the
canonical-kernel normalization of \citet{NualartTindel2011}.
Put \(K_R^{\mathrm{NT}}(t,u)=K_H(R+t,R+u)\), and let
\(X^{R,\mathrm{NT}}\) be the corresponding valley construction
(6)--(7) with this translated kernel.

\begin{proposition}[Failure of uncorrected receding origins]
\label{prop:NT-receding}
Let \(H=1/N\), \(N\geq4\), and \(w=1^{N-1}2\).
For all sufficiently large \(R\),
\begin{equation}\label{eq:NT-receding-variance}
 \operatorname{Var}(X^{R,\mathrm{NT},w}_{0,1})
                                      \geq c\log R-C.
\end{equation}
For any deterministic real numbers \(b_R\), the coordinates
\(X^{R,\mathrm{NT},w}_{0,1}-b_R\) are not tight.
Their first-level increment law on \([0,1]\) is the standard
fBm increment law for every \(R\).
\end{proposition}
\begin{proof}
Differentiating \eqref{eq:NT-canonical-kernel}, the derivative
of \((u/t)^\delta\) cancels the upper-limit contribution
from the integral. Thus
\[
 \partial_tK_H(t,u)
       =-c_H\delta(u/t)^\delta(t-u)^{-\delta-1}.
\]
In the highest colored chaos, let the unique color-\(2\)
variable be \(v=-x\), and restrict the color-\(1\) variables to
\[
 x_0\leq x\leq R/2,\qquad
             -x/2<u_j<-x/4,\quad 1\leq j\leq N-1,
\]
with fixed \(x_0\geq4\). Only the last valley contributes
on this region, since \(v\) is the unique minimum.
After symmetrization in color \(1\), its kernel is
\[
 \frac{(-1)^{N-1}}{(N-1)!}
 [K_H(R+1,R-x)-K_H(R,R-x)]
                   \prod_{j=1}^{N-1}K_H(R,R+u_j).
\]
The exact derivative and \((R-x)/(R+r)\geq1/3\),
\(0\leq r\leq1\), give a lower bound
\(c x^{H-3/2}\) for the absolute difference in brackets.
The positive first term in \eqref{eq:NT-canonical-kernel}
gives \(K_H(R,R+u_j)\geq c(-u_j)^{H-1/2}\).
Colored Wiener isometry and orthogonality to lower
Stratonovich contractions therefore give
\begin{align*}
 \operatorname{Var}(X^{R,\mathrm{NT},w}_{0,1})
 &\geq \frac{c}{(N-1)!}\int_{x_0}^{R/2}x^{2H-3}
       \left(\int_{-x/2}^{-x/4}(-u)^{2H-1}\dd u\right)^{N-1}
                                                    \dd x\\
 &=c'\int_{x_0}^{R/2}x^{2NH-3}\dd x
   =c'\log R-C'.
\end{align*}
For \(Q_R=X^{R,\mathrm{NT},w}_{0,1}-b_R\), finite-chaos
hypercontractivity gives
\(\|Q_R\|_4\leq3^{N/2}\|Q_R\|_2\), while
\(\|Q_R\|_2^2\geq c\log R-C\).
Paley--Zygmund applied to \(Q_R^2\) yields
\[
 \Prob\{|Q_R|\geq\|Q_R\|_2/\sqrt2\}
                              \geq (4\cdot3^{2N})^{-1}.
\]
The threshold diverges, proving non-tightness uniformly
over deterministic centerings. Finally, first-level
stationarity follows from the covariance of the fBm
represented by \(K_H\).
\end{proof}

\begin{proposition}[The finite-origin reciprocal endpoint]
\label{prop:NT-endpoint}
Let \(H=1/N\), \(N\geq4\), and let \(X^{\mathrm{NT}}\) denote
the finite-origin construction (6)--(7) of
\citet{NualartTindel2011}, with kernel
\eqref{eq:NT-canonical-kernel}. For \(1\leq q<\infty\), \(1\leq m\leq N\),
\(0\leq s<t\leq T\), and \(h=t-s\),
\begin{equation}\label{eq:NT-endpoint-upper}
 \|\pi_m X^{\mathrm{NT}}_{s,t}\|_{L^q}
 \leq C_{H,m,q,T}
 \begin{cases}
 h^{mH},&mH<1,\\
 h\{1+\log(T/h)\},&mH=1.
 \end{cases}
\end{equation}
For every fixed \(s_0\in(0,T)\), the word \(w=2\,1^{N-1}\)
and sufficiently small \(h>0\) satisfy
\begin{equation}\label{eq:NT-endpoint-lower}
 \operatorname{Var}(X^{\mathrm{NT},w}_{s_0,s_0+h})
                  \geq c h^2\log(r_0/h),
\end{equation}
where \(0<r_0<\min(s_0/4,T-s_0)\) is fixed.
Consequently a uniform \(Ch^{2NH}=Ch^2\) second-moment bound
does not hold for this finite-origin construction at the endpoint.
\end{proposition}
\begin{proof}
The case \(m=1\) follows from Wiener isometry and Gaussian moments.
For \(m\geq2\), use the finite-domain valley identity and
Lemma~\ref{lem:minimum-weight}, which applies to arbitrary repeated
colors. Put \(\beta=H-1/2\),
\[
 \psi(v)=|K_H(t,v)-K_H(s,v)|,\quad
 \phi(v)=K_H(s,v)+K_H(t,v),\quad
 M(v)=\int_v^t\phi(u)^2\dd u.
\]
The canonical kernel bound
\(K_H(t,u)\leq C\{(t-u)^\beta+u^\beta\}\) on \(0<u<t\)
follows from \eqref{eq:NT-canonical-kernel} by scaling the integral
with \(v=u(1+y)\). It gives \(M(v)\leq C(t-v)^{2H}\): in the origin term use
\(u^{2H-1}\leq(u-v)^{2H-1}\) for \(u>v\).
When \(s\leq2h\), mixed-weight signature estimates and
\(\|K_H(t,\cdot)\|_2=t^H\) bound the entire inverse-signature
product by \(Ch^{mH}\). Assume therefore \(s>2h\).
On \(s-h<v<t\), the two terms of the minimum estimate are
bounded by \(Ch^{mH}\), since \(\psi\leq\phi\),
\(M\leq Ch^{2H}\), and \(\int_{s-h}^t\phi^2\leq Ch^{2H}\).

On \(0<v<s-h\), put \(r=s-v\). The exact kernel derivative
above yields \(\psi(v)\leq Ch r^{H-3/2}\), while
\(\phi(v)\leq C(r^\beta+v^\beta)\) and \(M(v)\leq Cr^{2H}\).
The squared It\^o contribution is at most
\[
 Ch^2\int_h^s r^{2mH-3}\dd r.
\]
The contraction contribution is at most
\[
 Ch\int_h^s r^{mH-2}\dd r
 +Ch\int_0^{s-h}(s-v)^{(m-1)H-3/2}v^{H-1/2}\dd v.
\]
The last term is bounded by \(Ch^{mH}\) for \(mH\leq1\).
Indeed split at \(v=s/2\); if the remaining power of \(s-v\)
is \(-1\), use
\(\sup_{z\geq1}z^{-(1/2-H)}\log z<\infty\).
The other two integrals give \(Ch^{mH}\) when \(mH<1\), and
at most \(Ch\sqrt{1+\log(T/h)}\) and
\(Ch\{1+\log(T/h)\}\), respectively, when \(mH=1\).
Hypercontractivity gives
every finite \(q\).

For the lower bound project onto the colored chaos with one
color-\(2\) variable \(v=s_0-r\), \(h<r<r_0\), and
\(N-1\) color-\(1\) variables in
\((s_0-r/2,s_0-r/4)\). The color-\(2\) variable is the
unique minimum, so only the first valley contributes. Symmetrizing
color \(1\) gives the factor \(1/(N-1)!\). Colored isometry
therefore yields
\begin{align*}
 \operatorname{Var}(X^{\mathrm{NT},w}_{s_0,s_0+h})
 &\geq \frac1{(N-1)!}\int_h^{r_0}
 |K_H(s_0+h,s_0-r)-K_H(s_0,s_0-r)|^2\\
 &\hspace{13mm}\times
 \left(\int_{s_0-r/2}^{s_0-r/4}K_H(s_0+h,u)^2\dd u\right)^{N-1}
 \dd r\\
 &\geq c h^2\int_h^{r_0}r^{2NH-3}\dd r
 =c h^2\log(r_0/h).
\end{align*}
Here the exact derivative gives the bracket a lower bound
\(ch r^{H-3/2}\); the positive first kernel term gives the
inner integral a lower bound \(cr^{2H}\).
Orthogonality to all lower Stratonovich chaoses justifies using
this projection as a variance lower bound.
\end{proof}

\Needspace{9\baselineskip}
\smallskip\noindent\emph{The case \(H=1/4\), \(w=2111\).}
Fix \(0<h<r_0<\min(s_0/4,T-s_0)\), put \(t=s_0+h\),
\(I_r=(s_0-r/2,s_0-r/4)\), and write \(c=c_{1/4}\).
On \(v=s_0-r\), \(h<r<r_0\), \(u_j\in I_r\),
the fourth-chaos kernel, symmetric in the three color-\(1\)
variables, is exactly
\[
 F_h(v;u_1,u_2,u_3)
 =\frac{K_{1/4}(t,v)-K_{1/4}(s_0,v)}{3!}
       \prod_{j=1}^{3}K_{1/4}(t,u_j).
\]
Only the first valley contributes, with prefactor \(+1\);
the kernel difference is negative. The ratios
\(v/(s_0+a)\) and \(u_j/t\) exceed \(1/2\) for
\(0\leq a\leq h\), while \(r+a\leq2r\) and
\(t-u_j\leq2r\). The derivative above and the positive
first kernel term consequently give
\[
 |K_{1/4}(t,v)-K_{1/4}(s_0,v)|
       \geq c2^{-7/2}hr^{-5/4},\qquad
 \int_{I_r}K_{1/4}(t,u)^2\dd u
       \geq\frac{c^2}{8}r^{1/2}.
\]
The colored isometry has factor \(3!\), so orthogonality
to lower chaoses yields the explicit bound
\[
 \operatorname{Var}(X^{\mathrm{NT},2111}_{s_0,s_0+h})
 \geq\frac{c^8}{3!\,2^{16}}h^2
             \int_h^{r_0}\frac{\dd r}{r}
 =\frac{c^8}{3!\,2^{16}}h^2\log(r_0/h).
\]

The sharp endpoint estimate in
\citet[Proposition~4.1]{NualartTindel2011}
cannot be used at \(NH=1\). The lower bound uses only the first
valley, present in formulas (7) and (28) of that work.
The logarithmic bound
\eqref{eq:NT-endpoint-upper} still gives every group H\"older
exponent strictly below \(H\), so the finite-origin rough-path
construction remains valid.

Proposition~\ref{prop:NT-receding} concerns translation of the
finite-origin prescription. The critical construction in
Section~\ref{sec:critical-main}
uses random calibrations of the separate auxiliary-kernel
family. Algebraic character reconstruction in
\citet{UnterbergerGeneric2010} and the connection explained
in \citet[Section~5.2]{NualartTindel2011} accommodate
Volterra choices. The remote-cutoff and zero-lag estimates
establish convergence for the stationary causal selection.

\section{Fixed central slope calibrations and their obstructions}
\label{sec:slope-obstructions}
We retain the reciprocal-index cutoff family and fixed-reference
calibration of Section~\ref{sec:critical-construction}.

\begin{proposition}[Freedom within central slope subtractions]
\label{prop:central-scheme}
Fix the cutoff family in Theorem~\ref{thm:critical-causal}.
Let \(C_R\) be any random degree-\(N\) Lie element.
Then \(X^R_{s,t}\exp(-(t-s)C_R)\) converges in probability in the
same rough-path topologies if and only if \(C_R-Z_R\) converges in
probability in the finite-dimensional space of degree-\(N\) Lie
elements. If the latter limit is \(D\), the resulting enhancement is
\[
 \widehat X_{s,t}\exp(-(t-s)D).
\]
This classifies time-linear central subtractions of this cutoff
family, not all enhancements of fBm.
\end{proposition}
\begin{proof}
The difference of the highest coordinates of the two renormalized
families is exactly \(-(t-s)(C_R-Z_R)\).
At \((s,t)=(0,1)\), convergence of both families forces convergence
of \(C_R-Z_R\). Conversely a convergent slope gives convergence
in the degree-\(N\) H\"older seminorm because
\(\sup_{0<h\leq T}h^{1-N\alpha}<\infty\).
Lower coordinates are identical.
\end{proof}

\begin{proposition}[Scaling obstruction for central slope calibrations]
\label{prop:critical-scaling-obstruction}
Let \(d\geq2\), \(H=1/N\), and \(N\geq4\). Write
\(w=1^{N-1}2\). For every random degree-\(N\) Lie element
\(D\) with square-integrable coordinates, set
\[
 Y_{s,t}=\widehat X_{s,t}\exp(-(t-s)D).
\]
There are constants \(c>0,C<\infty\), independent of \(h,D\),
such that for \(0<h<1/4\),
\begin{equation}\label{eq:critical-scaling-lower}
 h^{-1}\|Y_{0,h}^w\|_2
 \geq c\sqrt{\log(1/h)}-C-\|D^w\|_2 .
\end{equation}
In particular no such calibration has exact degree-\(NH=1\)
self-similarity in law with a square-integrable highest level.
This applies to every convergent time-linear central
subtraction in Proposition~\ref{prop:central-scheme}
whose limiting model is square-integrable.
\end{proposition}
\begin{proof}
Let \(L_h=(J_0^{-1}J_h)^w\), where \(J\) is the local factor
over \([-1,t]\). Project onto the colored Wiener chaos of
degrees \(N-1\) in color \(1\) and \(1\) in color \(2\).
Denote the color-\(2\) variable by \(v\).
On
\[
 -\tfrac12<v<-h,\qquad v<u_\ell<0
                  \quad(1\leq\ell<N),
\]
only the last valley in \eqref{eq:stationary-valley}
contributes: every other valley has a color-\(1\) minimum.
After symmetrizing the color-\(1\) variables, this kernel is
\[
 \frac{(-1)^{N-1}}{(N-1)!}
      (\kappa_h-\kappa_0)(v)
                    \prod_{\ell=1}^{N-1}\kappa_0(u_\ell).
\]
Colored Wiener isometry therefore yields
\[
 \|L_h\|_2^2\geq
 \frac1{(N-1)!}\int_{-1/2}^{-h}
    |\kappa_h(v)-\kappa_0(v)|^2
       \left(\int_v^0\kappa_0(u)^2\dd u\right)^{N-1}\dd v .
\]
For \(h\leq r\leq1/2\), the mean-value theorem gives
\[
 |\kappa_h(-r)-\kappa_0(-r)|
       \geq c\,h r^{H-3/2},\qquad
 \int_{-r}^0\kappa_0(u)^2\dd u
       =\frac{\mathfrak c_H^2}{2H}r^{2H}.
\]
Since \(NH=1\), the last integral lower bound is
\(c h^2\int_h^{1/2}r^{-1}\dd r\).
Lower-chaos Stratonovich contractions cannot cancel this
orthogonal highest-chaos contribution.

The local--remote product estimate in the proof of
Theorem~\ref{thm:critical-causal} gives
\(\|\widehat X_{0,h}^w-L_h\|_2\leq Ch\).
Centrality gives \(Y_{0,h}^w=\widehat X_{0,h}^w-hD^w\).
The triangle inequality proves
\eqref{eq:critical-scaling-lower}, without any independence
assumption on \(D\).
Exact scaling would imply
\(Y_{0,h}^w\stackrel{\mathrm{law}}=hY_{0,1}^w\),
contradicting the lower bound.
Finally, if a limiting central calibration has
square-integrable coordinates, its difference at \((0,1)\)
from \(\widehat X\) identifies its slope \(D\) as an
\(L^2\) random variable.
\end{proof}

\begin{proposition}[Stationarity obstruction for central slope calibrations]
\label{prop:critical-stationarity-obstruction}
Under the hypotheses of
Proposition~\ref{prop:critical-scaling-obstruction},
there is \(c>0\), independent of \(s,D\), such that
\[
 \|Y_{s,s+1}^w\|_2
       \geq c\sqrt{\log s}-\|D^w\|_2,\qquad s\geq2.
\]
Consequently no square-integrable time-linear central
calibration of this cutoff family has stationary increments
on the positive half-line.
\end{proposition}
\begin{proof}
Let \(Q_+\) be orthogonal projection onto the colored
\((N-1,1)\) Wiener chaos generated only by white noise on
\((0,\infty)\). In kernel form this first takes the
colored chaos projection and then restricts all its
variables to \((0,\infty)\).
Since \(Z_R\) is measurable from noise before \(-1\),
independence of the two noise subspaces gives \(Q_+Z_R^w=0\).
At degree \(N\), central multiplication subtracts \(Z_R^w\)
exactly for an interval of length one. Thus
\(Q_+(\widehat X_{s,s+1}^w)\) is the same positive-noise
kernel as \(Q_+((X^R_{s,s+1})^w)\).
This identity passes to the limit by the pointwise \(L^2\)
convergence proved in Theorem~\ref{thm:critical-causal}.

Restrict that kernel further to
\[
 0<v<s-1,\qquad v<u_\ell<s\quad(1\leq\ell<N),
\]
with \(v\) the color-\(2\) variable.
Only the last valley contributes. The same symmetrization
and isometry as above give
\[
 \|\widehat X_{s,s+1}^w\|_2^2
 \geq\frac1{(N-1)!}\int_0^{s-1}
 |\kappa_{s+1}(v)-\kappa_s(v)|^2
       \left(\int_v^s\kappa_s(u)^2\dd u\right)^{N-1}\dd v.
\]
Put \(r=s-v\). For \(r\geq1\) the mean-value estimate is
\(|\kappa_{s+1}(v)-\kappa_s(v)|\geq c r^{H-3/2}\);
the inner integral equals
\(\mathfrak c_H^2r^{2H}/(2H)\).
The right-hand side is therefore at least
\(c\int_1^s r^{-1}\dd r=c\log s\).
No constant from a remote-path estimate on \([0,s+1]\)
is used in this lower bound.
Finally \(Y_{s,s+1}^w=\widehat X_{s,s+1}^w-D^w\);
the triangle inequality proves the assertion.
Stationarity would make this \(L^2\) norm independent of
\(s\), which is impossible.
\end{proof}

These obstructions are specific to time-linear central calibrations
of this cutoff family. Time-dependent calibrations, as in
Section~\ref{sec:moving-lag}, and lifts constructed by other methods
fall outside their scope.

\section{Ergodic separation of stationary enhanced laws}
\label{sec:ergodic-laws}

Fix \(d,N\). Let \(\mathcal E\) and \(\mathcal E_+\) be the
spaces of continuous \(G^N(\R^d)\)-valued Chen increment families
on the whole real line and the positive half-line, respectively,
with their local uniform topologies. Their Borel sigma-fields
are generated by increment evaluations at rational endpoints.
Write
\[
 (T_1\mathbf x)_{s,t}=\mathbf x_{s+1,t+1}.
\]

\begin{lemma}[White-noise factors and ergodic separation]
\label{lem:ergodic-singularity}
The unit shift \(\theta_1\) of two-sided white noise is mixing.
If a measurable map \(\Phi\) from that noise into \(\mathcal E\)
satisfies \(\Phi(\theta_1W)=T_1\Phi(W)\) almost surely, then its
whole-line law and its positive-half-line restriction are stationary
and mixing, hence ergodic, under \(T_1\).

More generally, any two distinct \(T_1\)-stationary ergodic laws
\(\mu,\nu\) on \(\mathcal E\) have mutually singular restrictions
to \(\mathcal E_+\). In particular their whole-line laws are
mutually singular.
\end{lemma}
\begin{proof}
For two square-integrable white-noise observables, take their
conditional expectations on the noise in \([-A,A]\). These
approximations converge in \(L^2\) as \(A\to\infty\), since the
finite-window sigma-fields generate the full noise sigma-field.
For fixed \(A\), sufficiently distant integer shifts put the two
approximations in disjoint noise windows and make them independent.
The \(L^2\) approximation errors are uniform under shifts.
Cauchy--Schwarz therefore gives
\[
 \E[f(W)g(\theta_nW)]\longrightarrow \E[f(W)]\E[g(W)],
                    \qquad n\to\infty.
\]
This proves mixing. The equivariance identity, intersected over
integer iterates, makes \(\Phi\) a factor of this shift. Composing
observables with \(\Phi\) proves stationarity and mixing of its
law. Restriction to nonnegative endpoints is another factor.
Mixing implies ergodicity for each of these laws.

For the final assertion, finite increment cylinders form a
generating class for the Borel sigma-field. Distinct laws therefore
admit a bounded measurable finite increment cylinder \(F\) with
different expectations. Choose an integer \(k\geq0\) so that every
endpoint used by \(g=F\circ T_1^k\) is nonnegative. Stationarity
gives
\[
 a=\int g\dd\mu\ne b=\int g\dd\nu.
\]
Birkhoff's theorem and ergodicity imply
\[
 \lim_{m\to\infty}\frac1m\sum_{n=0}^{m-1}g(T_1^n\mathbf x)
 =\begin{cases}
   a,&\mu\text{-almost surely},\\
   b,&\nu\text{-almost surely}.
  \end{cases}
\]
Each summand depends only on nonnegative-time increments. The two
different limits thus define disjoint Borel full-measure sets in
\(\mathcal E_+\), proving singularity of the restricted laws.
Their preimages under restriction separate the whole-line laws.
\end{proof}

\section{Nonuniqueness of complete natural enhancements}\label{sec:central-gauges}
At nonreciprocal indices, complete-past causality, stationarity and
exact self-similarity leave freedom in the highest level.
We exhibit this freedom through an explicit central perturbation
that preserves every lower level and uses only the original Gaussian past.
We combine the homogeneous Wiener-kernel construction of
\citet{BaiTaqqu2014} with central path multiplication
\citep{FrizOberhauser2009,BrunedEtAl2019} to obtain a full
finite-depth enhancement and an explicit change in its integrals.

\subsection{A homogeneous ordered kernel}

Assume $0<H\leq1/4$ and $1/H\notin\mathbb N$.  Put
\[
 N=\lfloor1/H\rfloor\geq4,\qquad
 \gamma=NH\in\left(\frac{N}{N+1},1\right),\qquad
 \rho=\gamma-\frac N2.
\]
Let
\[
 \mathcal C=\{y\in\R^{N-1}:y_1>\cdots>y_{N-1}>0\}
\]
and choose a nonzero real $g\in C_c^\infty(\mathcal C)$, extended
by zero outside $\mathcal C$.
Compact support inside this open cone already separates the support
from the ordering diagonals and the boundary $y_{N-1}=0$.
Set
\begin{align}
 P_t(u_1,\ldots,u_N)
  &=\ind_{\{u_1<\cdots<u_N<t\}}(t-u_N)^\rho
      g\left(\frac{u_N-u_1}{t-u_N},\ldots,
              \frac{u_N-u_{N-1}}{t-u_N}\right),\label{cg:eq:P}\\
 D_{st}&=P_t-P_s.
 \label{cg:eq:D}
\end{align}
No separate Wiener integral of $P_t$ will be used.

\begin{lemma}\label{cg:lem:kernel}
For $s<t$, $D_{st}\in L^2(\R^N)$ and
\begin{equation}\label{cg:eq:Dnorm}
 \norm{D_{st}}_2^2=c_{H,N,g}|t-s|^{2\gamma},\qquad
 \frac{\norm g_2^2}{2\gamma}\leq c_{H,N,g}<\infty.
\end{equation}
In particular the constant is strictly positive.
\end{lemma}

\begin{proof}
Kernel translation and dilation are exact:
\[
 P_{t+b}(u+b\mathbf1)=P_t(u),\qquad
 P_{\lambda t}(\lambda u)=\lambda^\rho P_t(u).
\]
Thus it suffices to prove the assertion for $(s,t)=(0,1)$.
For $v=u_N<t$, write $x=t-v$ and $y_j=v-u_j$.
The vector $y$ lies in $\mathcal C$, and integration in these
variables yields
\begin{equation}\label{cg:eq:slice}
 \int_{u_1<\cdots<u_{N-1}<v}|P_t(u_1,\ldots,u_{N-1},v)|^2
                         du_1\cdots du_{N-1}
       =x^{2\rho+N-1}\norm g_2^2
       =x^{2\gamma-1}\norm g_2^2.
\end{equation}
This is integrable near $x=0$ because $\gamma>0$.
On $-1\leq v<1$, the estimate
$|P_1-P_0|^2\leq2P_1^2+2P_0^2$ and
\eqref{cg:eq:slice} give a finite integral.  Both kernels vanish for
$v\geq1$.

For $v<-1$, put $x=-v>1$ and $y_j=v-u_j$, and define
\[
 Q(x,y)=x^\rho g(y/x),\qquad
 G(y)=\rho g(y)-y\cdot\nabla g(y).
\]
The time derivative here is taken on a remote region with no moving
boundary.  It satisfies
\[
 \partial_xQ(x,y)=x^{\rho-1}G(y/x),\qquad
 \norm{\partial_xQ(x,\cdot)}_{L^2(\mathcal C)}^2
                         =x^{2\gamma-3}\norm G_2^2.
\]
Minkowski's inequality gives
\begin{align*}
 \norm{D_{01}\ind_{\{u_N<-1\}}}_2
 &\leq\int_0^1
    \norm{\partial_xQ(x+r,y)}_{L^2(x>1,y\in\mathcal C)}dr\\
 &\leq\frac{\norm G_2}{\sqrt{2-2\gamma}}<\infty.
\end{align*}
This proves the infrared cancellation precisely when $\gamma<1$.
On $0<v<1$, $P_0$ vanishes, and this region alone contributes
\[
 \norm g_2^2\int_0^1(1-v)^{2\gamma-1}dv
                           =\frac{\norm g_2^2}{2\gamma}>0.
\]
Finally the transformation $u=s\mathbf1+hU$, $h=t-s$, multiplies
the squared norm by $h^{2\rho+N}=h^{2\gamma}$.
This proves \eqref{cg:eq:Dnorm}.
\end{proof}

\begin{remark}[Explicit constant and the reciprocal restriction]
The constant is
\[
 c_{H,N,g}=\frac{\norm g_2^2}{2\gamma}
 +\int_0^\infty\int_{\mathcal C}
 \left|(x+1)^\rho g\left(\frac y{x+1}\right)
                 -x^\rho g\left(\frac yx\right)\right|^2dy\,dx.
\]
An individual $P_t$ has squared norm
$\norm g_2^2\int_0^\infty x^{2\gamma-1}dx=\infty$.
At a reciprocal index with $NH=1$, the remote derivative power is
$x^{-1}$, so this argument does not define the proposed top gauge.
The nonreciprocal assumption is an analytic condition, not only a
choice of rough-path depth.
\end{remark}

\subsection{The ordered It\^o tensor and repeated-color orthogonality}

Let $W$ have $d\geq2$ independent two-sided Brownian coordinates.
Take a deterministic square-integrable kernel $D$ on the ordered
simplex $\{u_1<\cdots<u_N\}$. For a word $w=i_1\cdots i_N$,
define the iterated It\^o integral
\[
 T^w(D)=\int_{u_1<\cdots<u_N}
          D(u_1,\ldots,u_N)
          \,dW_{u_1}^{i_1}\cdots dW_{u_N}^{i_N}.
\]
Completion from elementary ordered kernels defines this integral.
Its associated tensor is
\[
 T(D)=\sum_{|w|=N}T^w(D)e_w.
\]

\begin{lemma}\label{cg:lem:orthogonality}
For any two deterministic ordered $L^2$ kernels $D,E$ and any
two length-$N$ words $v,w$,
\begin{equation}\label{cg:eq:word-orth}
 \E[T^v(D)T^w(E)]
        =\ind_{\{v=w\}}\langle D,E\rangle_{L^2(\R^N)}.
\end{equation}
In particular this holds for repeated colors.
Every $T^w(D)$ belongs to the $N$th homogeneous Wiener chaos.
\end{lemma}

\begin{proof}
Apply cross It\^o isometry to the last integration variable.  If the
last colors differ, the expectation is zero even though their
integrands are random.  If the last colors agree, the expectation
is the integral of the covariance of the two prefix integrals, with
the last variable fixed as a deterministic kernel parameter.
Induct on the word length.  This recursively forces equality of
every color in its ordered position and gives the right-hand side
of \eqref{cg:eq:word-orth}.  The argument first applies to elementary
kernels and extends by the isometry.
An iterated It\^o integral of a deterministic ordered kernel is a
multiple Wiener integral of order $N$, including when a color
repeats; no Stratonovich contractions are inserted.
\end{proof}

We apply this definition to $D_{st}$ and write $T_{st}=T(D_{st})$.
It is additive as a tensor increment, by linearity in $D_{st}$.

\subsection{The right Dynkin projection and an exact variance formula}

The bracket is the tensor commutator $[a,b]=a\otimes b-b\otimes a$.
For homogeneous tensors of degree $N$, define
\begin{equation}\label{cg:eq:Dynkin}
 \operatorname{Dyn}_N(v_1\otimes\cdots\otimes v_N)
   =\frac1N[v_1,[v_2,\ldots,[v_{N-1},v_N]\ldots]].
\end{equation}
The direction of nesting in this formula is part of the definition.

\begin{lemma}[Normalization and naturality]\label{cg:lem:Dynkin}
The map $\operatorname{Dyn}_N$ is a projection onto $\mathfrak L_N(\R^d)$ and commutes
with every linear coordinate map.  Its normalization is $1/N$,
not $1/N!$.
\end{lemma}

\begin{proof}
We derive the Dynkin--Specht--Wever identity for the stated nesting
convention. Give each generator of the tensor Hopf algebra primitive coproduct,
let $S$ be its antipode, and let $\mathcal Y$ be the grading
derivation, equal to multiplication by $m$ on degree $m$.
The convolution $D_{\rm right}=\mathcal Y*S$ obeys
\[
 D_{\rm right}(v)=v,\qquad
 D_{\rm right}(vw)=[v,D_{\rm right}(w)]
                    \quad\text{for nonempty }w.
\]
Indeed, expanding $\Delta(vw)=(v\otimes1+1\otimes v)\Delta w$
and using the grading derivation rule gives
\[
 D_{\rm right}(vw)
   =v\,m(\mathrm{id}\otimes S)\Delta w
                  +vD_{\rm right}(w)-D_{\rm right}(w)v.
\]
The first term is $v\varepsilon(w)=0$ for nonempty $w$.
Thus $D_{\rm right}$ is exactly the unnormalized right-nested
commutator in \eqref{cg:eq:Dynkin}, and its range in degree $N$ lies
in $\mathfrak L_N$.
Every homogeneous Lie polynomial $L\in\mathfrak L_N$ is primitive, so
\[
 (\mathcal Y*S)(L)=N L.
\]
Consequently $D_{\rm right}/N$ is the identity on $\mathfrak L_N$,
and is a projection onto it.  Naturality follows directly from the
commutator formula, or from the naturality of the Hopf operations.
This projection need not be orthogonal for the Euclidean tensor
inner product; orthogonality is not used.
\end{proof}

Define the Lie-valued increment
\begin{equation}\label{cg:eq:Z}
 Z_{st}=\operatorname{Dyn}_N T_{st}\in\mathfrak L_N.
\end{equation}

\begin{proposition}[Exact nondegeneracy in dimension two]
\label{cg:prop:nondeg}
Let $w=1^{N-1}2$ and $w'=1^{N-2}21$.  Then
\begin{equation}\label{cg:eq:Zcomponent}
 Z^w_{st}=\frac1N\left(T^w(D_{st})-T^{w'}(D_{st})\right),
 \qquad
 \E|Z^w_{st}|^2
      =\frac{2c_{H,N,g}}{N^2}|t-s|^{2\gamma}>0.
\end{equation}
\end{proposition}

\begin{proof}
Commutators preserve the multiset of coordinate labels.  Thus only
input words containing $N-1$ copies of $1$ and one copy of $2$
can contribute to output word $w$.  If that unique $2$ is not in
one of the last two positions, the innermost bracket is $[1,1]=0$.
The two remaining input words are $w$ and $w'$, and
\[
 \operatorname{Dyn}_N(e_w)=\frac1N\operatorname{ad}_{e_1}^{N-1}(e_2),
 \qquad
 \operatorname{Dyn}_N(e_{w'})=-\frac1N\operatorname{ad}_{e_1}^{N-1}(e_2).
\]
The coefficient of $e_w$ in this iterated adjoint is $1$.
This proves the first identity in \eqref{cg:eq:Zcomponent}.
The last colors of $w,w'$ differ, so their ordered It\^o integrals
are orthogonal by Lemma~\ref{cg:lem:orthogonality}; each has variance
$\norm{D_{st}}_2^2$.  Lemma~\ref{cg:lem:kernel} proves the second
identity.  No repeated-color symmetrization can cancel it.
\end{proof}

\subsection{Regularity, causality, and the joint noise symmetries}

Let
\[
 (\theta_bW)(A)=W(A+b),\qquad
 (S_cW)(A)=c^{-1/2}W(cA),\qquad c>0.
\]
These preserve the law of independent Brownian white noises.

\begin{proposition}\label{cg:prop:Zproperties}
The process $Z$ is additive, centered, belongs to the $N$th
homogeneous Wiener chaos, and is measurable from Wiener noise up to
time $t$ for $s<t$.  For every finite $p$,
\[
 \norm{Z_{st}}_{L^p}\leq C_{p,N,d,H,g}|t-s|^\gamma.
\]
It has a continuous version with increment H\"older exponent
$\eta$ for every $\eta<\gamma$.
Its symmetries hold as noise-functional identities:
\begin{align}
 Z_{s+b,t+b}(W)&=Z_{st}(\theta_bW),\label{cg:eq:Zshift}\\
 Z_{cs,ct}(W)&=c^\gamma Z_{st}(S_cW),\label{cg:eq:Zscale}\\
 Z_{st}(AW)&=A^{\otimes N}Z_{st}(W)
                       \quad\text{for fixed linear }A.
                       \label{cg:eq:Zlinear}
\end{align}
In particular the process has jointly stationary increments and
exact $\gamma$ self-similarity, with the corresponding orthogonal
spatial covariance.
\end{proposition}

\begin{proof}
Additivity follows from $D_{st}=D_{su}+D_{ut}$ and linearity of the
It\^o tensor and Dynkin map.  The chaos and mean assertions follow
from the deterministic iterated It\^o definition.
Lemma~\ref{cg:lem:orthogonality}, the finite-dimensional Dynkin map,
and \eqref{cg:eq:Dnorm} give the $L^2$ bound.
Hypercontractivity in the $N$th chaos gives all $L^p$ bounds.
Arbitrarily high moments and Kolmogorov's criterion give each
$\eta<\gamma$, simultaneously with additive continuous increments.
The endpoint-ordered support of $D_{st}$ has all its variables
strictly before $t$, proving causality in the completed past
filtration.  As usual, a path anchored at parameter time zero need
not itself be adapted at negative times; the assertion concerns
increments with $s<t$.

Kernel translation gives \eqref{cg:eq:Zshift}.  Kernel scaling supplies
$c^\rho$ and the $N$ Brownian integrators supply $c^{N/2}$, giving
the factor $c^{\rho+N/2}=c^\gamma$ in \eqref{cg:eq:Zscale}.
Iterated It\^o integration is multilinear in its coordinate
integrators.  Lemma~\ref{cg:lem:Dynkin} therefore gives
\eqref{cg:eq:Zlinear}.  All identities are first identities of finite
or elementary-kernel approximations and then hold by $L^p$
completion.  These functional identities ensure joint symmetry with
other equivariant functionals of the same noise.
\end{proof}

\subsection{Complete geometric enhancements}

Let $\bX$ be the nonreciprocal infinite-past enhancement constructed
in Section~\ref{sec:causal-construction} from the same white noises.  It has, jointly with
the noise, the functional identities
\[
 X_{s+b,t+b}(W)=X_{st}(\theta_bW),\qquad
 X_{cs,ct}(W)=\delta_{c^H}X_{st}(S_cW),
\]
and commutes with fixed linear spatial maps.  The finite-cutoff
identities and convergence through arbitrary real cutoffs prove
these functional versions, not merely their separate marginal laws.
Its coordinates have every finite moment.

\begin{theorem}[Full top-central gauge]\label{cg:thm:fullgauge}
For every fixed $\lambda\in\R$, define
\begin{equation}\label{cg:eq:gaugedX}
 \bX^{[\lambda]}_{st}=\bX_{st}\exp(\lambda Z_{st})
                   \quad\text{in }G^N(\R^d).
\end{equation}
For every $1/(N+1)<\alpha<H$, this is a complete geometric
$\alpha$-H\"older step-$N$ rough path over the same fBm first
level.  It has the same levels $1,\ldots,N-1$ as $\bX$ and
\[
 \pi_N\bX^{[\lambda]}_{st}=\pi_N\bX_{st}+\lambda Z_{st}.
\]
It is causal for the complete Wiener past, jointly stationary,
exactly degree-$kH$ self-similar at every retained level $k$, and
natural under fixed linear coordinate maps.
By Proposition~\ref{prop:complete-past}, it is also causal for the
complete fBm increment past.
\end{theorem}

\begin{proof}
The top homogeneous Lie space $\mathfrak L_N$ is central in the step-$N$
Lie algebra.  Its exponential is group-like and, at this truncation,
$\exp(\lambda Z)=1+\lambda Z$.
Centrality and additivity give
\begin{align*}
 \bX^{[\lambda]}_{su}\bX^{[\lambda]}_{ut}
 &=\bX_{su}\bX_{ut}\exp(\lambda(Z_{su}+Z_{ut}))
   =\bX_{st}\exp(\lambda Z_{st}).
\end{align*}
Thus all Chen and shuffle/group identities hold at every retained
degree, including all the intermediate levels needed for integration.
The only altered coordinate degree is $N$.

For the desired $\alpha$, choose $\eta$ with
$N\alpha<\eta<\gamma=NH$.
Proposition~\ref{cg:prop:Zproperties} gives
$|Z_{st}|\lesssim|t-s|^\eta$ on compact intervals, which is more
regular than the required degree-$N$ estimate $|t-s|^{N\alpha}$.
All lower rough estimates are those of $\bX$.
Moreover $N=\lfloor1/\alpha\rfloor$, so these are exactly all
levels required for the $\alpha$ rough path.
Working first at a slightly larger exponent below $H$ and applying
the standard signature approximation at the desired smaller
exponent proves geometricity, as in the existing construction.

Causality passes to \eqref{cg:eq:gaugedX} from its two factors.
The noise-functional shift relations are joint and hence preserve
stationarity after multiplication.
Since $\gamma=NH$, the factor $c^\gamma$ in
\eqref{cg:eq:Zscale} is exactly the top-degree action of
$\delta_{c^H}$, proving self-similarity at every level.
The linear naturality follows from \eqref{cg:eq:Zlinear}, the naturality
of $\bX$, and the functoriality of the Lie exponential.
\end{proof}

\subsection{Different integrals and uncountably many different laws}

\begin{corollary}[A concrete one-form witness]\label{cg:cor:integral}
For the smooth polynomial one-form
$F(x)=(x^1)^{N-1}dx^2$, the canonical rough integrals satisfy
\begin{equation}\label{cg:eq:integral-difference}
 \int_s^tF(B_r)\,d\bX_r^{[\lambda]}
       -\int_s^tF(B_r)\,d\bX_r
          =\lambda(N-1)!Z^w_{st},\qquad w=1^{N-1}2.
\end{equation}
In particular the selected integration maps, not just their
enhancement coordinates, are different.
\end{corollary}

\begin{proof}
All compensated-sum terms below degree $N$ are unchanged.
At degree $N$ the sole nonzero coefficient is
$\partial_1^{N-1}F_2=(N-1)!$ on word $w$.
Thus on every partition the difference of the compensated sums is
\[
 \lambda(N-1)!\sum_{[u,v]}Z^w_{uv}
                      =\lambda(N-1)!Z^w_{st}.
\]
Pass to the rough-integral limits. Localization to the compact
sample-path range supplies the bounded derivatives required for
integration of this smooth polynomial one-form.
Nondegeneracy follows from \eqref{cg:eq:Zcomponent}.
\end{proof}

\begin{proposition}[Uncountably many distinct laws]\label{cg:prop:laws}
The models in \eqref{cg:eq:gaugedX} are pairwise different as coupled
random models.  There is a deterministic half-line of parameters
on which their whole-line rough-path laws are pairwise different.
\end{proposition}

\begin{proof}
For $\lambda\ne\mu$, the top coordinate difference at $(0,1)$ is
$(\lambda-\mu)Z^w_{01}$, with nonzero variance.
For the laws, put
\[
 a=\E|Z^w_{01}|^2=\frac{2c_{H,N,g}}{N^2}>0,\qquad
 b=\E[X^w_{01}Z^w_{01}],\qquad
 m(\lambda)=\E|X^w_{01}+\lambda Z^w_{01}|^2.
\]
All quantities are finite, without an independence assumption.
Then
\[
 m(\lambda)=\E|X^w_{01}|^2+2b\lambda+a\lambda^2
\]
is strictly increasing on $[-b/a,\infty)$.
Every two different parameters on this half-line therefore give
different second moments, and hence different rough-path laws.
The half-line is uncountable.
\end{proof}

\begin{corollary}[Mutual singularity under the white-noise factor]
\label{cg:cor:singular}
On the canonical two-sided Wiener space, the stationary whole-line
rough-path laws are pairwise mutually singular for parameters in the
monotone half-line of Proposition~\ref{cg:prop:laws}.
\end{corollary}

\begin{proof}
The enhancement and the gauge are measurable noise functionals and
obey the joint shift identity, so each gauged stationary law is a
factor of the unit white-noise shift. The identities hold jointly
for integer shifts and rational endpoints on a common full-measure
set, and extend by continuity.
Lemma~\ref{lem:ergodic-singularity} gives ergodicity of each law.
For two different parameters in the monotone half-line,
Proposition~\ref{cg:prop:laws} distinguishes the expectations
$m(\lambda)$ of the finite-increment square observable.
A bounded truncation of that square separates the two laws, so
the same lemma proves their mutual singularity, already after
restriction to the positive half-line.
\end{proof}

\subsection{Relation to generalized Hermite kernels}

The process meets Bai--Taqqu's general kernel criterion
\cite[Proposition 3.1]{BaiTaqqu2014} with order $N$ and self-similarity
exponent $\gamma$:
\[
 (P_{ct}-P_0)(cu)=c^{\gamma-N/2}(P_t-P_0)(u),\qquad
 P_{t+r}(u)-P_t(u)=(P_r-P_0)(u-t\mathbf1).
\]
There is also a more precise identification in the present parameter
range.  Recall $G(y)=\rho g(y)-y\cdot\nabla g(y)$ and define
\[
 q(x_1,\ldots,x_N)
  =\ind_{\{x_1>\cdots>x_N>0\}}x_N^{\rho-1}
      G\left(\frac{x_1-x_N}{x_N},\ldots,
              \frac{x_{N-1}-x_N}{x_N}\right).
\]
For fixed strictly ordered $u$, $r\mapsto P_r(u)$ is smooth even
across $r=u_N$: its angular support forces it to vanish in a
neighborhood of that point.  Direct differentiation therefore gives
\begin{equation}\label{cg:eq:Hermiteform}
 P_t(u)-P_s(u)=\int_s^t q(r\mathbf1-u)dr.
\end{equation}
The kernel $q$ is nonzero and homogeneous of degree
\[
 \zeta=\rho-1=\gamma-\frac N2-1
              \in\left(-\frac{N+1}{2},-\frac N2\right),
\]
where the lower strict inequality uses
$\gamma>N/(N+1)>1/2$.
If $G$ vanished identically, Euler's homogeneous-function equation
would force $g$ to be homogeneous of degree $\rho$ along rays;
compact support and nonzero $g$ rule that out.

The generalized-Hermite integrability condition is
\[
 \int_{(0,\infty)^N}|q(x)q(x+\mathbf1)|dx<\infty.
\]
Here it follows directly from the angular support.  On the support
of $q(x)$ all $x_j$ are comparable with $x_N$, and some ratio
$(x_j-x_N)/x_N$ is bounded below by a positive constant.
For sufficiently small $x_N$, all ratios
$(x_j-x_N)/(1+x_N)$ are too small to lie in the angular support
of $q(x+\mathbf1)$.  Thus the product vanishes near the only
possible origin singularity.  At infinity its integrated slice is
bounded by
\[
 C x_N^{2\zeta+N-1}=C x_N^{2\gamma-3},
\]
which is integrable because $\gamma<1$.
On the intermediate compact region it is bounded.

Consequently $q$ satisfies precisely the generalized Hermite kernel
conditions of \cite[Definition 3.2]{BaiTaqqu2014}, and
\eqref{cg:eq:Hermiteform} is the kernel form of their Theorem 3.5.
Using labeled independent colors produces the tensor-valued
version by the same Wiener-integral proof; applying the classical
Dynkin projection produces the Lie-valued process needed here.

\Needspace{15\baselineskip}
\section{Closed one-forms and enhancement independence}
\label{sec:closed-forms}
\begin{proposition}\label{prop:closed-forms}
Fix \(0<\alpha<1/2\) and let \(F\) be a smooth one-form on \(\R^d\),
with enough derivatives for \(\alpha\)-H\"older rough integration.
The integral is the canonical one-form rough integral for the
specified enhancement; its controlled derivatives are induced
by the spatial derivatives of \(F\).
The following are equivalent.
\begin{enumerate}
\item The rough integral of \(F\) is the same for every geometric
 \(\alpha\)-H\"older enhancement of every first-level path.
\item \(F\) is closed: \(\partial_iF_j=\partial_jF_i\) for all \(i,j\).
\end{enumerate}
\end{proposition}
\begin{proof}
If \(F\) is closed, the Poincar\'e lemma on \(\R^d\) gives \(F=\dd U\).
The geometric rough chain rule yields
\(\int_s^tF(x)\dd\bX=U(x_t)-U(x_s)\), independently of the lift.

Conversely fix \(z\in\R^d\) and \(i\ne j\).
The constant path \(x_t=z\) has both the trivial enhancement and
the pure-area enhancement
\[
 \bX_{s,t}=\exp\bigl((t-s)[e_i,e_j]\bigr).
\]
Its even levels have regularity at least \(k/2\), its odd levels
vanish, and it is geometric at every \(\alpha<1/2\).
One can use standard geometric approximation at an exponent
slightly larger than \(\alpha\); equivalently, shrinking repeated
planar loops give the same lift.
At level two its coordinates are
\((t-s)(e_i\otimes e_j-e_j\otimes e_i)\).
The compensated integral equals
\[
 (t-s)\bigl(\partial_iF_j(z)-\partial_jF_i(z)\bigr).
\]
Higher even levels contribute zero: for \(k\geq4\) the first
antisymmetric pair lies among the \(k-1\) symmetric derivative
indices of \(D^{k-1}F\).
The trivial lift gives zero, forcing the curl to vanish at every \(z\).
\end{proof}

\section{Gaussian information and directed lifts}\label{app:information}
\subsection{Gaussian projections and intersections}
\begin{lemma}\label{lem:Gaussian-projection}
For two closed subspaces \(M_1,M_2\) of a real Gaussian
Hilbert space, let \(\mathcal G(M_j)\) be their completed
generated sigma-fields. Then
\[
 \mathcal G(M_1)\cap\mathcal G(M_2)
                  =\mathcal G(M_1\cap M_2)
                         \quad\text{modulo completion}.
\]
\end{lemma}
\begin{proof}
For a closed Gaussian subspace \(M\), conditional expectation
acts on the \(n\)-th Wiener chaos as
\[
 \E[I_n(f)\mid\mathcal G(M)]
                         =I_n(P_M^{\otimes n}f).
\]
To verify this, split \(h=P_Mh+(1-P_M)h\) in a Wick
exponential. The two Gaussian terms are independent, and
conditioning gives
\[
 \E[e^{uW(h)-u^2\|h\|^2/2}\mid\mathcal G(M)]
                  =e^{uW(P_Mh)-u^2\|P_Mh\|^2/2}.
\]
Compare Hermite coefficients, then use polarization and
the Wiener isometry.

Put \(K=M_1\cap M_2\). A tensor belonging to both
\(M_1^{\otimes n}\) and \(M_2^{\otimes n}\), contracted
against arbitrary vectors in all but one slot, has its
remaining vector in \(K\). Thus the \(K\)-projection in
each slot fixes it, and
\[
 M_1^{\otimes n}\cap M_2^{\otimes n}=K^{\otimes n}.
\]
Intersect with symmetric tensors. A square-integrable
variable measurable in both fields has every chaos
kernel in \(K^{\odot n}\), so is \(\mathcal G(K)\)-measurable.
Apply this to bounded indicators to obtain the full
completed sigma-field identity.

\end{proof}

\subsection{The interval-locality obstruction}\label{sec:locality}
The mixed-chaos argument of
\citet[Theorem~3.7(a)]{ChevyrevFerrucci2026} rules out square-integrable
interval-local area even when each coordinate may depend on every spatial component.
The proof below therefore applies without the spatial-index locality
condition in their definition.

\subsubsection{The fractional increment space and interval support}
Put $\delta=1/2-H\in[1/4,1/2)$. The scalar increment Hilbert
space $\cH$ is the closed span of interval indicators in the
concrete Fourier space
\[
 \dot H^\delta(\R)=\left\{f:\int_\R
        |\widehat f(\xi)|^2|\xi|^{2\delta}\dd\xi<\infty\right\}.
\]
Here \(\widehat f(\xi)=\int_\R e^{-iu\xi}f(u)\dd u\),
its representative satisfies the reality symmetry, and the fixed
normalization constant is chosen so that
$\|\ind_{(a,b)}\|_{\cH}=|b-a|^H$.
The isonormal map $\mathcal B$ sends
$\ind_{(a,b)}$ to $B_b^H-B_a^H$. This follows directly
from the first-level spectral representation.
The displayed space consists of weighted Fourier $L^2$
functions, interpreted as distributions; no quotient by
constants is used. For an interval $I=[s,t]$, let
$\cH_I$ be the closed span of indicators contained in $I$.

\begin{lemma}[Local support and an ordinary $L^2$ embedding]
\label{lem:local-support}
The Fourier space embeds continuously and injectively into
$\mathcal S'(\R)$. Moreover
\[
 \cH_I=\{f\in\dot H^\delta(\R):\supp f\subset I\},\qquad
 \|f\|_{L^2(I)}\leq C_\delta |I|^\delta
                                       \|f\|_{\cH},\quad f\in\cH_I.
\]
Consequently $\cH_I\otimes\cH_I$ embeds continuously into
$L^2(I^2)$, with its natural distributional realization.
\end{lemma}
\begin{proof}
For a Schwartz test function $\varphi$, weighted
Cauchy--Schwarz bounds the distributional pairing by
\[
 C\|f\|_{\dot H^\delta}
     \left(\int|\widehat\varphi(\xi)|^2
                                  |\xi|^{-2\delta}\dd\xi\right)^{1/2}.
\]
The inverse weight is integrable near zero because
$2\delta<1$, and the test function controls infinity.
The Fourier representative is locally integrable, so this
map is injective. In particular Hilbert limits preserve
support in a fixed closed set.

Conversely, take $f$ in the Fourier space with support in $I$.
With $m=(s+t)/2$ and $0<a<1$, compress it to
$f_a(x)=f(m+(x-m)/a)$, in the distributional sense.
Its support is strictly inside $I$ and
\[
 \widehat f_a(\xi)=a e^{-i(1-a)m\xi}\widehat f(a\xi),\qquad
 \|f_a\|_{\dot H^\delta}=a^{1/2-\delta}\|f\|_{\dot H^\delta}.
\]
Translations and dilations are strongly continuous in this
weighted $L^2$ realization: first verify the assertion on
smooth Fourier functions compactly supported away from zero,
then use density and the displayed locally bounded norms.
Hence $f_a\to f$ as $a\uparrow1$.
For fixed $a$, convolution with a sufficiently small compact
smooth mollifier stays inside $I$ and converges in the
Fourier norm by dominated convergence. This approximates $f$
by $C_c^\infty((s,t))$.
For such a smooth $\varphi$, choose $s<a<b<t$ enclosing
its support. The identity
\[
 \varphi(x)=\int_a^b\varphi'(r)\ind_{(r,b)}(x)\dd r
\]
is a Bochner integral in $\cH$: the indicator-valued map
is continuous since the norm of a difference of nested
indicators is the $H$-th power of the intervening length.
Its Riemann sums place $\varphi$ in $\cH_I$.
This proves the support characterization.

Finally, for smooth $f$ supported in an interval of length
$L$, one has $\|\widehat f\|_\infty\leq C L^{1/2}\|f\|_2$.
Thus the Fourier energy on $|\xi|\leq r$ is at most
$CrL\|f\|_2^2$. Choose $r=c/L$ to make this at most
$\|f\|_2^2/2$. The remaining energy is bounded by
$r^{-2\delta}\|f\|_{\dot H^\delta}^2$; absorption gives
the asserted embedding. The density just proved extends
it to $\cH_I$. Tensoring this bounded injective map gives
the last assertion, consistently with distributional tests.
\end{proof}

\begin{theorem}[Full-vector interval-local area obstruction]
\label{thm:locality}
For independent $d$-dimensional fBm, $d\geq2$,
$0<H\leq1/4$, and $i\ne j$, there is no square-integrable
family $A^{ij}_{s,t}$ on all finite intervals satisfying
\[
 A^{ij}_{s,t}=A^{ij}_{s,u}+A^{ij}_{u,t}
                 +(B_u^{H,i}-B_s^{H,i})(B_t^{H,j}-B_u^{H,j})
\]
that is measurable, on each interval, from all the fBm
increments on that interval. No H\"older estimate or
spatial-index locality is required for this conclusion.
\end{theorem}
\begin{proof}
Project onto the Gaussian chaos with one variable of color
$i$, one of color $j$, and none of any other color.
This chaos is isometric to $\cH\otimes\cH$ through
$f\otimes g\mapsto\mathcal B^i(f)\mathcal B^j(g)$.
Conditional expectation onto a Gaussian subspace acts on
chaos kernels by the corresponding tensor projection;
this identity is proved in Lemma~\ref{lem:Gaussian-projection}.
Interval measurability therefore places the projected
kernel $F_{s,t}$ in $\cH_{[s,t]}\otimes\cH_{[s,t]}$.
Other spatial colors cannot cancel this orthogonal block.
Projection of Chen's relation gives
\[
 F_{s,t}=F_{s,u}+F_{u,t}
                   +\ind_{(s,u)}\otimes\ind_{(u,t)}.
\]
Iterating over a finite partition $\pi$ of $[s,t]$ yields
\[
 F_{s,t}=\sum_{I\in\pi}F_I+
                           \sum_{I<J\text{ in }\pi}\ind_I\otimes\ind_J.
\]
By Lemma~\ref{lem:local-support}, these are ordinary
$L^2$ functions. The first sum is supported in the union of
the cell squares. Outside those squares the second sum
is the ordered triangle $\ind_{\{s<x<y<t\}}$.
Use a countable refining sequence with mesh tending to
zero. The remaining cell squares shrink to the diagonal,
a Lebesgue-null set. Thus $F_{s,t}$ equals that triangle
almost everywhere. No smallness estimate for the first
sum is needed.

After translation and dilation, consider
$T(x,y)=\ind_{\{0<x<y<1\}}$ and put
$q(z)=\int_0^1e^{-izy}\dd y$. Direct integration gives
\[
 \widehat T(\xi,\eta)=\frac{q(\eta)-q(\xi+\eta)}{i\xi}.
\]
Choose $\varepsilon>0$ with $|q(z)|\geq3/4$ on
$|z|\leq\varepsilon$. On the strip
$\eta=-\xi+z$, $|z|\leq\varepsilon$, and sufficiently
large positive $\xi$, the other term has
$|q(\eta)|\leq1/4$. Hence $|\widehat T|\geq1/(2\xi)$.
Its tensor Fourier norm is bounded below by
\[
 C\int_{R_0}^\infty\xi^{4\delta-2}\dd\xi=\infty,
                  \qquad \delta\geq1/4.
\]
The product inverse weight is locally integrable, so the
tensor Fourier space also embeds injectively into
$\mathcal S'(\R^2)$. The same triangle cannot have a
second Fourier representative of finite norm.
This contradicts $F_{s,t}\in\cH_I\otimes\cH_I$.
\end{proof}

The stronger obstruction in \citet{Han2026} removes moment
assumptions on positive-probability Borel domains.
Its proof uses a Gaussian four-corner argument
and oscillatory Cameron--Martin packets. The square-integrable
form above suffices for the arguments that follow.
The complete-past construction allows dependence on information
before the interval and hence lies outside this locality obstruction.

\begin{corollary}[Finite observed history and stationarity]
\label{cor:finite-history-stationarity}
Let \(T>0\), \(d\geq2\) and \(0<H\leq1/4\).
An enhancement on \([0,T]\) with square-integrable
off-diagonal second coordinates satisfying Chen's relation
cannot be both adapted to
\[
 \mathcal F_t^{B,0}
   =\overline{\sigma\{B_r-B_0:0\leq r\leq t\}}
\]
and jointly stationary in its enhanced increments under all
time shifts staying in \([0,T]\). The joint stationarity
includes the first level, not merely individual area marginals.
\end{corollary}
\begin{proof}
Fix \(0<h\leq T\) and an off-diagonal coordinate.
Completed \(\mathcal F_h^{B,0}\)-measurability gives a Borel
function \(f_h\) on \(C_0([0,h];\R^d)\) with
\[
 A_{0,h}=f_h((B_r-B_0)_{0\leq r\leq h})
                                      \quad\text{almost surely}.
\]
Joint stationarity transports this probability-one Borel graph
to \(((B_{s+r}-B_s)_{r\leq h},A_{s,s+h})\).
Thus each second coordinate is measurable from the full vector
of increments on its own interval. The proof of
Theorem~\ref{thm:locality} uses only countably many nested
partitions of a fixed interval; applied inside \([0,T]\),
it gives the same contradiction. No whole-line extension is
needed.
\end{proof}

In particular, a finite-origin adapted selection as in
Remark~\ref{rem:finite-history} cannot have jointly stationary
enhanced increments in this range. Its full enhanced law on
any \([0,T]\) therefore differs from that of the stationary
selections constructed here. This excludes the stated
finite-observation filtration, without identifying a unique
possible enlargement of the initial information.

\subsection{Past--future information and time reversal}
\label{sec:time-reversal}
The locality obstruction also constrains time reversal of a
causal enhancement. We first identify the information shared
by a complete past and a complete future. This is an
intersection statement, with no conditional independence
or Markov assumption.

Write \(\mathcal F^-_t=\mathcal F_t^B\) for the completed
increment past. For \(I=[s,t]\), define
\begin{align*}
 \mathcal F^+_s&=\overline{\sigma\{B_b^i-B_a^i:
                         s\leq a<b,\ 1\leq i\leq d\}},\\
 \mathcal F_I&=\overline{\sigma\{B_b^i-B_a^i:
                  s\leq a<b\leq t,\ 1\leq i\leq d\}}.
\end{align*}

\begin{lemma}[Gaussian information intersections]
\label{lem:Gaussian-intersection}
For the fBm increment fields above,
\begin{equation}\label{eq:past-future-intersection}
             \mathcal F^-_t\cap\mathcal F^+_s=\mathcal F_{[s,t]}.
\end{equation}
\end{lemma}
\begin{proof}
Apply Lemma~\ref{lem:Gaussian-projection}.
For fBm, let \(\cH_-(t)\) and \(\cH_+(s)\) be the
closed spans of interval indicators respectively ending
before \(t\) and starting after \(s\). In the Fourier
realization of Lemma~\ref{lem:local-support}, their
intersection has distributional support in \([s,t]\).
That lemma's compression, mollification and indicator
approximation give
\[
             \cH_-(t)\cap\cH_+(s)=\cH_{[s,t]}.
\]
No characterization of every half-line-supported
distribution is required. Take the direct sum over the
independent spatial colors and apply the general Gaussian
identity. This proves \eqref{eq:past-future-intersection}.
\end{proof}

\begin{theorem}[Irreversibility of complete-past causal lifts]
\label{thm:causal-irreversibility}
Let \(d\geq2\) and \(0<H\leq1/4\). A continuous whole-line
\(G^2(\R^d)\)-valued multiplicative functional over fBm
that is causal for its
complete increment past and has square-integrable
second-level coordinates cannot have a law invariant
under the reversal
\begin{equation}\label{eq:enhanced-reversal}
       (\mathcal R\mathbf x)_{s,t}
                       =(\mathbf x_{-t,-s})^{-1}.
\end{equation}
Stationarity and self-similarity are not needed for this
obstruction. It applies to the first two levels of every
full geometric enhancement.
\end{theorem}
\begin{proof}
Anchor the raw path by \(B_0=0\), and put
\((\rho B)_r=B_{-r}\). Its Gaussian law \(\mu\) is invariant
under \(\rho\). A causal enhancement is determined by the
complete raw path: countably many increment coordinates
and continuity give a Borel map \(F\) with
\(\mathbf X=F(B)\) almost surely. One may use the Polish
space of continuous increment families on compact sets;
no nonseparable endpoint H\"older space is needed.
Let \(p\) extract the anchored first-level path, so
\(pF(B)=B\), and define
\[
           G(B)=\mathcal R F(\rho B).
\]
Then \(pG(B)=B\). If the enhanced law were invariant under
\(\mathcal R\), invariance of \(\mu\) under \(\rho\) would give
\(F_*\mu=G_*\mu\).
The Borel graph event
\(\{\mathbf x:\mathbf x=F(p\mathbf x)\}\) has full
\(F_*\mu\)-measure. The equality of laws therefore implies
\(G(B)=F(pG(B))=F(B)\) almost surely.
This uses the whole enhanced law, including its first level.

For \(s<t\),
\[
 G(B)_{s,t}=(F(\rho B)_{-t,-s})^{-1}
\]
is measurable from the complete future of \(B\) beginning
at \(s\): it is the past of \(\rho B\) ending at \(-s\).
The equality \(G(B)=F(B)\) makes each coordinate measurable
both in \(\mathcal F^-_t\) and in \(\mathcal F^+_s\).
Lemma~\ref{lem:Gaussian-intersection} makes it interval-local.
An off-diagonal second coordinate would then be the
square-integrable Chen family excluded by
Theorem~\ref{thm:locality}.
\end{proof}

The threshold is sharp for simultaneous causality and
time-reversal invariance: when \(H>1/4\), the canonical
Gaussian enhancement of \citet{CoutinQian2002}, or the
Young enhancement when applicable, is interval-local and
square-integrable. Piecewise linear interpolation on a
symmetric dyadic grid commutes with reversal and group
inversion. Passing to its canonical rough limit gives
a time-reversal-invariant law.

\begin{corollary}[Two directed selections at every small Hurst index]
\label{cor:forward-backward-models}
For every \(0<H\leq1/4\), let \(F(B)\) be the stationary,
exactly self-similar causal enhancement constructed here.
On the same complete input define
\[
 X^\rightarrow(B)=F(B),\qquad
 X^\leftarrow(B)=\mathcal R F(\rho B).
\]
Both are full geometric enhancements with the prescribed
first level, stationary increments, exact graded
self-similarity and linear coordinate naturality.
Their complete-history enhanced laws are different.
The forward selection uses the complete past; the backward
selection uses the complete future.
Their canonical integrals of a recentered linear one-form
\((x^i)\dd x^j\), \(i\ne j\), differ on some interval with
positive probability.
\end{corollary}
\begin{proof}
Reversal and group inversion preserve the geometric
increment identities. In particular
\[
 X^\leftarrow_{s,u}X^\leftarrow_{u,t}
  =\bigl(F(\rho B)_{-t,-u}
                   F(\rho B)_{-u,-s}\bigr)^{-1}
  =X^\leftarrow_{s,t}.
\]
The first level is
\(-((\rho B)_{-s}-(\rho B)_{-t})=B_t-B_s\).
Reversal conjugates a time shift to its opposite and
commutes with positive dilations and linear spatial maps.
Thus the reversed law inherits the stated undirected
symmetries and moment bounds.
Theorem~\ref{thm:causal-irreversibility} distinguishes the
two whole-line laws.

For any \(i\ne j\), agreement of their entire
\(ij\)-coordinate families on this same input would already
give an interval-local Chen area by
\eqref{eq:past-future-intersection}. Theorem~\ref{thm:locality}
excludes this. Continuity reduces equality to countably
many rational intervals, so on at least one interval
the difference is nonzero with positive probability.
The identity
\[
 \int_s^t(B_r^i-B_s^i)\dd_{\mathbf X}B_r^j=X_{s,t}^{ij}
\]
then gives the claimed integral witness.
\end{proof}

The backward model uses a different information class, and the
two models cannot share all second levels. This comparison
therefore lies outside the fixed-lower-level causal class of
Theorem~\ref{thm:critical-rigidity}.
The symmetries and factor property of the selected models
give a stronger comparison of their laws.

\begin{corollary}[Finite-horizon distinction and infinite-time singularity]
\label{cor:directed-laws}
The forward and backward models in
Corollary~\ref{cor:forward-backward-models} have different
enhanced laws on every \([0,T]\), \(T>0\).
Their enhanced laws restricted to \([0,\infty)\) are
mutually singular.
\end{corollary}
\begin{proof}
Suppose their restrictions to some \([0,T]\) had the
same law. Exact graded \(H\)-self-similarity of each model
transports this equality by the deterministic map
\[
 \mathbf x_{s,t}\longmapsto
               \delta_{\lambda^H}\mathbf x_{s/\lambda,t/\lambda}
\]
to every \([0,\lambda T]\), \(\lambda>0\).
Stationarity then transports equality to every finite
whole-line window. Finite increment cylinders generate
the Borel sigma-field of continuous whole-line enhanced
paths, so their whole-line laws would agree. This
contradicts Theorem~\ref{thm:causal-irreversibility}.

For singularity, the forward law is a measurable factor of the unit
white-noise shift by its noise-functional shift identity, hence is
ergodic by Lemma~\ref{lem:ergodic-singularity}. Its reverse is also
ergodic because \(\mathcal R\vartheta_1=\vartheta_{-1}\mathcal R\).
The two whole-line laws are distinct by the first part of the proof.
The same lemma therefore gives mutually singular restrictions to
\([0,\infty)\).
\end{proof}

Irreversibility follows from the interval-locality obstruction
and the Gaussian information identity. Reversible Markov rough lifts
\citep{BassHamblyLyons2002} and reversibility constraints on
deterministic area drifts in homogenization
\citep{GottwaldMelbourne2024} concern different assumptions
and objects.

\section{Comparison with a Fourier-normal-ordered lift}
\label{sec:fourier-comparison}
We compare the causal construction with a class of resonance-filtered
Fourier kernels containing the selection of \citet{Unterberger2010}.
We test causality using the full second-level kernel, including
both skeleton and boundary terms.

Use the unitary transform
\(\mathcal F_+f(\xi)=(2\pi)^{-1/2}\int e^{iu\xi}f(u)\dd u\)
in innovation-kernel variables. The fractional-derivative
multipliers in Appendix~\ref{sec:filtration} use the usual
negative-exponent convention. Write the first-level
spectral representation as
\[
 B_t-B_s=\int_{\R}
       \frac{e^{it\xi}-e^{is\xi}}{i\xi}a_H(\xi)
                                      \widehat W(\dd\xi),
 \qquad
 a_H(\xi)=c(H)(i\xi)^{1/2-H}.
\]
Here \(c(H)>0\) incorporates the unitary Fourier normalization
and is chosen consistently with \(\kappa_t\).
This is the Fourier transform of the prescribed causal
moving-average representation, including its outer phase:
\[
 \mathcal F_+(\kappa_r-\kappa_s)(\xi)
       =a_H(\xi)\frac{e^{ir\xi}-e^{is\xi}}{i\xi}.
\]
The harmonizable spectral noise is rotated by the
unit-modulus phase matching this first level. Every two-color
kernel then acquires the two corresponding phase factors.
Frequency restrictions, permutations and forest products
commute with this operation, so the reconstruction below uses
\(a_H\) in the prescribed causal coupling.
Put \(\Delta_{s,t}(\omega)=e^{it\omega}-e^{is\omega}\).
Our word \(12\) denotes
\(\int_s^t(B_r^1-B_s^1)\dd B_r^2\).
For \(\xi\) in color \(1\), \(\eta\) in color \(2\), and
\(\omega=\xi+\eta\), the formal canonical kernel is
\[
 C_{s,t}(\xi,\eta)=a_H(\xi)a_H(\eta)
 \left\{\frac{\Delta_{s,t}(\omega)}{(i\xi)(i\omega)}
       -\frac{e^{is\xi}\Delta_{s,t}(\eta)}{(i\xi)(i\eta)}\right\}.
\]
We will use its \(L^2\) slices, not assume that it belongs
to \(L^2(\R^2)\).

Choose a real Borel multiplier \(\chi\), finite almost everywhere and homogeneous
of degree zero on the punctured plane, with the symmetries
\(\chi(\xi,\eta)=\chi(\eta,\xi)=\chi(-\xi,-\eta)\).
Put \(E=\{|\eta|\leq|\xi|\}\).
No boundedness or cone-support assumption is imposed on \(\chi\).
Null boundaries do not affect the kernels. The two-vertex
restriction in \citet[Definition~3.2]{Unterberger2010}, for \(0<c<1\), is the hard member
\(\chi=\ind_{\{|\xi+\eta|>c\max(|\xi|,|\eta|)\}}\).
For general \(\chi\) define the complete two-color kernel by
\begin{equation}\label{eq:FNO-full-kernel}
 F_{s,t}=C_{s,t}+D_{s,t},
\end{equation}
where
\begin{equation}\label{eq:FNO-difference-kernel}
\begin{split}
 D_{s,t}(\xi,\eta)
  ={}&-a_H(\xi)a_H(\eta)(1-\chi(\xi,\eta))
                    \frac{\Delta_{s,t}(\xi+\eta)}{i(\xi+\eta)}
       \left\{\frac{\ind_E}{i\xi}
                         -\frac{\ind_{E^c}}{i\eta}\right\}.
\end{split}
\end{equation}
The quotient at \(\xi+\eta=0\) is interpreted continuously.
Formula \eqref{eq:FNO-full-kernel} includes the singleton
boundary and forest terms. To check it, on \(E\) use the
displayed formula for \(C\), and on \(E^c\) use integration
by parts to write the same kernel as
\[
 C_{s,t}=a_H(\xi)a_H(\eta)
 \left\{-\frac{\Delta_{s,t}(\xi+\eta)}{(i\eta)i(\xi+\eta)}
       +\frac{e^{it\eta}\Delta_{s,t}(\xi)}{(i\eta)(i\xi)}\right\}.
\]
Definitions~2.5 and 2.9 and Lemma~2.10 of
\citet{Unterberger2010} retain the singleton terms in these
two formulas; Definition~3.2 multiplies their two-vertex
skeleton terms by the hard \(\chi\). Applying the same
two-sector formula with the chosen multiplier defines our
degree-two comparison family. Subtraction gives
\eqref{eq:FNO-difference-kernel}.
That work uses the opposite word ordering; the inner color
\(1\) and outer color \(2\) here fix the convention.

\begin{lemma}[Exact angular integrability criterion]
\label{lem:Fourier-angular}
Put
\[
 J_H(\chi)=\int_{\R}
  \frac{|q|^{1-2H}|1-q|^{1-2H}}
       {\max(|q|,|1-q|)^2}|\chi(q,1-q)|^2\dd q.
\]
For \(0<H<1/2\), \(F_{s,t}\in L^2(\R^2)\) if and only if
\(J_H(\chi)<\infty\); then
\(\|F_{s,t}\|_2\leq C_{H,\chi}|t-s|^{2H}\).
For \(0<H\leq1/4\) and \(\chi\in C^1(\R^2\setminus\{0\})\),
this criterion is equivalent to
\(\chi(1,-1)=\chi(-1,1)=0\).
\end{lemma}
\begin{proof}
Write \(F=M+Q\), where \(M\) is the retained skeleton
and \(Q\) contains the two unchanged boundary terms.
On \(E\), integrate the squared boundary first over
\(|\xi|\geq|\eta|\); the other sector is symmetric. This gives
\[
 \|Q\|_2^2\leq C_H\int_0^\infty
       r^{-1-4H}\min\{h^2r^2,4\}\dd r\leq C_Hh^{4H},
                         \qquad h=t-s.
\]
Thus \(F\in L^2\) if and only if \(M\in L^2\).
Put \(\omega=\xi+\eta\), \(q=\xi/\omega\).
The Jacobian is \(|\omega|\), and the modulus of the sector
braces in \eqref{eq:FNO-difference-kernel} is
\(1/\max(|\xi|,|\eta|)\). Homogeneity and simultaneous
sign symmetry give the exact factorization
\[
 \|M\|_2^2=c(H)^4J_H(\chi)
       \int_{\R}|\Delta_{s,t}(\omega)|^2
                            |\omega|^{-1-4H}\dd\omega.
\]
Tonelli applies also when the angular integral is infinite.
The radial integral is finite, positive and proportional to
\(h^{4H}\). This proves the criterion and bound.

The angular weight is locally integrable and is comparable
to \(|q|^{-4H}\) at infinity. By homogeneity,
\(\chi(q,1-q)=\chi(1,-1+1/q)\) for \(q>0\), and
the other tail approaches the opposite resonant ray.
A nonzero resonant value makes the angular integral diverge
when \(4H\leq1\). If both values vanish, \(C^1\) regularity
bounds the multiplier on the tails by \(C|q|^{-1}\), which
makes the integrand integrable.
\end{proof}

When this criterion holds, the inverse transform of \(F\)
defines a real mixed second Wiener integral.
The defect \(D_{s,t}\) is additive in \([s,t]\), while
the canonical frequency identity supplies Chen's cross term.
Under exchange of \(\xi,\eta\), the sector braces change sign
and \(D\) is antisymmetric. Hence these coordinates, with
diagonals \(\tfrac12(B_t^i-B_s^i)^2\), satisfy the degree-two
Chen and shuffle identities. Finite-chaos hypercontractivity
gives \(L^p\) second-level bounds of order \(|t-s|^{2H}\)
for every finite \(p\). Metric Kolmogorov estimates and the
Fourier shift and dilation identities give continuous,
stationary, graded-self-similar degree-two data. This does
not construct the higher levels required at small \(H\).

\begin{lemma}[Reflection of the Fourier selection]\label{lem:Fourier-reversal}
For \(J_H(\chi)<\infty\), the complete first-and-second-level
law just constructed is invariant under \(\mathcal R\) of
\eqref{eq:enhanced-reversal}.
\end{lemma}
\begin{proof}
On spectral white noise define the real unitary reflection
\[
 \widehat W^\rho(\dd\xi)
   =r(\xi)\widehat W(-\dd\xi),\qquad
 r(\xi)=-\frac{a_H(-\xi)}{a_H(\xi)}.
\]
The phase has modulus one and \(r(-\xi)=\overline{r(\xi)}\),
so this is again real white noise of the same law, independently
in every color. Its first level is \(B_t(W^\rho)=B_{-t}(W)\).
Put
\[
 L_{s,t}(\xi,\eta)=a_H(\xi)a_H(\eta)
      \frac{\Delta_{s,t}(\xi)\Delta_{s,t}(\eta)}{(i\xi)(i\eta)}.
\]
The formal canonical frequency identity under reflection is
\[
 r(-\xi)r(-\eta)C_{-t,-s}(-\xi,-\eta)=L_{s,t}-C_{s,t}.
\]
For the defect, \(\Delta_{-t,-s}(-\omega)/(-i\omega)
=\Delta_{s,t}(\omega)/(i\omega)\), while the sector braces
change sign. Simultaneous sign symmetry of \(\chi\) gives
\(r(-\xi)r(-\eta)D_{-t,-s}(-\xi,-\eta)=-D_{s,t}\).
Adding these pointwise identities yields
\[
 r(-\xi)r(-\eta)F_{-t,-s}(-\xi,-\eta)=L_{s,t}-F_{s,t}.
\]
Only this combined identity is integrated: \(F,L\) are jointly
\(L^2\), although \(C,D\) separately need not be.
Thus for distinct colors
\(X^{ij}_{-t,-s}(W^\rho)=B^i_{s,t}B^j_{s,t}-X^{ij}_{s,t}(W)\).
The diagonal identity is immediate. Group inversion and the
first-level relation now give
\(\mathcal RX(W^\rho)=X(W)\), first at rational endpoints
and then everywhere by continuity. Invariance of the noise law
proves the assertion.
\end{proof}

\begin{proposition}[Noncausality of resonance-filtered coordinates]
\label{prop:FNO-noncausal}
For \(d\geq2\), \(0<H\leq1/4\), \(s<t\), and a multiplier
in the stated class with \(J_H(\chi)<\infty\), the off-diagonal
coordinate with kernel \(F_{s,t}\) is not
\(\mathcal F_t^W\)-measurable in the prescribed causal
innovation representation.
Any full enhancement having this coordinate therefore has
a complete-history law different from that of any
complete-past causal enhancement of this fBm.
\end{proposition}
\begin{proof}
Suppose one fixed off-diagonal coordinate on one interval were
past-measurable. For a pure mixed second Wiener integral this
means that its time kernel is supported in
\(({-\infty},t]^2\). Homogeneity and translation of the Fourier
kernel give, for \(b\in\R\) and \(c>0\),
\[
 K_{cs+b,ct+b}(u,v)
   =c^{2H-1}K_{s,t}((u-b)/c,(v-b)/c).
\]
This transports past support to every interval. The common
kernel prescription and color permutations give the same
conclusion for every off-diagonal coordinate; the diagonal
coordinates are already past-measurable. Hence the whole
continuous \(G^2\)-valued functional would be causal.
Proposition~\ref{prop:complete-past} identifies this past with
the complete fBm increment past. Its reversible law from
Lemma~\ref{lem:Fourier-reversal} contradicts
Theorem~\ref{thm:causal-irreversibility}, which requires only
these two levels.

Finally the completed increment past is generated by
countably many rational fBm increments. Measurability from
that past is a property of their joint law with the second
coordinate. A complete-history enhanced law contains all
these first-level increments, so it cannot agree with a
complete-past causal law.
\end{proof}

For a smooth example without an open zero cone, fix \(0<c<1\) and take
\(\chi(\xi,\eta)=\psi((\xi+\eta)^2/(\xi^2+\eta^2))\),
where \(\psi:[0,2]\to[0,1]\) is smooth, vanishes only at zero,
and equals one on \([c^2/2,2]\). It satisfies the angular
criterion and is noncausal. Smooth inner-gap cutoffs are
also included. Thus the sharp boundary of the
hard restriction is not responsible for the obstruction.

The criterion applies to the homogeneous skeleton multipliers above.
The algebraic framework of \citet[Section~2.5]{Unterberger2010}
allows more general linear skeleton maps, including choices outside
this class. Applying the full-enhancement conclusion requires a
separate construction of the higher levels.

\subsection{The two-vertex kernel-BPHZ prescription}
\label{sec:kernel-bphz}
The preceding class also contains an unbounded multiplier
arising from the kernel-level BPHZ subtraction of
\citet{Unterberger2013}.

Put \(\delta=1/2-H\). In a two-vertex normal-ordered tree,
let \(\xi\) be the inner momentum and \(\eta\) the outer
momentum, with \(|\eta|\leq|\xi|\). After solving momentum
conservation, the outgoing momentum is
\(\omega=\xi+\eta\). The amputated root numerator is
\(|\omega-\xi|^\delta|\xi|^\delta\), with the internal
propagator \((i\xi)^{-1}\). Taylor evaluation at the outgoing
momentum zero replaces it by \(|\xi|^{2\delta}\).
The retained numerator is therefore
\begin{equation}\label{eq:BPHZ-numerator}
             (|\eta|^\delta-|\xi|^\delta)|\xi|^\delta.
\end{equation}
This is the minimal subdiagram \((\xi_2,\zeta_3,\xi_3)\)
in \citet[Section~4, equations~(4.9)--(4.10)]{Unterberger2013}, with
\(\omega=\zeta_2\), \(\xi=\zeta_3\).
It implements the global forest rule (3.4), rather than
the useful-only formula (3.18), and its
zero-external-momentum operation (3.1).
It does not Taylor-evaluate the overall momentum-conservation
delta distribution. The singleton Taylor numerator is
\(|0|^\delta=0\), so the chosen singleton skeletons and
their boundary terms remain unchanged, consistently with
the first-level increment normalization in Definition~1.10(a).
At two vertices there is only this one nonzero subtraction.
The other frequency sector is obtained
by the same exchange and integration-by-parts reconstruction
used above. No noise contraction occurs for two distinct colors.

\begin{proposition}[Kernel-BPHZ comparison]\label{prop:kernel-BPHZ}
For this two-vertex prescription, the complete cross kernel
is \eqref{eq:FNO-full-kernel}--\eqref{eq:FNO-difference-kernel}
with
\begin{equation}\label{eq:BPHZ-multiplier}
 \chi_{\mathrm B}(\xi,\eta)
     =1-\left(\frac{\max(|\xi|,|\eta|)}
                         {\min(|\xi|,|\eta|)}\right)^{1/2-H}.
\end{equation}
Axis values are immaterial. For every \(0<H<1/2\) this
kernel belongs to \(L^2(\R^2)\), and the \(L^p\) norms of its
cross coordinates are bounded by \(C_p|t-s|^{2H}\),
\(1\leq p<\infty\).
For \(0<H\leq1/4\), each such coordinate is noncausal for
the prescribed complete innovation past. Any full enhancement
retaining this coordinate has a complete-history law different
from the causal selections constructed here.
\end{proposition}
\begin{proof}
On \(E\), divide \eqref{eq:BPHZ-numerator} by the unmodified
numerator \(|\eta\xi|^\delta\), giving
\(1-(|\xi|/|\eta|)^\delta\). Exchange of sectors gives
\eqref{eq:BPHZ-multiplier}. Rotating each harmonizable noise
by the first-level causal phase produces the same \(a_H\)
factors as before. In particular it rotates both terms in
the retained numerator; no counterterm phase is discarded.
The unchanged singleton and forest terms give precisely
the full \(C+D\) kernel.

For its angular integral, put
\(m=\max(|q|,|1-q|)\), \(l=\min(|q|,|1-q|)\).
The integrand simplifies exactly to
\[
                    m^{2\delta-2}(m^\delta-l^\delta)^2.
\]
It is bounded near \(q=0,1\). At infinity,
\(m-l=1\), so the mean-value theorem bounds it by
\(C|q|^{4\delta-4}=C|q|^{-2-4H}\).
Hence \(J_H(\chi_{\mathrm B})<\infty\).
The multiplier is real, homogeneous, exchange-symmetric and
even under simultaneous sign reversal. Lemma~\ref{lem:Fourier-angular}
and Proposition~\ref{prop:FNO-noncausal} apply.
\end{proof}

The multiplier's axis singularities and noncompact frequency
modification are covered by the finite-angular-integral criterion.
Theorem~0.1 of \citet{Unterberger2013} assumes \(1/H\notin\mathbb N\);
its following remark treats reciprocal indices with a loss
of regularity. The two-vertex \(L^2\) calculation above holds
also at those indices. At a reciprocal index the full-law
comparison applies to any extension retaining this two-vertex
kernel; existence of such an extension requires separate estimates
at higher levels. The calculation concerns this specific BPHZ
prescription.

\subsection{The causal selection in Fourier variables}
\label{sec:causal-Fourier}
The preceding obstruction uses real, simultaneously even
multipliers. The causal Volterra coordinate has a different
phase. We identify it explicitly. Throughout this subsection
\(0<H<1/2\), \(\delta=1/2-H\), and
\(\gamma=H+1/2=1-\delta\). The normalization in
Appendix~\ref{sec:fourier-comparison} is
\[
 c(H)=\frac{\mathfrak c_H\Gamma(\gamma)}{\sqrt{2\pi}}.
\]
Write \([iz+0]^{-2H}=\lim_{\epsilon\downarrow0}
(\epsilon+iz)^{-2H}\), with the principal branch, and set
\begin{equation}\label{eq:causal-Fourier-density}
\begin{split}
 \mathfrak d_H(\xi,\eta)
 ={}&\frac{\delta c(H)^2\Gamma(2H)}{\Gamma(\gamma)^2}
  \int_0^1 r^{-\delta}
  \left\{[i(r\xi+\eta)+0]^{-2H}
               -[i(\xi+r\eta)+0]^{-2H}\right\}\dd r.
\end{split}
\end{equation}
The integral defines a locally integrable function of
\((\xi,\eta)\). Values on null sets are immaterial.

\begin{proposition}[Fourier representation of the causal coordinate]
\label{prop:causal-Fourier}
The second-level coordinate of Theorem~\ref{thm:noncritical-causal}
has the complete spectral kernel
\begin{equation}\label{eq:causal-complete-Fourier}
 F^{\mathrm c}_{s,t}(\xi,\eta)
  =C_{s,t}(\xi,\eta)
    +\frac{\Delta_{s,t}(\xi+\eta)}{i(\xi+\eta)}
                                      \mathfrak d_H(\xi,\eta).
\end{equation}
The identity first holds in tempered distributions; its sum
belongs to \(L^2(\R^2)\). Define, away from the axes,
\begin{equation}\label{eq:causal-Fourier-multiplier}
 \chi_{\mathrm c}(\xi,\eta)
 =1+\frac{\mathfrak d_H(\xi,\eta)}
  {a_H(\xi)a_H(\eta)
       \{\ind_E/(i\xi)-\ind_{E^c}/(i\eta)\}}.
\end{equation}
It is homogeneous of degree zero, exchange-symmetric and
Hermitian, in the sense that
\(\chi_{\mathrm c}(-\xi,-\eta)=
\overline{\chi_{\mathrm c}(\xi,\eta)}\).
It has \(J_H(\chi_{\mathrm c})<\infty\), but is not real or
simultaneously even. In particular, for \(0<r<1\),
\begin{equation}\label{eq:causal-phase-witness}
                  \operatorname{Im}\chi_{\mathrm c}(1,-r)>0.
\end{equation}
At the resonant rays its continuous representative satisfies
\(\chi_{\mathrm c}(\xi,-\xi)=0\), \(\xi\ne0\).
\end{proposition}
\begin{proof}
The physical mixed kernel of the Volterra coordinate is
\[
 K_{s,t}(u,v)=\ind_{\{u<v\}}
       (\kappa_t-\kappa_s)(u)\kappa_t(v)
  -\ind_{\{v<u\}}\kappa_s(u)(\kappa_t-\kappa_s)(v).
\]
Its \(L^2\) integrability follows from
Lemma~\ref{lem:valley-tails}, since \(2H<1\).
Set
\[
 \begin{aligned}
 T_t(u,v)&=\ind_{\{u<v\}}\kappa_t(u)\kappa_t(v),\\
 K_{s,t}&=T_t-T_s-\kappa_s\otimes(\kappa_t-\kappa_s).
 \end{aligned}
\]
The formal canonical kernel has time derivative
\((\kappa_t-\kappa_s)\otimes\partial_t\kappa_t\).
Subtracting it leaves the ordinary locally integrable density
\begin{equation}\label{eq:causal-defect-density}
 b_t(u,v)=\ind_{\{u<v\}}\kappa'_t(u)\kappa_t(v)
             -\ind_{\{v<u\}}\kappa_t(u)\kappa'_t(v).
\end{equation}
Only the derivative at the smaller auxiliary time occurs in
each term. Near \(u=v=t\), its absolute integral over a ball
of radius \(\epsilon\) is \(O(\epsilon^{2H})\).

Here is a distributional justification of this subtraction.
Replace \(x_+^{-\delta}\) by
\(f_\epsilon(x)=(x+\epsilon)^{-\delta}\ind_{\{x>0\}}\).
In \(\partial_tT_t^\epsilon-
\kappa_t^\epsilon\otimes\partial_t\kappa_t^\epsilon\),
the boundary term at \(v=t\) cancels exactly. The boundary
at \(u=t<v\) vanishes by causal support. The remaining two
wedge densities are bounded in absolute value by those in
\eqref{eq:causal-defect-density}. Dominated convergence against
Schwartz tests therefore gives \(b_t\), with no corner term.
For the canonical part, the negative-exponent transforms of
\(f_\epsilon\) are bounded by \(C|\xi|^{-\gamma}\), uniformly
in \(\epsilon\): split at \(|\xi|^{-1}\), then integrate the
decreasing tail by parts. With the unitary normalization,
\[
 \frac1{\sqrt{2\pi}}\int_0^\infty
       e^{-ix\xi}f_\epsilon(x)\dd x
 \longrightarrow
       \frac{\Gamma(\gamma)}{\sqrt{2\pi}}[i\xi+0]^{-\gamma}.
\]
The negative exponent comes from \(u=t-x\) in the
positive-exponent innovation transform.
This proves convergence to the displayed \(C_{s,t}\) in
tempered distributions. Thus the complete defect is
\(\int_s^t b_r\dd r\).

To transform \(b_0\), set \(u=-x,v=-rx\) on \(u<v<0\).
The radial power, including the Jacobian, is \(x^{2H-1}\).
Abel damping in \(\max(-u,-v)\), followed by the Gamma
integral in \(x\), gives
\eqref{eq:causal-Fourier-density}; the other wedge exchanges
\(\xi,\eta\). Passage to the Abel limit is valid locally in
\(L^1\), since \(2H<1\) and \(\delta<1\).
Translation gives \(\mathcal F_+^{\otimes2}b_r
=e^{ir(\xi+\eta)}\mathfrak d_H\). Integration in \(r\)
proves \eqref{eq:causal-complete-Fourier}.

The density is homogeneous of degree \(-2H\), antisymmetric
under exchange and Hermitian under simultaneous sign reversal.
The denominator in \eqref{eq:causal-Fourier-multiplier} has
the same homogeneity and these latter two symmetries.
This gives the stated properties of \(\chi_{\mathrm c}\).
The angular factorization in Lemma~\ref{lem:Fourier-angular}
uses sign symmetry only through absolute values, so its
integrability criterion also applies to a complex Hermitian
multiplier. Since \eqref{eq:causal-complete-Fourier} is the
transform of \(K_{s,t}\in L^2\), it gives
\(J_H(\chi_{\mathrm c})<\infty\).

For the strict phase assertion, write \(\xi=1,\eta=-r\),
\(0<r<1\). The denominator is \(-ic(H)^2r^\delta\).
For \(0<q<1\),
\(|q-r|<1-qr\). Taking real parts in
\eqref{eq:causal-Fourier-density} shows
\(\operatorname{Re}\mathfrak d_H(1,-r)>0\), because the
common real phase factor is \(\cos(\pi H)>0\).
Division by the displayed denominator proves
\eqref{eq:causal-phase-witness}.
Finally the Beta integral and
\(\Gamma(\gamma)\Gamma(\delta)=\pi/\cos(\pi H)\) give
\[
 \mathfrak d_H(\xi,-\xi)
       =ic(H)^2\sgn(\xi)|\xi|^{-2H}.
\]
The sector denominator there is its negative. Continuity
at these rays follows from integrability of the endpoint
singularity of order \(2H<1\), proving the last assertion.
\end{proof}

This representation places the causal coordinate in a complex
Hermitian Fourier multiplier family, outside the real, simultaneously
even class of Proposition~\ref{prop:FNO-noncausal}.
General character-based normal ordering
accommodates broad choices of skeleton data
\citep{UnterbergerGeneric2010}; the present formula identifies
the chosen coordinate, while the higher-level and critical
estimates are supplied by the Volterra construction.

\section{The word regularity-structure model}\label{sec:word-model}
The word construction realizes a chosen enhancement as a
regularity-structure model \citep{Hairer2014,BrunedEtAl2019}.
We give the model and identify its reconstruction with the
controlled integral of Section~\ref{sec:controlled-integration}.

Fix \(N\) and \(1/(N+1)<\alpha<1/N\) as in
Section~\ref{sec:controlled-integration}. Let \(\tau_w\) denote a
formal symbol indexed by a word \(w\) of length at most \(N\)
in \(\{1,\ldots,d\}\), with \(\tau_\varnothing=\mathbf1\).
Give it degree \(|w|\alpha\), and put
\[
 A=\{0,\alpha,\ldots,N\alpha\},\qquad
 T=\bigoplus_{m=0}^N
        \operatorname{span}\{\tau_w:|w|=m\}.
\]
There is a truncated shuffle product:
\(\tau_u\tau_v=\sum_{w\in\operatorname{Sh}(u,v)}\tau_w\)
when \(|u|+|v|\leq N\), with multiplicities retained.
For a chosen geometric rough path \(\bX\), define
\begin{align}
 (\Pi_s\tau_w)(t)&=X_{s,t}^w,\label{eq:word-Pi}\\
 \Gamma_{rs}\tau_w
   &=\sum_{w=uv}X_{s,r}^{u}\tau_v.\label{eq:word-Gamma}
\end{align}
The sum includes empty words. Inverses define increments
when the endpoints are reversed.
Explicitly the structure group is
\(G=\{\mathcal G_g:g\in G^N(\R^d)\}\), where
\(\mathcal G_g\tau_w=\sum_{w=uv}g^u\tau_v\).
It obeys \(\mathcal G_h\mathcal G_g=\mathcal G_{gh}\);
this order agrees with \eqref{eq:word-Gamma}.

\begin{proposition}\label{prop:word-model}
The maps \((\Pi,\Gamma)\) define a model on \((A,T,G)\),
locally in the interior of the parameter-time domain.
At a boundary of a forward-only construction, use
one-sided tests.
They satisfy
\[
 \Pi_s=\Pi_r\Gamma_{rs},\qquad
 \Gamma_{qr}\Gamma_{rs}=\Gamma_{qs},
\]
and the analytic estimates
\[
 |\langle\Pi_s\tau_w,\varphi_s^\lambda\rangle|
       \leq C\lambda^{|w|\alpha},\qquad
 \|\Gamma_{rs}\tau_w\|_{m\alpha}
       \leq C|r-s|^{(|w|-m)\alpha},
 \quad m<|w|.
\]
Here \(\varphi_s^\lambda=\lambda^{-1}
\varphi((\cdot-s)/\lambda)\), where
\(\varphi\in C_c^\infty((-1,1))\) and
\(\|\varphi\|_{C^1}\leq1\); constants are local in time.
The model depends continuously on \(\bX\) in the corresponding
H\"older model seminorms.
\end{proposition}
\begin{proof}
The first identity is Chen's identity, with deconcatenation
\(w=uv\); applying it twice gives the group identity.
Every nonidentity term in \eqref{eq:word-Gamma} has lower
homogeneity. The estimates follow respectively from
\(|X_{s,t}^w|\leq C|t-s|^{|w|\alpha}\) and
\(|X_{s,r}^u|\leq C|r-s|^{|u|\alpha}\).
The same estimates for differences of rough paths prove
continuity. The shuffle identities ensure that \(\Pi_s\)
respects every product that is defined above.
\end{proof}

To include integration, extend the regularity structure to
\[
 \overline T=T\oplus
       \operatorname{span}\{\dot\tau_w:1\leq|w|\leq N\},
 \qquad
 \overline A=A\cup\{m\alpha-1:1\leq m\leq N\}.
\]
Set \(D\mathbf1=0\), \(D\tau_w=\dot\tau_w\), and
\[
 \Gamma_{rs}\dot\tau_w
   =\sum_{\substack{w=uv\\v\ne\varnothing}}
                      X_{s,r}^{u}\dot\tau_v.
\]
This gives the group action on the negative sector and
commutes with \(D\). Distributional differentiation defines
\[
 \Pi_s\dot\tau_w=\partial_tX_{s,t}^{w},\qquad
 |\dot\tau_w|=|w|\alpha-1.
\]
In particular \(\dot\tau_i\) represents
\(\Xi_i=\partial_tB^i\), and
\(\dot\tau_{ij}\) represents the chosen meaning of
\((B^i-B_s^i)\partial_tB^j\).
It is defined by differentiating the enhanced coordinate;
it is not obtained by multiplying two uncontrolled
distributions. Integration by parts against a test function
gives
\begin{equation}\label{eq:derivative-model-bound}
 |\langle\Pi_s\dot\tau_w,\varphi_s^\lambda\rangle|
       \leq C\lambda^{|w|\alpha-1}.
\end{equation}
Differentiating the identity for \(\Pi\) gives the compatible
reexpansion identity. Every negative homogeneity exceeds
\(-1\), so the uniform \(C^1\) test class suffices.
The negative reexpansion estimates follow from the same
increment bounds as the positive ones.

\begin{proposition}[Reconstruction of controlled integrals]
\label{prop:word-reconstruction}
Let \(c\) be an integration jet over \(\bX\), and let \(J\) be its
integral from Theorem~\ref{thm:controlled-integration}.
Set \(\rho=(N+1)\alpha\), \(\gamma=\rho-1>0\), and
\[
 f(s)=\sum_{1\leq|w|\leq N}c_w(s)\dot\tau_w.
\]
For the extended word model, \(f\in\mathcal D^\gamma\) and
\(\mathcal Rf=\partial_tJ\). At a boundary of a forward-only
construction, the reconstruction estimate uses one-sided tests.
\end{proposition}
\begin{proof}
The \(\dot\tau_v\) coefficient of
\(f(t)-\Gamma_{ts}f(s)\) is precisely the expression in
\eqref{eq:integration-jet}. Its exponent is
\((N+1-|v|)\alpha=\gamma-(|v|\alpha-1)\), proving the
modelled-distribution estimate.
For uniformly \(C^1\)-bounded tests, integration by parts
and \eqref{eq:word-integral-jet} imply
\[
 |\langle\partial_tJ-\Pi_sf(s),\varphi_s^\lambda\rangle|
 \leq C\lambda^{-2}\int_{|t-s|\leq\lambda}|t-s|^\rho\dd t
 \leq C\lambda^\gamma .
\]
Uniqueness of reconstruction for \(\gamma>0\) identifies
\(\mathcal Rf=\partial_tJ\).
\end{proof}

For the models of Sections~\ref{sec:causal-construction},
\ref{sec:critical-construction}, \ref{sec:moving-lag} and
\ref{sec:zero-lag}, continuity of the word-model construction transfers
enhanced-cutoff convergence to convergence in the model seminorms.
At a reciprocal index, the fixed-reference central subtraction
changes a highest-word coordinate by \(-(t-s)Z_R^w\);
its distributional derivative changes by the time-constant
random distribution \(-Z_R^w\).
In the lagged model the corresponding derivative correction is
the continuous random density \(-z_R^\ell(t)^w\), before
taking the renormalized limit.
The joint zero-lag limit also converges in these model seminorms.
Its integrated central correction need not have a continuous
density; its distributional derivative is obtained by continuity
in the negative model sector.
The negative-sector reexpansion maps are unchanged by
this central top-degree calibration, since all their
nonconstant coefficients have degree below \(N\).

The calibrations are random and history-dependent, rather than
deterministic BPHZ counterterms for a local Gaussian model.
Random kernel renormalizations also
occur in \citet{Unterberger2013}.
Section~\ref{sec:kernel-bphz} compares the resulting two-vertex
kernel with the causal selection.
General stochastic BPHZ estimates, such as those obtained by the
spectral-gap approach of \citet{HairerSteele2024}, have separate
analytic and probabilistic hypotheses. Here the model bounds
follow from the enhancement estimates.

At each Hurst index considered here, choose
\(N=\lfloor1/H\rfloor\) and \(1/(N+1)<\alpha<H\).
The condition \((N+1)\alpha>1\) makes sewing applicable at this
finite truncation. Sewing and reconstruction then yield the same
controlled integral for the chosen enhancement.

\bibliographystyle{spr-ims-nameyear}
\bibliography{main}
\end{document}